Zuse Institute Berlin

Takustr. 7
14195 Berlin
Germany

Timo Berthold  Peter J. Stuckey  Jakob Witzig[1]

# Local Rapid Learning for Integer Programs

[1] 0000-0003-2698-0767

# Local Rapid Learning for Integer Programs

Timo Berthold[1], Peter J. Stuckey[2], and Jakob Witzig[3]

[1]Fair Isaac Germany GmbH, Takustr. 7, 14195 Berlin, Germany
timoberthold@fico.com

[2]Monash University, Melbourne, Australia
Peter.Stuckey@monash.edu

[3]Zuse Institute Berlin, Takustr. 7, 14195 Berlin, Germany
witzig@zib.de

**Abstract**

Conflict learning algorithms are an important component of modern MIP and CP solvers. But strong conflict information is typically gained by depth-first search. While this is the natural mode for CP solving, it is not for MIP solving. Rapid Learning is a hybrid CP/MIP approach where CP search is applied at the root to learn information to support the remaining MIP solve. This has been demonstrated to be beneficial for binary programs. In this paper, we extend the idea of Rapid Learning to integer programs, where not all variables are restricted to the domain $\{0, 1\}$, and rather than just running a rapid CP search at the root, we will apply it repeatedly at local search nodes within the MIP search tree. To do so efficiently, we present six heuristic criteria to predict the chance for local Rapid Learning to be successful. Our computational experiments indicate that our extended Rapid Learning algorithm significantly speeds up MIP search and is particularly beneficial on highly dual degenerate problems.

## 1 Introduction

Constraint programming (CP) and integer programming (IP) are two complementary ways of tackling discrete optimization problems. Hybrid combinations of the two approaches have been used for many years, see, e.g., [2, 9, 10, 17, 37, 42, 22]. Both technologies have incorporated *conflict learning* capabilities [21, 27, 38, 1, 35] that derive additional valid constraints from the analysis of infeasible subproblems extending methods developed by the SAT community [33].

Conflict learning is a technique that analyzes infeasible subproblems encountered during a tree search algorithm. In a tree search, each subproblem can be identified by its local variable bounds, i.e., by local bound changes that come from branching decisions and propagation at the current node and its ancestors. If propagation detects infeasibility, conflict learning will traverse this chain of decisions and deductions reversely, reconstructing which bound changes led to which other bound changes. In this way, conflict learning identifies explanations for the infeasibility. If it can be shown that a small subset of the bound changes suffices to prove infeasibility, a so-called conflict constraint is generated that can be exploited in the remainder of the search to prune parts of the tree.

In the context of constraint programming, conflict constraints are also referred to as *no-goods*. For binary programs (BPs), i.e., mixed integer (linear) programs for which all variables have domain $\{0, 1\}$, conflict constraints will have the form of *set covering* constraints. These are linear constraints of the form "sum of variables (or their negated form) is greater than or equal to one".

Rapid Learning [13] is a heuristic algorithm for BPs that searches for valid conflict constraints, global bound reductions, and primal solutions. It is based on the observation that a CP solver can typically perform an incomplete search on a few thousand nodes in a fraction of the time that a MIP solver needs for processing the root node. In addition, CP solvers make use of depth-first search, as opposed to the hybrid best-first/depth-first search of MIP solvers, which more rapidly generates strong no-goods. Typically CP solvers do not differentiate the root node from other nodes. They apply fast (at least typically) propagation algorithms to infer new information about the possible values variables can take, and then take branching decisions. In contrast, a MIP solver invests a substantial amount of time at the root node to gather global information about the problem and to initialize statistics that can help for the search. A significant portion of root node processing time comes from the computational effort needed to solve the initial LP relaxation from scratch. Further aspects are the LP resolves during cutting plane generation, strong branching [7] for branching statistic evaluation, and primal heuristics, see, e.g., [11].

The idea of Rapid Learning is to apply a fast CP depth-first branch-and-bound search for a few hundred or thousand nodes, generating and collecting valid conflict constraints at the root node of a MIP search. Using this, the MIP solver is already equipped with the valuable information of which bound changes will lead to an infeasibility, and can avoid them by propagating the derived constraints. Just as important, the partial CP search might find primal solutions, thereby acting as a primal heuristic. Furthermore, the knowledge of conflict constraints can be used to initialize branching statistics, just like strong branching. In this paper, we will extend Rapid Learning to integer programs and to nodes beyond the root.

The remainder of the paper is organized as follows. In Section 2, we provide more background on conflict learning for MIPs, in particular the extension to general integer variables, which is important for our extended Rapid Learning algorithm. In Section 3, we describe details of the Rapid Learning algorithm for general integer programs, extending the work of Berthold et al. [13]. In Section 4, we discuss what special considerations have to be taken when applying Rapid Learning repeatedly at local subproblems during the MIP tree search instead of using it as a onetime global procedure. We introduce six criteria to predict the benefit of local Rapid Learning. Section 5 presents our computational study, in which we apply our extended Rapid Learning algorithm to a set of integer programs from the well-known benchmark sets of Miplib 3, Miplib 2003, and Miplib 2010 [28]. The experiments have been conducted with the constraint integer programming solver SCIP [24] and indicate that a significant speed-up can be achieved for (pure) integer programs, when using Rapid Learning locally. In Section 6, we conclude.

## 2 Conflict Learning in Integer Programming

A mixed integer program is a mathematical optimization problem defined as follows.

**Definition 1 (mixed integer program)** *Let $m, n \in \mathbb{Z}_{\geq 0}$. Given a matrix $A \in \mathbb{R}^{m \times n}$, a right-hand-side vector $b \in \mathbb{R}^m$, an objective function vector $c \in \mathbb{R}^n$, a lower and an upper bound vector $l \in (\mathbb{R} \cup \{-\infty\})^n$, $u \in (\mathbb{R} \cup \{+\infty\})^n$ and a subset $\mathcal{I} \subseteq \mathcal{N} = \{1, \ldots, n\}$, the corresponding* mixed integer program (MIP) *is given by*

$$\begin{aligned} \min \quad & c^{\mathsf{T}} x \\ s.t. \quad & Ax \leq b \\ & l_j \leq x_j \leq u_j \quad \text{for all } j \in \mathcal{N} \\ & x_j \in \mathbb{R} \quad \text{for all } j \in \mathcal{N} \setminus \mathcal{I} \\ & x_j \in \mathbb{Z} \quad \text{for all } j \in \mathcal{I}. \end{aligned} \tag{1}$$

Mixed integer programs can be categorized by the classes of variables that are part of their formulation:

- If $\mathcal{N} = \mathcal{I}$, problem (1) is called a *(pure) integer program (IP)*.
- If $\mathcal{N} = \mathcal{I}$, $l_j = 0, j \in \mathcal{N}$ and $u_j = 1, j \in \mathcal{N}$, problem (1) is called a *(pure) binary program (BP)*.
- If $\mathcal{I} = \emptyset$, problem (1) is called a *linear program (LP)*.

Conflict analysis techniques were originally developed by the artificial intelligence research community [40] and, later extended by the SAT community [33]; they led to a huge increase in the size of problems modern SAT solvers can handle [31, 33, 43]. The most successful SAT learning approaches use so-called *one-level first unique implication point (1-UIP)* [43] learning which in some sense captures the conflict constraint "closest" to the infeasibility. Conflict analysis also is successfully used in the CP community [25, 26, 35] (who typically refer to it as no-good learning) and the MIP world [1, 21, 38, 41]. Nowadays, commercial MIP solvers like FICO Xpress [23] employ conflict learning by default.

Constraint programming and mixed integer programming are two complementary ways of tackling discrete optimization problems. Because they have different strengths and weaknesses hybrid combinations are attractive. One notable example, the software SCIP [3], is based on the idea of *constraint integer programming* (CIP) [2, 6]. CIP is a generalization of MIP that supports the notion of general constraints as in CP. SCIP itself follows the idea of a very low-level integration of CP, SAT, and MIP techniques. All involved algorithms operate on a single search tree and share information and statistics through global storage of, e.g., solutions, variable domains, cuts, conflicts, the LP relaxation and so on. This allows for a very close interaction amongst CP and MIP (and other) techniques.

There is one major difference between BPs and IPs in the context of Rapid Learning: in IP, the problem variables are not necessarily binary. To deal with this, the concept of a *conflict graph* needs to be extended. A conflict graph gets constructed whenever infeasibility is detected in a local search node; it represents the logic of how the set of branching decisions led to the detection of infeasibility.

More precisely, the conflict graph is a directed acyclic graph in which the vertices[1] represent bound changes of variables, e.g., $x_i \leq \lambda_i$ or $x_i \geq \mu_i$. The conflict graph is built such that when the solver infers a bound change $v$ as a consequence of a set of existing bound changes $U$, i.e., $U \rightarrow v$, then we have an arc $(u, v)$ from each $u \in U$ to $v$. Bound changes caused by branching decisions are vertices without incoming edges. Finally the conflict graph includes a dummy vertex *false* representing failure which is added when the solver infers unsatisfiability.

Given a conflict graph, each cut that separates the branching decisions from the artificial infeasibility vertex *false* gives rise to a valid conflict constraint. A *unique implication point (UIP)* is an (inner) vertex of the conflict graph which is traversed by all paths from the branching vertices to the conflict vertex. Or, how Zhang et al. [43] describe it: "Intuitively, a UIP is the *single* reason that implies the conflict at [the] current decision level." UIPs are natural candidates for finding small cuts in the conflict graph. The *1-UIP* is the first cut separating the conflict vertex from the branching decisions when traversing in reverse assignment order.

For integer programs, conflict constraints can be expressed as so-called *bound disjunction* constraints:

**Definition 2** *For an IP, let $\mathcal{L} \subseteq \mathcal{I}, \mathcal{U} \subseteq \mathcal{I}$ be disjoint index sets of variables, let $\lambda \in \mathbb{Z}^{\mathcal{L}}$ with $l_i \leq \lambda_i \leq u_i$ for all $i \in \mathcal{L}$, and $\mu \in \mathbb{Z}^{\mathcal{U}}$ with $l_i \leq \mu_i \leq u_i$ for all $i \in \mathcal{U}$. Then, a constraint of the form*

$$\bigvee_{i \in \mathcal{L}} (x_i \geq \lambda_i) \vee \bigvee_{i \in \mathcal{U}} (x_i \leq \mu_i)$$

*is called a* bound disjunction *constraint.*

[1] For disambiguation, we will use the term *vertex* for elements of the conflict graph, as opposed to *nodes* of the search tree.

For details on bound disjunction constraints, see Achterberg [1]. If all involved conflict values $\lambda, \mu$ correspond to global bounds of the variables, the bound disjunction constraint can be equivalently expressed as a knapsack constraint of form

$$\sum_{i \in \mathcal{U}} x_i - \sum_{i \in \mathcal{L}} x_i \leq \sum_{i \in \mathcal{U}} u_i - \sum_{i \in \mathcal{L}} l_i - 1. \tag{2}$$

Note that for BPs all conflicts only involve global bounds.

The power of conflict learning arises because often branch-and-bound based algorithms implicitly repeat the same search in a slightly different context in another part of the tree. Conflict constraints help to avoid redundant work in such situations. As a consequence, the more search is performed by a solver and the earlier conflicts are detected, the greater the chance for conflict learning to be beneficial. Note that conflict generation has a positive interaction with depth-first search. Depth-first search leads to the creation of no-goods that explain why a whole subtree contains no solutions, and hence the no-goods generated by depth-first search are likely to prune more of the subsequent search.

# 3 Rapid Learning for Integer Programs

The principle motivation for Rapid Learning [13] is the fact that a CP solver can typically search hundreds or thousand of nodes in a fraction of the time that a MIP solver needs for processing the root node of the search tree. Rapid Learning applies a fast CP search[2] for a few hundred or thousand nodes, before starting the MIP search. Using this approach, conflict constraints can be learnt before, and not only during, MIP search. Very loosely speaking: while the aim of conflict learning is to avoid making mistakes a second time, Rapid Learning tries to avoid making them the first time (during MIP search).

Rapid Learning is related to large neighborhood search heuristics, such as RINS and RENS [12, 20]. But, rather than doing an incomplete search on a subproblem using the same (MIP search) algorithm, Rapid Learning performs an incomplete search on the same problem using a much faster algorithm (CP search). Rapid Learning differs from primal heuristics in that it aims at improving the dual bound by collecting information on infeasibility rather than searching for feasible solutions.

Each piece of information collected in a rapid CP search can be used to guide the MIP search or even deduce further reductions during root node processing. Since the CP solver is solving the same problem as the MIP solver

- each generated conflict constraint is valid for the MIP search,
- each global bound change can be applied at the MIP root node,
- each feasible solution can be added to the MIP solver's solution pool,

[2]By CP search we mean applying a depth-first search using only propagation for reasoning, no LP relaxation is solved during the search.

- the branching statistics can initialize a hybrid MIP branching rule, see [4], and
- if the CP solver completely solves the problem, the MIP solver can abort.

All five types of information may be beneficial for a MIP solver, and are potentially generated by our algorithm which we now describe more formally.

The Rapid Learning algorithm is outlined in Figure 1. Here, $l(P)$ and $u(P)$ are lower and upper bound vectors, respectively, of the problem at hand, $P$. For the moment we assume $P$ is the root problem, in the next section we will examine the use of Rapid Learning at subproblem nodes. The symbol $C$ refers to a single globally valid conflict constraint explaining the infeasibility of the current subproblem. Rapid Learning is an incomplete CP search: a branch-and-bound algorithm which traverses the search space in a depth-first manner (Line 3), using propagation (Line 4) and conflict analysis (Line 7), but no LP relaxation. Instead, the *pseudo-solution* [2], i.e., an optimal solution of a relaxation consisting only of the variable bounds (Line 5), is used for the bounding step.

Propagation of linear constraints is conducted by the bound strengthening technique of Brearley et al. [18] which uses the residual activity of linear constraints within the local bounds. For special cases of linear constraints, SCIP implements special, more efficient propagators. Knapsack constraints use efficient integer arithmetic instead of floating point arithmetic, and sort by coefficient values to propagate each variable only once. SCIP also features methods to extract clique information about the binary variables of a problem. A clique is a set of binary variables of which at most one variable can take the value 1 in a feasible solution. Clique information can be used to strengthen the propagation of knapsack constraints. Set cover constraints are propagated by the highly efficient two-watched literal scheme [33], which is based on the fact that the only domain reduction to be inferred from a set cover constraint is to fix a variable to 1 if all other variables have already been fixed to 0.

Variable and value selection takes place in Line 14; inference branching [2] is used as branching rule. Inference branching maintains statistics about how often the fixing of a variable led to fixings of other variables, i.e., it is a history rule, its essentially a MIP equivalent of impact-based search [29, 36]. Since history rules are often weak in the beginning of the search, we seed the CP solver with statistics that the MIP solver has collected in probing [39] during MIP presolving.

We assume that the propagation routines in Line 4 may also deduce global bound changes and modify the global bound vectors $l(P)$ and $u(P)$. Single-clause conflicts are automatically upgraded to global bound changes in Line 9. Note that it suffices to check constraint feasibility in Line 11, since the pseudo-solution $\bar{x}$ (see Line 5) will always take the value of one of the (integral) bounds for each variable.

Our implementation of the Rapid Learning heuristic uses a secondary SCIP instance to perform the CP search. Only a few parameters need to be altered from their default values to turn SCIP into a CP solver, an overview is given in

Figure 1: Rapid Learning algorithm

---

**Input** : IP $P$ as in (1) (with $\mathcal{R} = \emptyset$),
node limit $\lim_{\text{node}}$,
primal bound $\bar{c}$ for $P$ (might be $\infty$)

**Output:** set of valid conflict constraints $\mathcal{L}_C$ for $P$,
valid global domain box $[l, u]$ for $P$,
feasible solution $\tilde{x}$ for $P$ or $\emptyset$

$\mathcal{L} \leftarrow \{P\}$, $\text{n}_{\text{node}} \leftarrow 0$, $\mathcal{L}_C \leftarrow \emptyset$, $\tilde{x} \leftarrow \emptyset$;
**while** $\mathcal{L} \neq \emptyset \wedge \text{n}_{\text{node}} < \lim_{\text{node}}$ **do**
    $\tilde{P} \leftarrow$ `select_dfs(`$\mathcal{L}$`)`, $\mathcal{L} \leftarrow \mathcal{L} \setminus \tilde{P}$, $\text{n}_{\text{node}} \leftarrow \text{n}_{\text{node}} + 1$;
    $[l(\tilde{P}), u(\tilde{P})] \leftarrow$ `propagate(`$[l(\tilde{P}), u(\tilde{P})]$`)`;
    $\bar{x} \leftarrow \text{argmin}\{c^{\mathsf{T}} x \mid x \in [l(\tilde{P}), u(\tilde{P})]\}$;

    `/* analyze infeasible subproblem, potentially store globally valid conflict constraint */`
    **if** $[l(\tilde{P}), u(\tilde{P})] = \emptyset$ **or** $c(\bar{x}) \geq \bar{c}$ **then**
        $C \leftarrow$ `analyze(`$\tilde{P}$`)`;
        **if** $C \neq \emptyset$ **then** $\mathcal{L}_C \leftarrow \mathcal{L}_C \cup \{C\}$;
        **if** $|C| = 1$ **then** `tighten(`$[l(P), u(P)]$`)`;
        **continue**;

    `/* check for new incumbent solution */`
    **if** $A\bar{x} \leq b$ **and** $c^T\bar{x} < \bar{c}$ **then**
        $\tilde{x} \leftarrow \bar{x}$, $\bar{c} \leftarrow c^T\bar{x}$;
        **continue**;

    $(x_i, v) \leftarrow$ `select_infer(`$\tilde{P}$, $\bar{x}$`)`;
    $\tilde{P}_l \leftarrow \tilde{P} \cup \{x_i \leq v\}$, $\tilde{P}_r \leftarrow \tilde{P} \cup \{x_i \geq v\}$;
    $\mathcal{L} \leftarrow \mathcal{L} \cup \{\tilde{P}_l, \tilde{P}_r\}$;
**return** $(\mathcal{L}_C, [l(P), u(P)], \tilde{x})$;

---

Table 1. Most importantly, we disabled the LP relaxation and use a pure depth-first search with inference branching (but without any additional tie breakers). Further, we switch from All-UIP to 1-UIP in order to generate only one conflict per infeasibility. This is a typical behavior of CP solvers, but not for MIP solvers. Expensive feasibility checks and propagation of the objective function as a constraint are also avoided.

In order to avoid spending too much time in Rapid Learning, the number of nodes explored during the CP search is limited to at most 5000. The actual number of allowed nodes is determined by the number of simplex iterations

Table 1: Settings for Rapid Learning sub-SCIP.

| parameter name | value | effect |
|---|---|---|
| lp/solvefreq | -1 | disable LP |
| conflict/fuiplevels | 1 | use 1-UIP |
| nodeselection/dfs/stdpriority | INT_MAX/4 | use DFS |
| branching/inference/useweightedsum | FALSE | pure inference, no VSIDS |
| constraints/disableenfops | TRUE | no extra checks |
| propagating/pseudoobj/freq | -1 | no objective propagation |
| conflict/maxvarsfac | 0.05 | only short conflicts |
| history/valuebased | TRUE | extensive branch. statistics |

$\text{iter}_{\text{LP}}$ performed so far in the main SCIP but at least 500, i.e.,

$$\lim_{\text{node}} = \min\{5000, \max\{500, \text{iter}_{\text{LP}}\}\}.$$

The idea is to restrict Rapid Learning more rigorously for problems where processing of a single MIP node is cheap already. The number of simplex iterations is a deterministic estimate for node processing cost.

We aim to generate short conflict constraints, since these are most likely to frequently trigger propagations in the upcoming MIP search. Thus, we only collect conflicts that contain at most 5 % of the problem variables. Finally, we adapt the collection of branching statistics such that history information on general integer variables are collected per value in the domain rather than having one counter for down- and one for up-branches regardless of the value on which was branched. This can be essential for performing an efficient CP search on general integer variables, and was a building block that enabled us to use Rapid Learning on IPs rather than solely on BPs, as in [13].

In addition to the particular parameters listed in Table 1, we set the emphasis[3] for presolving to "fast". Emphasis settings for cutting are not necessary, since no LP relaxation is solved, from the armada of primal heuristics only a few are applied that do not require an LP relaxation, see [5]. Note that since Rapid Learning will be called at the end of the MIP root node, or even locally, see next Section, the problem that the CP solver considers has already been presolved, might contain cutting planes as additional linear constraints and have an objective cutoff constraint if a primal solution has been found by a primal heuristic during root node processing.

# 4 Local Rapid Learning

The original Rapid Learning algorithm [13] was used as part of a root preprocessing, i.e., for every instance it was run exactly once at the end of the root node. But only running Rapid Learning at the root limits its effectiveness. We

[3]In SCIP, emphasis settings correspond to a group of individual parameters being changed.

now discuss the factors that arise when we allow Rapid Learning to be run at local nodes inside the search tree

When running in the root only all information returned by the CP solver is globally valid, and the overhead to maintain the information gathered by Rapid Learning is negligible [13]. In contrast, when applying Rapid Learning at a local node within the tree conflicts and bound changes will only be locally valid in general. Since Rapid Learning uses a secondary SCIP instance to perform the CP search, all local information of the current node becomes part of the initial problem formulation for the CP search. Thus, conflicts gathered by Rapid Learning do not include bound changes made along the path from the root to the current node, they are simply considered as valid for this local node. As a consequence, these conflicts will only be locally valid and hence only applied to the current node of the MIP search. Using an assumption interface [34], local conflicts could be lifted to be globally valid. However, this is subject to future investigation and not considered in the current implementation of Rapid Learning.

In practice, all local information needs to be maintained when switching from one node of the tree to another. In CP solvers, switching nodes is typically very cheap, because depth-first search is used. However, a MIP solver frequently "jumps" within the tree. Therefore, two consecutively processed nodes can be quite different. In what follows, we will refer to the time spent for moving from one node to another node as *switching time*. The switching time can be used as an indicator to quantify the overhead introduced by all locally added information found by Rapid Learning.

To ensure that the amount of locally added information does not increase the switching time too much, we apply Rapid Learning very rarely by using a exponentially decreasing frequency of execution. Rapid Learning is executed at every node of depth $d$ with

$$\log_{\beta}(d/f) \in \mathbb{Z}, \tag{3}$$

where $\beta$ and $f$ are two parameters to control the speed of decrease. For example, if $\beta = 1$ Rapid Learning is executed at every depth $d = i \cdot f$ with $i \in \mathbb{Z}_{+}$.

Unfortunately, the amount of locally valid information produced by Rapid Learning still leads to an increase of switching time by 21 %. Consequently, the overall performance decreased by 20 % in our first experiments. At the same time the number of explored branch-and-bound nodes decreased by 16 %. This indicates the potential gains possible using local Rapid Learning.

To control at which subproblem Rapid Learning is applied we propose six criteria to forecast the potential of Rapid Learning. These criteria aim at identifying one of two situations. The first is to estimate whether the (sub)problem is infeasible or a pure feasibility problem. In these cases propagating conflicts is expected to be particularly beneficial. The second is to estimate the dual degeneracy of a problem. In this case, VSIDS branching statistics are expected to be particularly beneficial. The VSIDS [31] (variable state independent decaying sum) statistics takes the contribution of every variable (and its negated complement) in conflict constraints found so far into account. For every variable, the number of clauses (in MIP speaking: conflict constraints) the variable is part of is counted. In the remainder of the search the VSIDS are periodically scaled by a predefined constant. By this, the weight of older clauses is reduced over time and more recent observations have a bigger impact.

A basic solution of an LP is called *dual degenerate*, when it has nonbasic variables with zero reduced costs. One can define the dual degeneracy of a MIP as the average number of nonbasic variables with zero reduced costs appearing in a basic solution of its LP relaxation. The higher the dual degeneracy, the higher the chance that the LP objective will not change by branching and hence many of the costs involved in the pseudo-cost computation are zero. Therefore, for highly dual degenerate problems, using other branching criteria, such as VSIDS or inference scores, is crucial for solving the problem.

We now describe the six criteria we use to identify infeasible or dual degenerate problems, already using the criteria abbreviations from the tables in Section 5:

**Criterion I: Dual Bound Improvement.** During the tree search a valid lower bound for each individual subproblem is given by the respective LP solution. A globally valid lower bound is given by the minimum over all individual lower bounds. This global bound is called the *dual bound*. If the dual bound has not changed after processing a certain number of nodes, i.e., the dual bound is equal to the lower bound of the root node, it might be the case that the MIP lies inside a level plane of the objective, i.e., all feasible LP (and MIP) solutions will have the same objective. In other words, the instance might be a feasibility instance for which Rapid Learning was already shown to be very successful [13]. Feasibility instances are typically highly dual degenerate. The `dualbound` criterion means to call local Rapid Learning if the dual bound never changed during the MIP search.

**Criterion II: Leaves Pruned by Infeasibility or Exceeding the Cutoff bound.** During the tree search every leaf node either provides a new incumbent solution (the rare case), is proven to be infeasible or to exceed the current cutoff bound which is given by the incumbent solution. The ratio of the latter two cases is used in SCIP's default branching rule. *Hybrid branching* [4] combines pseudo-costs, inference scores, and conflict information into one single branching score. The current implementation in SCIP puts a higher weight on conflict information, e.g., VSIDS [31], and a lower weight on pseudo-costs when the ratio of infeasible and cutoff nodes is larger than a predefined threshold. The `leaves` criterion means to call local Rapid Learning if the ratio of infeasible leaves over those exceeding the cutoff bound is larger than 10. The rationale is that we expect (local) conflicts to be most beneficial, when infeasibility detection appears to be the main driver for pruning the tree.

**Criterion III: LP Degeneracy.** As mentioned above, the more nonbasic variables are dual degenerate, the less information can be gained during strong branching or pseudo-cost computation. As a consequence, Berthold et al. [14] introduced a modification to strong branching that considers the dual degeneracy of the LP solution. In rough terms, if either the share of dual degenerate nonbasic variables or the variable-constraint ratio of the optimal face exceed certain thresholds, strong branching will be deactivated. We adapt this idea of using the dual degeneracy of the current LP solution. The `degeneracy` criterion means to call local Rapid Learning if more than 80 % of the nonbasic variables are degenerate or the variable-constraint ratio of the optimal face is larger than 2, as proposed in [14]. In both cases we expect that "strong conflict generation" will be useful.

**Criterion IV: (Local) Objective Function.** If all variables with non-zero objective coefficients are fixed at the local subproblem, i.e., the objective is constant, Criteria I and II will apply: every LP solution is fully dual degenerate and the only possibility to prune a leaf node is by infeasibility. If there are only very few unfixed variables with nonzero objective are left, the criteria might not apply. However, it is likely that the targeted situations occur frequently in the tree rooted at the current subproblem, at the latest, when all the variables occurring in the objective are fixed. The `obj` criterion means to call local Rapid Learning once the objective support is small enough, in anticipation of the current subproblem turning into a feasibility problem. In our implementation we apply this criterion very conservatively, and call Rapid Learning only if the local objective is zero.

**Criterion V: Number of Solutions.** The most obvious evidence, and indeed a necessary one, that a MIP instance is infeasible, is that no feasible solution has been found during the course of the MIP search. Note that for most (feasible) MIP instances, primal heuristics find a feasible solution at the root node [11] or at the latest during the first dive in the branch-and-bound. The `nsols` criterion means to call local Rapid Learning if no feasible solution has been found so far.

**Criterion VI: Strong Branching Improvements.** In the beginning of the tree search it is very unlikely that enough leaf nodes are explored to reliably guess whether the actual MIP is a feasibility instance. Therefore, we consider the subproblems evaluated during strong branching, which are concentrated at the top of the search tree. Similarly to Criterion II, we compute the ratio between the number of strong branching problems that gave no improvement in the objective or went infeasible to the number of strong branching problems where we observed an objective change. The `sblps` criterion means to call local Rapid Learning if this ratio exceeds a threshold of 10, hence strong branching does not appear to be efficient for generating pseudo-cost information.

In addition to the exponentially decreasing frequency and the six criteria above, we applied the following three changes to the original implementation of Rapid Learning used in [13].

- We limited the number of conflict constraints transferred from Rapid Learning back to the original search tree to ten. This corresponds to the SCIP parameter `conflict/maxconss` for the maximal allowed number of added conflicts per call of conflict analysis. We greedily use the shortest conflicts.

- We prefer conflict constraints that have a linear representation over bound disjunction constraints (see Definition 2).

- To exploit performance variability [19, 30] every CP search is initialized with a different pseudo-random seed.

# 5 Computational results

To evaluate how local Rapid Learning impacts IP solving performance we used the academic constraint integer programming solver SCIP 6.0 [24] (with SoPlex 4.0 as LP solver) and extended the existing code of Rapid Learning. The original implementation of Rapid Learning was already shown to significantly improve the performance of SCIP 1.2.0.5 on pure binary instances [13]. In this setting, Rapid Learning was applied exactly once at the root node. However, during the last eight years SCIP has changed in many places. In SCIP 6.0, Rapid Learning is deactivated by default, since it led to a big performance variability.

Therefore, we use SCIP without Rapid Learning (as it is the current default) as a baseline. We will refer to this setting as `default`. In our computational experiments we evaluate the impact of local Rapid Learning if one or more of the criteria described in Section 4 are fulfilled. In the following, we will refer to the criteria I–VI as `dualbound`, `leaves`, `degeneracy`, `obj`, `nsols`, and `sblps`, respectively. Within the tree, Rapid Learning is applied with an exponentially decreasing frequency (see Section 4). In our experiments, we used $f = 5$ and $\beta = 4$, i.e., Rapid Learning is called at depths $d$ with $log_4(d/5) \in \mathbb{Z}$, i.e., $d = 0, 5, 20, 80, 320\ldots$, if one of the six criteria is fulfilled.

As a test set we used all pure integer problems of Miplib 3 [16], Miplib 2003 [8] and the Miplib 2010 [28] benchmark set. This test set consists of 71 publicly available instances, which we will refer to as `MMM-IP`. To account for the effect of performance variability [19, 30] all experiments were performed with five different global random seeds. Every pair of instance and seed was treated as an individual observation, effectively resulting in a test set of 355 instances. We will use the term "instance" when actually referring to an instance-seed-combination. The experiments were run on a cluster of identical machines, each with an Intel Xeon E5-2690 with 2.6 GHz and 128 GB of RAM; a time limit of 3600 seconds was set.

In a first experiment we evaluated the efficacy of each individual criterion and global Rapid Learning as published in [13]. Aggregated results are shown in Table 2, section Exp.1. For detailed results see Table 4 in the appendix. For every setting, the table shows the number of solved instances out of 71

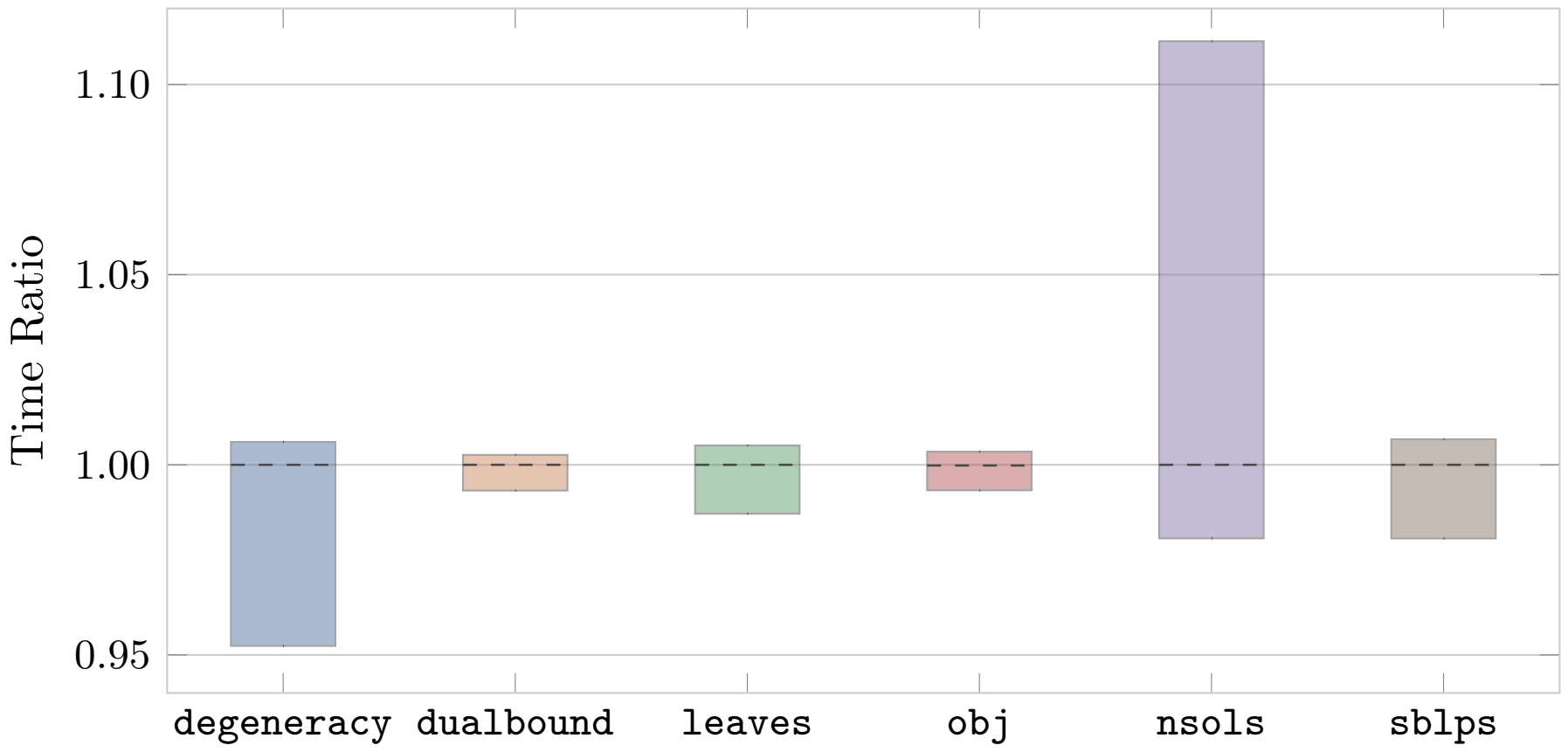

Figure 2: Box-plot of the performance ratios of the individual criteria compared to `default` on the set of affected instances.

(**solved**), shifted geometric means [3] of the absolute solving time in seconds (**time**, shift = 1) and number of explored nodes (**nodes**, shift = 100), as well as the relative solving time ($\mathbf{time_Q}$) and number of nodes ($\mathbf{nodes_Q}$) w.r.t. to `default` as a baseline. Local Rapid Learning without any of the presented criteria (`nochecks`) leads to a performance decrease of 20 % on the complete test set `MMM-IP` compared to `default`. Always applying Rapid Learning only at the root (`onlyroot`), which corresponds to Rapid Learning as published in [13], leads to slowdown of 7.8 % but solves three instance more. For this settings, we could observe a performance decrease of 38.1 % on the group of instances that are not affected[4] by Rapid Learning. To avoid a computational overhead and performance variability on instances where Rapid Learning is not expected to be beneficial, we apply the criteria `degeneracy`, `obj`, and `nsols` at the root node, too. Afterwards, the performance decrease of global Rapid Learning reduced to 0.9 %. The computational results indicate that almost all individual criteria are useful on their own. The solving time and generated nodes can be reduced by up to 6 % and 15 %, respectively, on the complete test set of 71. The exception is the `obj` criterion, which leads to a marginal slowdown of 0.6 %, but solves three more instance than `default`. On the group of affected instances the solving time and number of generated nodes can be reduced up to 8 % and 19 %, respectively, both by using `degeneracy`.

The impact of the individual criteria on the solving time is illustrated in Figure 2. For each criterion, the box plot [32] shows the median (dashed line), and the 1st and 3rd quartile (shaded box) of all observations. The plot shows that `degeneracy` performs best since it leads by far to the best improvement on the 1st quartile. In contrast to that, the 3rd quartile of `nsols` indicates this criterion leads to a deterioration of more than 10 % on 25 % of the instances.

[4]An instance is called affected when the solving path changes.

Table 2: Computational results for every individual heuristic criterion on MMM-IP.

| | | solved | time | nodes | time$_{\mathbf{Q}}$ | nodes$_{\mathbf{Q}}$ |
|---|---|---|---|---|---|---|
| Exp.1 | default | 304 | 50.92 | 2466 | – | – |
| | degeneracy | 306 | 48.02 | 2084 | 0.943 | **0.845** |
| | dualbound | 304 | 50.26 | 2369 | 0.987 | 0.961 |
| | leaves | 305 | 49.31 | 2347 | 0.968 | 0.952 |
| | obj | **307** | 51.23 | 2399 | 1.006 | 0.973 |
| | nsols | 304 | 47.94 | 2142 | **0.941** | 0.869 |
| | sblps | 304 | 50.38 | 2296 | 0.989 | 0.931 |
| | nochecks | 299 | 61.10 | 1925 | 1.200 | **0.781** |
| | onlyroot | **307** | 54.89 | 2404 | 1.078 | 0.975 |
| Exp.2 | degeneracy + leaves | 306 | 47.81 | 2073 | **0.939** | 0.841 |
| | leaves + obj | **307** | 49.80 | 2284 | 0.978 | 0.926 |
| | all6criterion | 303 | 48.23 | 2030 | 0.947 | **0.823** |

Grouping all instances of MMM-IP based on the degeneracy at the end of the root node shows the importance of this criterion. On the group of instances where at least 1 % of the variables is dual degenerate at the end of the root node Rapid Learning leads to a performance improvement of 9.1 %. On all instances where at least 80 % of the variable are dual degenerate at the root node, we could observe a reduction of solving time by 28.8 %. Note that this was one of the two thresholds for the degeneracy criterion.

In a second experiment (Table 2, section Exp.2; Table 5 in the appendix) we combined all individual criteria. Combining two or more criteria leads to more aggressive version of Rapid Learning since it runs if at least one of the chosen criteria is satisfied. The two (out of fifteen) best pairwise combinations as well as the (most aggressive) combination of all six criteria are shown in Table 2. Interestingly, no combined setting is superior to degeneracy. The combination of degeneracy and leaves, which were the two outstanding criteria in the individual test, performs almost the same as the degeneracy criterion alone.

For a final experiment we choose degeneracy as the best criterion, since it was one of two criteria that solved an additional instance, clearly showed the best search reduction, and was a close second to leaves with respect to running time. Our final experiment evaluates the impact of the individual information gained from local Rapid Learning. To this end, we individually deactivated transferring variable bounds, conflict constraints, inference information, and primal feasible solutions (see Table 3; Table 6 in the appendix). This experiment indicates that primal solutions are the most important information for the remainder of the MIP search. When ignoring solutions found during the CP search, the overall solving time increased by 9.9 % (primsols). When ignoring conflict constraints, the original motivation of Rapid Learning, solving time increased by 4.4 % (conflicts). Both transferring variable bounds and inference information proved beneficial, with a 1.4 % (variablebounds) and 0.6 % (infervals) impact on performance, respectively. It is not surprising that finding primal solutions has the largest effect. Firstly, they are applied globally, in contrast to bound changes and conflicts. Secondly, highly dual degenerate problems are known to be cumbersome not only for MIP branching but also for primal heuristics [11], which means that solution-generating procedures that do not rely on solving LPs are particularly promising for such problems.

Table 3: Performance impact of individual gained information on MMM-IP.

| | | solved | time | nodes | $\text{time}_Q$ | $\text{nodes}_Q$ |
|---|---|---|---|---|---|---|
| Exp.3 | degeneracy | 306 | 48.02 | 2084 | – | – |
| | variablebounds | 305 | 48.69 | 2089 | 1.014 | 1.002 |
| | conflicts | 306 | 50.11 | 2225 | 1.044 | 1.067 |
| | infervals | 305 | 48.29 | 2179 | 1.006 | 1.046 |
| | primsols | 305 | 52.77 | 2349 | 1.099 | 1.127 |

# 6 Conclusion

In this paper, we extended the idea of Rapid Learning [13]. Firstly, we generalized Rapid Learning to integer programs and described the details that were necessary for doing so: value-based inference branching, additional propagators and generalized conflict constraints, most of which were already available in SCIP. Secondly, we applied Rapid Learning repeatedly during the search. This generates a true hybrid CP/MIP approach, with two markedly different search strategies communicating information forth and back. To this end, we introduced six heuristic criteria to decide when to start local Rapid Learning. Those criteria are based on degeneracy information, branch-and-bound statistics, and the local structure of the problem. Our computational experiments showed a speed-up of up to 7 % when applying local Rapid Learning in SCIP. Calling local Rapid Learning depending on the local degree of dual degeneracy is the best strategy found in our experiments.

Interesting future work in this direction includes: extending the CP search to generate global conflicts at local nodes using an assumption interface, running the CP search in a parallel thread where whenever the MIP solver moves to a new node the CP search restarts from that node, and extending the method to handle problems that include continuous variables.

## Acknowledgments

The work for this article has been partly conducted within the Research Campus MODAL funded by the German Federal Ministry of Education and Research (BMBF grant number 05M14ZAM). We thank the anonymous reviewers for their valuable suggestions and helpful comments.

# References

[1] T. Achterberg. Conflict analysis in mixed integer programming. *Discrete Optimization*, 4(1):4–20, 2007.

[2] T. Achterberg. *Constraint Integer Programming.* PhD thesis, Technische Universität Berlin, 2007.

[3] T. Achterberg. SCIP: Solving Constraint Integer Programs. *Mathematical Programming Computation*, 1(1):1–41, 2009.

[4] T. Achterberg and T. Berthold. Hybrid branching. In W.-J. van Hoeve and J. N. Hooker, editors, *Integration of AI and OR Techniques in Constraint Programming for Combinatorial Optimization Problems, 6th International Conference, CPAIOR 2009*, volume 5547 of *LNCS*, pages 309–311. Springer Berlin Heidelberg, May 2009.

[5] T. Achterberg, T. Berthold, and G. Hendel. Rounding and propagation heuristics for mixed integer programming. In D. Klatte, H.-J. Lüthi, and K. Schmedders, editors, *Operations Research Proceedings 2011*, pages 71–76. Springer Berlin Heidelberg, 2012.

[6] T. Achterberg, T. Berthold, T. Koch, and K. Wolter. Constraint integer programming: A new approach to integrate CP and MIP. In L. Perron and M. A. Trick, editors, *Integration of AI and OR Techniques in Constraint Programming for Combinatorial Optimization Problems, 5th International Conference, CPAIOR 2008*, volume 5015 of *LNCS*, pages 6–20. Springer Berlin Heidelberg, May 2008.

[7] T. Achterberg, T. Koch, and A. Martin. Branching rules revisited. *Operations Research Letters*, 33(1):42–54, 2005.

[8] T. Achterberg, T. Koch, and A. Martin. Miplib 2003. *Operations Research Letters*, 34(4):361–372, 2006.

[9] E. Althaus, A. Bockmayr, M. Elf, M. Jünger, T. Kasper, and K. Mehlhorn. SCIL – symbolic constraints in integer linear programming. In *Algorithms – ESA 2002*, pages 75–87, 2002.

[10] I. D. Aron, J. N. Hooker, and T. H. Yunes. SIMPL: A system for integrating optimization techniques. In *Integration of AI and OR Techniques in Constraint Programming for Combinatorial Optimization Problems, CPAIOR 2004*, volume 3011 of *LNCS*, pages 21–36, 2004.

[11] T. Berthold. *Heuristic algorithms in global MINLP solvers.* PhD thesis, Technische Universität Berlin, 2014.

[12] T. Berthold. RENS – the optimal rounding. *Mathematical Programming Computation*, 6(1):33–54, 2014.

[13] T. Berthold, T. Feydy, and P. J. Stuckey. Rapid learning for binary programs. In A. Lodi, M. Milano, and P. Toth, editors, *Integration of AI and OR Techniques in Constraint Programming for Combinatorial Optimization Problems, 7th International Conference, CPAIOR 2010*, volume 6140 of *LNCS*, pages 51–55. Springer Berlin Heidelberg, June 2010.

[14] T. Berthold, G. Gamrath, and D. Salvagnin. Cloud Branching. In preparation.

[15] T. Berthold, P. J. Stuckey, and J. Witzig. Local Rapid Learning for Integer Programs. Technical Report 18-56, ZIB, Takustr. 7, 14195 Berlin, 2018.

[16] R. E. Bixby, S. Ceria, C. M. McZeal, and M. W. Savelsbergh. An updated mixed integer programming library: Miplib 3.0. Technical report, 1998.

[17] A. Bockmayr and T. Kasper. Branch-and-infer: A unifying framework for integer and finite domain constraint programming. *INFORMS Journal on Computing*, 10(3):287–300, 1998.

[18] A. Brearley, G. Mitra, and H. Williams. Analysis of mathematical programming problems prior to applying the simplex algorithm. *Mathematical Programming*, 8:54–83, 1975.

[19] E. Danna. Performance variability in mixed integer programming. Presentation slides from MIP 2008 workshop in New York City. `http://coral.ie.lehigh.edu/~jeff/mip-2008/program.pdf`, 2008.

[20] E. Danna, E. Rothberg, and C. L. Pape. Exploring relaxation induced neighborhoods to improve MIP solutions. *Mathematical Programming*, 102(1):71–90, 2004.

[21] B. Davey, N. Boland, and P. J. Stuckey. Efficient intelligent backtracking using linear programming. *INFORMS Journal of Computing*, 14(4):373–386, 2002.

[22] T. Davies, G. Gange, and P. J. Stuckey. Automatic logic-based Benders decomposition with minizinc. In *Proceedings of the 31st AAAI Conference on Artificial Intelligence (AAAI-17)*, pages 787–793. AAAI Press, 2017. https://aaai.org/ocs/index.php/AAAI/AAAI17/paper/view/14489.

[23] FICO Xpress Optimizer. `http://www.fico.com/en/Products/DMTools/xpress-overview/Pages/Xpress-Optimizer.aspx`.

[24] A. Gleixner, M. Bastubbe, L. Eifler, T. Gally, G. Gamrath, R. L. Gottwald, G. Hendel, C. Hojny, T. Koch, M. E. Lübbecke, S. J. Maher, M. Miltenberger, B. Müller, M. E. Pfetsch, C. Puchert, D. Rehfeldt, F. Schlösser, C. Schubert, F. Serrano, Y. Shinano, J. M. Viernickel, M. Walter, F. Wegscheider, J. T. Witt, and J. Witzig. The SCIP Optimization Suite 6.0. Technical Report 18-26, ZIB, Takustr. 7, 14195 Berlin, 2018.

[25] N. Jussien and V. Barichard. The PaLM system: explanation-based constraint programming. In *Proceedings of TRICS: Techniques foR Implementing Constraint programming Systems, a post-conference workshop of CP 2000*, pages 118–133, 2000.

[26] G. Katsirelos and F. Bacchus. Generalised nogoods in CSPs. In *Proceedings of AAAI-2005*, pages 390–396, 2005.

[27] G. Katsirelos and F. Bacchus. Generalized nogoods in csps. In *Proceedings, The Twentieth National Conference on Artificial Intelligence and the Seventeenth Innovative Applications of Artificial Intelligence Conference, July 9-13, 2005, Pittsburgh, Pennsylvania, USA*, pages 390–396. AAAI Press / The MIT Press, 2005.

[28] T. Koch, T. Achterberg, E. Andersen, O. Bastert, T. Berthold, R. E. Bixby, E. Danna, G. Gamrath, A. M. Gleixner, S. Heinz, A. Lodi, H. Mittelmann, T. Ralphs, D. Salvagnin, D. E. Steffy, and K. Wolter. MIPLIB 2010. *Mathematical Programming Computation*, 3(2):103–163, 2011.

[29] C. M. Li and Anbulagan. Look-ahead versus look-back for satisfiability problems. In *Proc. of CP*, pages 342–356, Autriche, 1997. Springer.

[30] A. Lodi and A. Tramontani. Performance variability in mixed-integer programming. In *Theory Driven by Influential Applications*, pages 1–12. INFORMS, 2013.

[31] J. P. Marques-Silva and K. A. Sakallah. GRASP: A search algorithm for propositional satisfiability. *IEEE Transactions of Computers*, 48:506–521, 1999.

[32] R. McGill, J. W. Tukey, and W. A. Larsen. Variations of box plots. *The American Statistician*, 32(1):12–16, 1978.

[33] M. H. Moskewicz, C. F. Madigan, Y. Zhao, L. Zhang, and S. Malik. Chaff: Engineering an efficient SAT solver. In *Proceedings of DAC'01*, pages 530–535, 2001.

[34] A. Nadel and V. Ryvchin. Efficient SAT solving under assumptions. In *Proc SAT*, volume 7317 of *LNCS*, pages 242–255. Springer, 2012.

[35] O. Ohrimenko, P. J. Stuckey, and M. Codish. Propagation via lazy clause generation. *Constraints*, 14(3):357–391, 2009.

[36] P. Refalo. Impact-based search strategies for constraint programming. In M. Wallace, editor, *Principles and Practice of Constraint Programming – CP 2004*, pages 557–571, Berlin, Heidelberg, 2004. Springer Berlin Heidelberg.

[37] R. Rodosek, M. G. Wallace, and M. T. Hajian. A new approach to integrating mixed integer programming and constraint logic programming. *Annals of Operations Research*, 86(1):63–87, 1999.

[38] T. Sandholm and R. Shields. Nogood learning for mixed integer programming. In *Workshop on Hybrid Methods and Branching Rules in Combinatorial Optimization, Montréal*, 2006.

[39] M. W. P. Savelsbergh. Preprocessing and probing techniques for mixed integer programming problems. *ORSA Journal on Computing*, 6:445–454, 1994.

[40] R. M. Stallman and G. J. Sussman. Forward reasoning and dependency-directed backtracking in a system for computer-aided circuit analysis. *Artificial Intelligence*, 9(2):135–196, 1977.

[41] J. Witzig, T. Berthold, and S. Heinz. Experiments with conflict analysis in mixed integer programming. In *Integration of AI and OR Techniques in Constraint Programming for Combinatorial Optimization Problems, 14th International Conference, CPAIOR 2017*, volume 10335 of *LNCS*, pages 211–220. Springer Berlin Heidelberg, May 2017.

[42] T. H. Yunes, I. D. Aron, and J. N. Hooker. An integrated solver for optimization problems. *Operations Research*, 58(2):342–356, 2010.

[43] L. Zhang, C. F. Madigan, M. H. Moskewicz, and S. Malik. Efficient conflict driven learning in a Boolean satisfiability solver. In *Proceedings of the 2001 IEEE/ACM international conference on Computer-aided design*, pages 279–285. IEEE Press, 2001.

Table 4: Detailed computational results on `MMM-IP` with five random seeds. Relative changes by at least 5% are highlighted in bold and blue (**improvement**) or italic and red (*deterioration*).

| | default | | degeneracy | | dualbound | | leaves | | obj | | nsols | | sblps | | nochecks | | onlyroot | |
|---|---|---|---|---|---|---|---|---|---|---|---|---|---|---|---|---|---|---|
| Instance | time | nodes | $time_Q$ | $nodes_Q$ | $time_Q$ | $nodes_Q$ | $time_Q$ | $nodes_Q$ | $time_Q$ | $nodes_Q$ | $time_Q$ | $nodes_Q$ | $time_Q$ | $nodes_Q$ | $time_Q$ | $nodes_Q$ | $time_Q$ | $nodes_Q$ |
| 10teams | 2.14 | 1 | *1.277* | 1.000 | 0.990 | 1.000 | 0.994 | 1.000 | 0.990 | 1.000 | *1.271* | 1.000 | 0.994 | 1.000 | *1.248* | 1.000 | *1.258* | 1.000 |
| | 20.00 | 1369 | **0.241** | **0.072** | *1.290* | **0.687** | *1.236* | 1.001 | 0.995 | 1.000 | **0.244** | **0.072** | 0.963 | **0.303** | **0.245** | **0.072** | **0.246** | **0.072** |
| | 15.57 | 287 | *1.607* | *2.729* | *2.062* | *4.067* | 0.998 | 1.000 | 1.006 | 1.000 | *1.308* | *2.729* | *2.080* | *4.238* | *1.341* | **0.829** | **0.672** | **0.323** |
| | 4.16 | 5 | *2.721* | *1.219* | 1.000 | 1.000 | 1.006 | 1.000 | 0.996 | 1.000 | *2.713* | *1.219* | 1.004 | 1.000 | *2.709* | *1.219* | *2.329* | *1.219* |
| | 4.94 | 5 | **0.909** | 1.000 | 1.000 | 1.000 | 1.002 | 1.000 | 1.007 | 1.000 | **0.909** | 1.000 | 0.995 | 1.000 | **0.912** | 1.000 | **0.914** | 1.000 |
| 30n20b8 | 314.13 | 195 | **0.605** | **0.553** | 1.002 | 1.000 | 1.001 | 1.000 | 0.990 | 1.000 | **0.289** | **0.346** | 0.999 | 1.000 | **0.285** | **0.346** | **0.289** | **0.346** |
| | 157.33 | 71 | **0.712** | **0.596** | 0.999 | 1.000 | 0.998 | 1.000 | 0.995 | 1.000 | **0.810** | **0.596** | 0.999 | 1.000 | **0.789** | **0.596** | **0.852** | **0.596** |
| | 87.25 | 4 | 1.001 | 1.000 | 1.000 | 1.000 | 1.003 | 1.000 | 0.998 | 1.000 | *1.191* | 1.000 | 0.997 | 1.000 | *1.180* | 0.981 | *1.188* | 1.000 |
| | 213.06 | 102 | *1.208* | *1.312* | 0.996 | 1.000 | 0.991 | 1.000 | 0.978 | 0.950 | *1.251* | **0.896** | 0.993 | 1.000 | 1.047 | **0.782** | *1.251* | **0.896** |
| | 470.09 | 290 | **0.811** | **0.723** | 0.998 | 1.000 | 1.003 | 1.000 | 0.998 | 1.000 | **0.238** | **0.423** | 0.996 | 1.000 | **0.241** | **0.421** | **0.238** | **0.423** |
| acc-tight5 | 125.35 | 172 | *1.114* | *1.171* | **0.387** | **0.295** | **0.423** | **0.302** | 1.011 | 0.965 | **0.651** | **0.735** | 0.991 | 1.000 | *1.087* | *1.171* | *1.115* | *1.622* |
| | 84.38 | 98 | *1.310* | *2.072* | *1.512* | *1.286* | **0.448** | **0.326** | **0.856** | **0.551** | *1.786* | *2.650* | 0.991 | 1.000 | *1.792* | *2.650* | **0.721** | 0.952 |
| | 50.45 | 480 | 0.973 | **0.929** | *2.656* | *2.631* | *1.753* | *1.966* | 1.022 | 0.993 | 0.973 | **0.929** | 1.002 | 1.000 | 0.974 | **0.929** | *3.299* | *6.017* |
| | 85.88 | 87 | **0.340** | **0.189** | *1.590* | *1.794* | **0.401** | **0.290** | **0.420** | **0.274** | **0.339** | **0.189** | 1.000 | 1.000 | **0.340** | **0.189** | *1.464* | *1.891* |
| | 73.38 | 864 | **0.417** | **0.161** | *1.440* | *1.091* | *2.064* | *1.326* | *1.488* | *1.309* | **0.417** | **0.161** | 0.996 | 1.000 | **0.419** | **0.161** | *1.256* | *1.088* |
| air03 | 1.70 | 1 | 1.011 | 1.000 | 1.015 | 1.000 | 1.015 | 1.000 | 1.015 | 1.000 | 1.004 | 1.000 | 1.007 | 1.000 | *6.433* | 1.000 | *6.515* | 1.000 |
| | 1.71 | 1 | 1.004 | 1.000 | 1.004 | 1.000 | 1.011 | 1.000 | 1.004 | 1.000 | 0.996 | 1.000 | 1.007 | 1.000 | *6.590* | 1.000 | *6.472* | 1.000 |
| | 1.74 | 1 | 0.993 | 1.000 | 1.011 | 1.000 | 1.004 | 1.000 | 0.989 | 1.000 | 1.000 | 1.000 | 0.985 | 1.000 | *6.365* | 1.000 | *6.372* | 1.000 |
| | 1.75 | 1 | 1.000 | 1.000 | 0.996 | 1.000 | 1.000 | 1.000 | 0.996 | 1.000 | 1.000 | 1.000 | 1.000 | 1.000 | *6.378* | 1.000 | *6.335* | 1.000 |
| | 1.75 | 1 | 0.996 | 1.000 | 0.993 | 1.000 | 0.985 | 1.000 | 0.985 | 1.000 | 0.993 | 1.000 | 1.004 | 1.000 | *6.491* | 1.000 | *6.378* | 1.000 |
| air04 | 49.46 | 110 | 1.003 | 1.000 | 1.007 | 1.000 | 1.008 | 1.000 | 1.005 | 1.000 | 1.003 | 1.000 | 1.006 | 1.000 | 0.978 | **0.776** | **0.911** | **0.652** |
| | 36.09 | 39 | 1.008 | 1.000 | 1.003 | 1.000 | 1.004 | 1.000 | 0.996 | 1.000 | 1.000 | 1.000 | 1.002 | 1.000 | *1.172* | *1.144* | *1.112* | *1.158* |
| | 52.88 | 118 | 1.005 | 1.000 | 1.000 | 1.000 | 1.002 | 1.000 | 0.999 | 1.000 | 0.997 | 1.000 | 0.990 | 1.000 | 0.993 | **0.665** | **0.949** | **0.665** |
| | 36.56 | 31 | 1.014 | 1.000 | 1.007 | 1.000 | 1.010 | 1.000 | 1.001 | 1.000 | 0.995 | 1.000 | 1.010 | 1.000 | *1.552* | *1.481* | *1.735* | *3.061* |
| | 48.35 | 122 | 1.006 | 1.000 | 1.002 | 1.000 | 1.004 | 1.000 | 1.001 | 1.000 | 1.004 | 1.000 | 1.003 | 1.000 | **0.803** | **0.500** | **0.801** | **0.500** |
| air05 | 24.82 | 189 | 0.990 | 1.000 | 0.996 | 1.000 | 0.999 | 1.000 | 1.006 | 1.000 | 1.003 | 1.000 | 0.990 | 1.000 | *2.135* | *1.388* | *1.270* | *1.547* |
| | 26.54 | 239 | 1.004 | 1.000 | 1.005 | 1.000 | 1.001 | 1.000 | 1.008 | 1.000 | 1.001 | 1.000 | 1.003 | 1.000 | *3.006* | *1.714* | *1.275* | *1.419* |
| | 31.37 | 438 | 0.992 | 1.000 | 0.992 | 1.000 | 0.997 | 1.000 | 0.996 | 1.000 | 0.994 | 1.000 | 0.997 | 1.000 | *1.869* | **0.690** | *1.132* | **0.946** |
| | 27.36 | 283 | 1.007 | 1.000 | 1.006 | 1.000 | 1.006 | 1.000 | 1.007 | 1.000 | 1.006 | 1.000 | 0.996 | 1.000 | *1.679* | *1.222* | *1.175* | *1.204* |
| | 27.14 | 279 | 1.003 | 1.000 | 1.003 | 1.000 | 1.002 | 1.000 | 1.002 | 1.000 | 1.004 | 1.000 | 0.996 | 1.000 | *2.084* | *1.343* | *1.128* | **0.934** |
| blend2 | 0.69 | 714 | 1.041 | 1.000 | 1.041 | 1.000 | 1.036 | 1.000 | 1.047 | 1.000 | 1.041 | 1.000 | 1.047 | 1.000 | 1.041 | 1.000 | 1.041 | 1.000 |
| | 1.08 | 973 | 0.986 | 1.000 | 1.005 | 1.000 | 0.971 | 1.000 | 1.019 | 1.000 | 1.005 | 1.000 | 0.971 | 1.000 | 1.010 | 1.000 | 0.986 | 1.000 |
| | 1.05 | 912 | 1.005 | 1.000 | 0.995 | 1.000 | 0.990 | 1.000 | 1.000 | 1.000 | 0.985 | 1.000 | 1.000 | 1.000 | 1.000 | 1.000 | 0.995 | 1.000 |
| | 0.67 | 543 | 0.994 | 1.000 | 1.006 | 1.000 | 1.012 | 1.000 | 1.030 | 1.000 | 1.000 | 1.000 | 1.012 | 1.000 | 1.012 | 1.000 | 1.018 | 1.000 |
| | 0.66 | 709 | 1.030 | 1.000 | 1.018 | 1.000 | 1.018 | 1.000 | 1.006 | 1.000 | 1.036 | 1.000 | 1.042 | 1.000 | 1.024 | 1.000 | 1.030 | 1.000 |
| bley_xl1 | 164.49 | 1 | 0.997 | 1.000 | 1.002 | 1.000 | 0.998 | 1.000 | 1.000 | 1.000 | 1.008 | 1.000 | 0.996 | 1.000 | 0.997 | 1.000 | 1.003 | 1.000 |
| | 187.04 | 3 | 1.012 | *1.146* | 1.007 | 1.000 | 1.005 | 1.000 | 1.002 | 1.000 | 1.005 | 1.000 | 1.003 | 1.000 | **0.906** | 0.981 | **0.907** | 0.981 |
| | 182.51 | 13 | **0.892** | **0.903** | 1.012 | 1.000 | 1.015 | 1.000 | 1.012 | 1.000 | 1.013 | 1.000 | 1.001 | 0.956 | **0.901** | **0.903** | **0.898** | **0.903** |
| | 133.90 | 1 | 1.006 | 1.000 | 1.008 | 1.000 | 1.003 | 1.000 | 1.005 | 1.000 | 1.005 | 1.000 | 1.011 | 1.000 | 0.982 | 1.000 | *1.118* | 1.000 |

cont. on next page

Table 4: Detailed computational results on `MMM-IP` with five random seeds. Relative changes by at least 5% are highlighted in bold and blue (**improvement**) or italic and red (*deterioration*).

| Instance | default | | degeneracy | | dualbound | | leaves | | obj | | nsols | | sblps | | nochecks | | onlyroot | |
|---|---|---|---|---|---|---|---|---|---|---|---|---|---|---|---|---|---|---|
| | time | nodes | $time_Q$ | $nodes_Q$ | $time_Q$ | $nodes_Q$ | $time_Q$ | $nodes_Q$ | $time_Q$ | $nodes_Q$ | $time_Q$ | $nodes_Q$ | $time_Q$ | $nodes_Q$ | $time_Q$ | $nodes_Q$ | $time_Q$ | $nodes_Q$ |
| | 135.56 | 1 | 1.006 | 1.000 | 1.006 | 1.000 | 1.001 | 1.000 | 1.000 | 1.000 | 1.005 | 1.000 | 1.008 | 1.000 | 1.014 | 1.000 | 1.023 | 1.000 |
| bnatt350 | 378.09 | 89 | **0.052** | **0.030** | **0.431** | **0.301** | **0.107** | **0.054** | *3.061* | *2.782* | **0.053** | **0.030** | **0.121** | **0.041** | **0.052** | **0.030** | *2.109* | *2.434* |
| | 521.76 | 326 | **0.033** | **0.022** | **0.276** | **0.183** | **0.144** | **0.084** | *1.081* | *1.382* | **0.033** | **0.022** | 1.011 | *1.431* | **0.033** | **0.022** | **0.487** | **0.420** |
| | 397.10 | 160 | **0.048** | **0.026** | **0.119** | **0.057** | **0.131** | **0.063** | **0.610** | **0.445** | **0.048** | **0.026** | **0.301** | **0.229** | **0.048** | **0.026** | **0.889** | **0.846** |
| | 155.46 | 4028 | **0.109** | **0.067** | **0.262** | **0.129** | *2.325* | *2.588* | *1.524* | *1.482* | **0.109** | **0.067** | **0.389** | **0.127** | **0.110** | **0.067** | *4.491* | *3.610* |
| | 994.16 | 2528 | **0.044** | **0.021** | **0.575** | **0.724** | **0.329** | **0.259** | **0.713** | **0.670** | **0.044** | **0.021** | **0.156** | **0.076** | **0.045** | **0.021** | **0.317** | **0.315** |
| cap6000 | 2.69 | 1830 | 1.003 | 1.000 | 1.019 | 1.000 | 1.003 | 1.000 | 1.019 | 1.000 | 1.016 | 1.000 | 0.989 | 1.000 | *1.412* | *1.131* | *1.098* | **0.866** |
| | 2.67 | 1837 | 0.986 | 1.000 | 0.984 | 1.000 | 0.989 | 1.000 | 0.986 | 1.000 | 0.975 | 1.000 | 1.003 | 1.000 | *1.455* | *1.256* | *1.155* | *1.142* |
| | 2.61 | 1959 | 0.992 | 1.000 | 1.000 | 1.000 | 0.997 | 1.000 | 1.006 | 1.000 | 0.997 | 1.000 | 0.992 | 1.000 | *1.280* | 0.982 | *1.150* | *1.066* |
| | 3.17 | 1924 | 1.007 | 1.000 | 0.990 | 1.000 | 1.007 | 1.000 | 0.993 | 1.000 | 0.995 | 1.000 | 0.995 | 1.000 | *1.422* | *1.166* | *1.086* | **0.894** |
| | 2.82 | 1762 | 1.008 | 1.000 | 1.008 | 1.000 | 0.995 | 1.000 | 0.987 | 1.000 | 0.995 | 1.000 | 1.008 | 1.000 | *1.387* | *1.325* | *1.157* | *1.258* |
| cov1075 | 3600.00 | 1445091 | 0.989 | *1.189* | 1.000 | 0.981 | 1.000 | 0.996 | 1.000 | 1.004 | 1.000 | 0.991 | 1.000 | 0.985 | 1.000 | **0.314** | 0.961 | *1.189* |
| | 3600.00 | 1787171 | 1.000 | 0.975 | 1.000 | 0.982 | 1.000 | 0.998 | 1.000 | 1.004 | 1.000 | 1.003 | 1.000 | 0.989 | 1.000 | **0.303** | 1.000 | 0.989 |
| | 3600.00 | 1528146 | 1.000 | *1.198* | 1.000 | 0.986 | 1.000 | 0.998 | 1.000 | 0.997 | 1.000 | 0.999 | 1.000 | 0.956 | 1.000 | **0.373** | 1.000 | *1.210* |
| | 3600.00 | 1808088 | 1.000 | 0.985 | 1.000 | 0.995 | 1.000 | 1.006 | 1.000 | 1.002 | 1.000 | 1.006 | 1.000 | 0.974 | 1.000 | **0.250** | 0.997 | 0.995 |
| | 3600.00 | 1201194 | 1.000 | *1.508* | 1.000 | 0.986 | 1.000 | 0.997 | 1.000 | 0.999 | 1.000 | 1.000 | 1.000 | **0.898** | 1.000 | **0.460** | 0.979 | *1.514* |
| csched010 | 3600.00 | 228369 | 1.000 | 1.003 | 1.000 | 1.001 | 1.000 | 0.963 | 1.000 | 0.999 | **0.758** | **0.644** | 1.000 | 1.001 | 1.000 | **0.571** | **0.907** | **0.925** |
| | 3600.00 | 175223 | 1.000 | 1.001 | 1.000 | 1.001 | 1.000 | **0.864** | 1.000 | 1.000 | 0.972 | 1.038 | 1.000 | 0.999 | 1.000 | **0.756** | **0.850** | 1.037 |
| | 2891.69 | 214060 | 1.000 | 1.000 | 0.951 | **0.941** | *1.147* | *1.206* | 0.997 | 1.000 | *1.245* | *1.145* | 1.002 | 1.000 | *1.245* | **0.689** | *1.245* | *1.280* |
| | 3600.00 | 195746 | 1.000 | 0.998 | 1.000 | 0.996 | 0.988 | **0.820** | 1.000 | 0.998 | 1.000 | 0.960 | 1.000 | 0.997 | 1.000 | **0.578** | 1.000 | 1.028 |
| | 3168.37 | 226703 | 1.001 | 1.000 | 0.999 | 1.000 | *1.136* | *1.303* | 1.000 | 1.000 | *1.136* | *1.352* | 1.001 | 1.000 | *1.136* | **0.682** | *1.136* | *1.307* |
| eil33-2 | 69.44 | 583 | 1.003 | 1.000 | 1.034 | 1.000 | 1.005 | 1.000 | 1.034 | 1.000 | 1.031 | 1.000 | 1.008 | 1.000 | *1.772* | **0.909** | *1.094* | **0.909** |
| | 81.92 | 739 | 1.007 | 1.000 | 1.007 | 1.000 | 1.010 | 1.000 | 0.997 | 1.000 | 1.009 | 1.000 | 1.007 | 1.000 | *1.671* | **0.852** | 1.028 | **0.852** |
| | 84.87 | 847 | 0.989 | 1.000 | 0.999 | 1.000 | 1.001 | 1.000 | 1.001 | 1.000 | 1.008 | 1.000 | 1.002 | 1.000 | *1.904* | **0.852** | *1.195* | **0.852** |
| | 73.60 | 543 | 0.988 | 1.000 | 1.007 | 1.000 | 1.000 | 1.000 | 1.001 | 1.000 | 1.010 | 1.000 | 1.003 | 1.000 | *2.006* | *1.509* | *1.260* | *2.728* |
| | 70.01 | 583 | 0.987 | 1.000 | 1.007 | 1.000 | 1.002 | 1.000 | 1.002 | 1.000 | 1.007 | 1.000 | 1.003 | 1.000 | *2.176* | *1.433* | *1.602* | *1.416* |
| eilB101 | 378.99 | 10582 | 0.998 | 1.000 | 1.021 | 1.000 | 1.002 | 1.000 | 1.023 | 1.000 | 0.999 | 1.000 | 1.009 | 1.000 | *1.083* | **0.724** | **0.705** | **0.724** |
| | 365.27 | 8427 | 0.984 | 1.000 | 0.995 | 1.000 | 0.993 | 1.000 | 0.999 | 1.000 | 1.005 | 1.000 | 0.997 | 1.000 | *1.391* | *1.150* | **0.916** | *1.050* |
| | 401.63 | 9358 | 0.995 | 1.000 | 1.000 | 1.000 | 0.998 | 1.000 | 0.993 | 1.000 | 1.008 | 1.000 | 0.997 | 1.000 | *1.147* | **0.780** | **0.819** | **0.780** |
| | 397.46 | 11465 | 0.996 | 1.000 | 0.999 | 1.000 | 0.998 | 1.000 | 0.999 | 1.000 | 1.002 | 1.000 | 0.995 | 1.000 | *1.304* | 1.025 | **0.905** | 1.025 |
| | 365.79 | 10006 | 0.989 | 1.000 | 0.998 | 1.000 | 0.989 | 1.000 | 0.994 | 1.000 | 0.996 | 1.000 | 1.001 | 1.000 | *1.201* | **0.894** | **0.818** | **0.894** |
| enigma | 0.50 | 809 | 1.000 | **0.825** | 1.020 | **0.716** | 1.000 | 1.000 | 1.000 | 0.998 | 1.000 | **0.542** | *1.140* | **0.686** | 1.000 | **0.542** | 1.000 | **0.856** |
| | 0.50 | 449 | 1.000 | **0.497** | 1.000 | **0.814** | 1.000 | 1.000 | 1.000 | 1.000 | 1.000 | **0.470** | 1.000 | **0.778** | 1.000 | **0.470** | 1.000 | *1.616* |
| | 0.50 | 541 | 1.000 | **0.254** | 1.000 | **0.510** | 1.000 | 1.000 | 1.000 | 0.977 | 1.000 | **0.237** | *1.060* | **0.721** | 1.000 | **0.237** | 1.000 | **0.665** |
| | 0.50 | 679 | *1.073* | **0.537** | 1.000 | **0.397** | 1.000 | 1.000 | 1.000 | 1.000 | 1.000 | **0.151** | *1.133* | *1.113* | 1.000 | **0.151** | 1.000 | **0.244** |
| | 0.50 | 445 | 1.020 | **0.543** | 1.000 | **0.833** | 1.000 | 1.000 | 1.000 | 1.006 | 1.000 | **0.543** | 1.000 | **0.806** | 1.013 | **0.543** | 1.000 | *1.675* |
| enlight13 | 0.50 | 1 | 1.000 | 1.000 | 1.000 | 1.000 | 1.000 | 1.000 | 1.000 | 1.000 | 1.000 | 1.000 | 1.000 | 1.000 | 1.000 | 1.000 | 1.000 | 1.000 |
| | 0.50 | 1 | 1.000 | 1.000 | 1.000 | 1.000 | 1.000 | 1.000 | 1.000 | 1.000 | 1.000 | 1.000 | 1.000 | 1.000 | 1.000 | 1.000 | 1.000 | 1.000 |
| | 0.50 | 1 | 1.000 | 1.000 | 1.000 | 1.000 | 1.000 | 1.000 | 1.000 | 1.000 | 1.000 | 1.000 | 1.000 | 1.000 | 1.000 | 1.000 | 1.000 | 1.000 |

cont. on next page

Table 4: Detailed computational results on `MMM-IP` with five random seeds. Relative changes by at least 5% are highlighted in bold and blue (**improvement**) or italic and red (*deterioration*).

| | default | | degeneracy | | dualbound | | leaves | | obj | | nsols | | sblps | | nochecks | | onlyroot | |
|---|---|---|---|---|---|---|---|---|---|---|---|---|---|---|---|---|---|---|
| Instance | time | nodes | time$_Q$ | nodes$_Q$ | time$_Q$ | nodes$_Q$ | time$_Q$ | nodes$_Q$ | time$_Q$ | nodes$_Q$ | time$_Q$ | nodes$_Q$ | time$_Q$ | nodes$_Q$ | time$_Q$ | nodes$_Q$ | time$_Q$ | nodes$_Q$ |
| | 0.50 | 1 | 1.000 | 1.000 | 1.000 | 1.000 | 1.000 | 1.000 | 1.000 | 1.000 | 1.000 | 1.000 | 1.000 | 1.000 | 1.000 | 1.000 | 1.000 | 1.000 |
| | 0.50 | 1 | 1.000 | 1.000 | 1.000 | 1.000 | 1.000 | 1.000 | 1.000 | 1.000 | 1.000 | 1.000 | 1.000 | 1.000 | 1.000 | 1.000 | 1.000 | 1.000 |
| enlight14 | 0.50 | 1 | 1.000 | 1.000 | 1.000 | 1.000 | 1.000 | 1.000 | 1.000 | 1.000 | 1.000 | 1.000 | 1.000 | 1.000 | 1.000 | 1.000 | 1.000 | 1.000 |
| | 0.50 | 1 | 1.000 | 1.000 | 1.000 | 1.000 | 1.000 | 1.000 | 1.000 | 1.000 | 1.000 | 1.000 | 1.000 | 1.000 | 1.000 | 1.000 | 1.000 | 1.000 |
| | 0.50 | 1 | 1.000 | 1.000 | 1.000 | 1.000 | 1.000 | 1.000 | 1.000 | 1.000 | 1.000 | 1.000 | 1.000 | 1.000 | 1.000 | 1.000 | 1.000 | 1.000 |
| | 0.50 | 1 | 1.000 | 1.000 | 1.000 | 1.000 | 1.000 | 1.000 | 1.000 | 1.000 | 1.000 | 1.000 | 1.000 | 1.000 | 1.000 | 1.000 | 1.000 | 1.000 |
| | 0.50 | 1 | 1.000 | 1.000 | 1.000 | 1.000 | 1.000 | 1.000 | 1.000 | 1.000 | 1.000 | 1.000 | 1.000 | 1.000 | 1.000 | 1.000 | 1.000 | 1.000 |
| fiber | 1.00 | 2 | 1.010 | 1.000 | 1.010 | 1.000 | 0.990 | 1.000 | 1.005 | 1.000 | 1.000 | 1.000 | 1.000 | 1.000 | *1.160* | 1.010 | *1.075* | 1.010 |
| | 1.42 | 4 | 0.988 | 1.000 | 1.004 | 1.000 | 0.988 | 1.000 | 1.000 | 1.000 | 0.992 | 1.000 | 1.008 | 1.000 | *1.384* | *1.240* | **0.938** | 0.990 |
| | 1.74 | 4 | 1.007 | 1.000 | 0.993 | 1.000 | 0.996 | 1.000 | 0.996 | 1.000 | 0.978 | 1.000 | 0.982 | 1.000 | *1.401* | 1.010 | **0.949** | 0.990 |
| | 1.53 | 4 | 1.008 | 1.000 | 1.000 | 1.000 | 0.996 | 1.000 | 0.996 | 1.000 | 0.988 | 1.000 | 1.000 | 1.000 | *1.194* | 0.990 | 1.024 | 1.000 |
| | 1.00 | 3 | 1.000 | 1.000 | 1.005 | 1.000 | 1.005 | 1.000 | 1.000 | 1.000 | 0.995 | 1.000 | 1.005 | 1.000 | *1.080* | 0.990 | 1.030 | 0.990 |
| flugpl | 0.50 | 1 | 1.000 | 1.000 | 1.000 | 1.000 | 1.000 | 1.000 | 1.000 | 1.000 | 1.000 | 1.000 | 1.000 | 1.000 | 1.000 | 1.000 | 1.000 | 1.000 |
| | 0.50 | 1 | 1.000 | 1.000 | 1.000 | 1.000 | 1.000 | 1.000 | 1.000 | 1.000 | 1.000 | 1.000 | 1.000 | 1.000 | 1.000 | 1.000 | 1.000 | 1.000 |
| | 0.50 | 1 | 1.000 | 1.000 | 1.000 | 1.000 | 1.000 | 1.000 | 1.000 | 1.000 | 1.000 | 1.000 | 1.000 | 1.000 | 1.000 | 1.000 | 1.000 | 1.000 |
| | 0.50 | 1 | 1.000 | 1.000 | 1.000 | 1.000 | 1.000 | 1.000 | 1.000 | 1.000 | 1.000 | 1.000 | 1.000 | 1.000 | 1.000 | 1.000 | 1.000 | 1.000 |
| | 0.50 | 1 | 1.000 | 1.000 | 1.000 | 1.000 | 1.000 | 1.000 | 1.000 | 1.000 | 1.000 | 1.000 | 1.000 | 1.000 | 1.000 | 1.000 | 1.000 | 1.000 |
| gt2 | 0.50 | 1 | 1.000 | 1.000 | 1.000 | 1.000 | 1.000 | 1.000 | 1.000 | 1.000 | 1.000 | 1.000 | 1.000 | 1.000 | 1.000 | 1.000 | 1.000 | 1.000 |
| | 0.50 | 1 | 1.000 | 1.000 | 1.000 | 1.000 | 1.000 | 1.000 | 1.000 | 1.000 | 1.000 | 1.000 | 1.000 | 1.000 | 1.000 | 1.000 | 1.000 | 1.000 |
| | 0.50 | 1 | 1.000 | 1.000 | 1.000 | 1.000 | 1.000 | 1.000 | 1.000 | 1.000 | 1.000 | 1.000 | 1.000 | 1.000 | 1.000 | 1.000 | 1.000 | 1.000 |
| | 0.50 | 1 | 1.000 | 1.000 | 1.000 | 1.000 | 1.000 | 1.000 | 1.000 | 1.000 | 1.000 | 1.000 | 1.000 | 1.000 | 1.000 | 1.000 | 1.000 | 1.000 |
| | 0.50 | 1 | 1.000 | 1.000 | 1.000 | 1.000 | 1.000 | 1.000 | 1.000 | 1.000 | 1.000 | 1.000 | 1.000 | 1.000 | 1.000 | 1.000 | 1.000 | 1.000 |
| harp2 | 691.56 | 1491878 | 1.005 | 1.000 | 1.001 | 1.000 | 1.001 | 1.000 | 1.007 | 1.000 | 1.000 | 1.000 | 1.003 | 1.000 | *5.200* | *5.260* | *4.215* | *4.834* |
| | 1159.66 | 2857677 | 1.004 | 1.000 | 0.998 | 1.000 | 0.994 | 1.000 | 0.998 | 1.000 | 0.993 | 1.000 | 0.994 | 1.000 | **0.410** | **0.228** | *2.318* | *2.253* |
| | 336.73 | 595854 | 0.992 | 1.000 | 1.007 | 1.000 | 1.007 | 1.000 | 1.005 | 1.000 | 1.002 | 1.000 | 0.999 | 1.000 | *1.604* | *1.232* | *1.348* | *1.300* |
| | 937.54 | 2046181 | 0.997 | 1.000 | 0.987 | 1.000 | 0.992 | 1.000 | 1.000 | 1.000 | 0.997 | 1.000 | 0.999 | 1.000 | **0.506** | **0.321** | *1.065* | *1.052* |
| | 425.68 | 923780 | 1.018 | 1.000 | 1.010 | 1.000 | 1.009 | 1.000 | 1.013 | 1.000 | 1.014 | 1.000 | 1.021 | 1.000 | *2.779* | *2.192* | *7.701* | *8.726* |
| iis-100-0-cov | 629.92 | 103302 | 0.993 | 1.000 | 0.993 | 1.000 | 1.014 | 1.000 | 0.991 | 1.000 | 0.997 | 1.000 | 0.992 | 1.000 | 0.999 | 0.999 | 0.965 | 0.999 |
| | 545.20 | 92393 | 0.999 | 1.000 | 0.997 | 1.000 | 1.001 | 1.000 | 1.002 | 1.000 | 0.996 | 1.000 | 0.998 | 1.000 | *1.128* | 1.046 | *1.053* | 1.046 |
| | 520.78 | 81639 | 0.999 | 1.000 | 1.000 | 1.000 | 1.002 | 1.000 | 1.001 | 1.000 | 1.000 | 1.000 | 0.997 | 1.000 | *1.151* | *1.128* | *1.097* | *1.128* |
| | 506.53 | 84800 | 1.002 | 1.000 | 0.997 | 1.000 | 0.995 | 1.000 | 0.996 | 1.000 | 0.998 | 1.000 | 0.998 | 1.000 | *1.117* | *1.080* | *1.077* | *1.080* |
| | 544.56 | 92313 | 0.996 | 1.000 | 0.999 | 1.000 | 1.004 | 1.000 | 1.001 | 1.000 | 1.002 | 1.000 | 0.998 | 1.000 | *1.089* | *1.052* | 0.990 | *1.052* |
| iis-bupa-cov | 1915.72 | 157619 | 1.001 | 1.000 | 1.009 | 1.000 | 1.002 | 1.000 | 1.000 | 1.000 | 1.007 | 1.000 | 1.003 | 1.000 | *1.704* | *1.376* | *1.344* | *1.376* |
| | 1944.83 | 164079 | 0.997 | 1.000 | 0.995 | 1.000 | 0.997 | 1.000 | 0.996 | 1.000 | 0.994 | 1.000 | 0.995 | 1.000 | *1.259* | *1.099* | *1.117* | *1.099* |
| | 1881.83 | 156888 | 0.998 | 1.000 | 0.999 | 1.000 | 0.993 | 1.000 | 0.998 | 1.000 | 0.996 | 1.000 | 0.996 | 1.000 | *1.189* | *1.093* | *1.081* | *1.093* |
| | 1959.19 | 161969 | 0.989 | 1.000 | 1.005 | 1.000 | 0.996 | 1.000 | 0.996 | 1.000 | 0.994 | 1.000 | 0.998 | 1.000 | *1.398* | *1.234* | *1.193* | *1.234* |
| | 1820.52 | 158284 | 0.990 | 1.000 | 0.997 | 1.000 | 0.994 | 1.000 | 0.995 | 1.000 | 0.995 | 1.000 | 0.994 | 1.000 | *1.378* | *1.219* | *1.241* | *1.219* |
| iis-pima-cov | 293.79 | 4824 | 0.995 | 1.000 | 0.997 | 1.000 | 1.002 | 1.000 | 1.002 | 1.000 | 0.993 | 1.000 | 0.997 | 1.000 | *1.072* | *1.109* | 1.016 | *1.109* |
| | 437.19 | 12189 | 0.995 | 1.000 | 0.997 | 1.000 | 0.994 | 1.000 | 1.001 | 1.000 | 1.003 | 1.000 | 0.999 | 1.000 | **0.746** | **0.445** | **0.712** | **0.445** |

cont. on next page

Table 4: Detailed computational results on `MMM-IP` with five random seeds. Relative changes by at least 5% are highlighted in bold and blue (**improvement**) or italic and red (*deterioration*).

| | default | | degeneracy | | dualbound | | leaves | | obj | | nsols | | sblps | | nochecks | | onlyroot | |
|---|---|---|---|---|---|---|---|---|---|---|---|---|---|---|---|---|---|---|
| Instance | **time** | **nodes** | **time$_Q$** | **nodes$_Q$** | **time$_Q$** | **nodes$_Q$** | **time$_Q$** | **nodes$_Q$** | **time$_Q$** | **nodes$_Q$** | **time$_Q$** | **nodes$_Q$** | **time$_Q$** | **nodes$_Q$** | **time$_Q$** | **nodes$_Q$** | **time$_Q$** | **nodes$_Q$** |
| | 246.99 | 4668 | 1.000 | 1.000 | 1.001 | 1.000 | 1.004 | 1.000 | 1.006 | 1.000 | 1.005 | 1.000 | 1.004 | 1.000 | *1.326* | *1.583* | *1.214* | *1.583* |
| | 284.88 | 6607 | 1.002 | 1.000 | 1.003 | 1.000 | 1.000 | 1.000 | 1.001 | 1.000 | 0.996 | 1.000 | 0.998 | 1.000 | *1.225* | *1.265* | *1.144* | *1.265* |
| | 556.78 | 17821 | 1.003 | 1.000 | 0.998 | 1.000 | 1.000 | 1.000 | 0.999 | 1.000 | 0.999 | 1.000 | 0.995 | 1.000 | **0.594** | **0.343** | **0.562** | **0.343** |
| l152lav | 5.51 | 23 | 0.994 | 1.000 | 0.989 | 1.000 | 0.994 | 1.000 | 0.994 | 1.000 | 0.995 | 1.000 | 1.005 | 1.000 | *2.278* | *1.520* | 1.035 | *1.244* |
| | 1.67 | 31 | 0.985 | 1.000 | 0.993 | 1.000 | 0.985 | 1.000 | 0.993 | 1.000 | 0.989 | 1.000 | 0.966 | 1.000 | *2.131* | *1.122* | *1.225* | *1.183* |
| | 2.93 | 50 | 1.000 | 1.000 | 0.995 | 1.000 | 0.992 | 1.000 | 1.000 | 1.000 | 1.003 | 1.000 | 1.008 | 1.000 | *3.046* | *1.367* | 0.962 | 0.967 |
| | 2.34 | 20 | 0.985 | 1.000 | 0.994 | 1.000 | 0.991 | 1.000 | 0.991 | 1.000 | 0.985 | 1.000 | 1.003 | 1.000 | *3.713* | *2.008* | *1.275* | *1.342* |
| | 3.02 | 52 | 0.995 | 1.000 | 1.000 | 1.000 | 0.995 | 1.000 | 0.993 | 1.000 | 0.998 | 1.000 | 1.002 | 1.000 | *2.167* | *1.204* | **0.891** | **0.921** |
| lectsched-4-obj | 27.04 | 805 | **0.753** | **0.215** | *1.223* | **0.740** | 1.004 | 1.000 | 1.009 | **0.924** | 1.005 | 1.000 | **0.730** | **0.314** | **0.748** | **0.215** | 0.978 | *1.608* |
| | 23.16 | 868 | *1.201* | **0.288** | **0.813** | **0.525** | 1.029 | *1.080* | **0.887** | **0.765** | 0.988 | 1.000 | **0.666** | **0.311** | *1.201* | **0.288** | 0.995 | *1.282* |
| | 86.08 | 1061 | **0.288** | **0.264** | **0.416** | **0.216** | **0.627** | **0.461** | 1.020 | 1.006 | **0.385** | **0.696** | **0.458** | **0.146** | **0.265** | **0.121** | **0.384** | **0.696** |
| | 26.83 | 619 | **0.327** | **0.082** | **0.768** | **0.367** | **0.928** | **0.540** | **0.926** | **0.772** | **0.257** | **0.082** | *2.333* | **0.477** | **0.330** | **0.082** | **0.257** | **0.082** |
| | 22.73 | 745 | **0.909** | **0.253** | 1.038 | **0.805** | 1.005 | 1.000 | 0.992 | 0.954 | 1.004 | 1.000 | *1.484* | **0.910** | **0.913** | **0.253** | **0.701** | **0.647** |
| lseu | 0.75 | 268 | 1.011 | 1.000 | 1.011 | 1.000 | 1.006 | 1.000 | 1.000 | 1.000 | 1.011 | 1.000 | 1.034 | 1.000 | *1.566* | *1.095* | 1.006 | 0.986 |
| | 0.50 | 4 | 1.000 | 1.000 | 1.000 | 1.000 | 1.000 | 1.000 | 1.000 | 1.000 | 1.000 | 1.000 | 1.000 | 1.000 | 1.000 | 1.010 | 1.000 | 1.000 |
| | 0.51 | 205 | 0.993 | 1.000 | 1.000 | 1.000 | 0.993 | 1.000 | 0.993 | 1.000 | 0.993 | 1.000 | 1.000 | 1.000 | *1.238* | **0.866** | 0.993 | **0.931** |
| | 0.50 | 91 | 1.000 | 1.000 | 1.000 | 1.000 | 1.000 | 1.000 | 1.000 | 1.000 | 1.000 | 1.000 | 1.000 | 1.000 | 1.000 | **0.738** | 1.000 | 0.969 |
| | 0.50 | 191 | 1.007 | 1.000 | 1.000 | 1.000 | 1.000 | 1.000 | 1.007 | 1.000 | 1.007 | 1.000 | 1.007 | 1.000 | *1.380* | **0.670** | 1.000 | **0.680** |
| m100n500k4r1 | 3600.00 | 2166065 | 1.000 | *1.118* | 1.000 | **0.786** | 1.000 | 0.998 | 1.000 | 1.001 | 1.000 | 1.001 | 1.000 | **0.777** | 1.000 | **0.780** | 1.000 | *1.171* |
| | 3600.00 | 2266017 | 1.000 | 1.035 | 1.000 | **0.712** | 1.000 | 0.995 | 1.000 | 0.999 | 1.000 | 0.998 | 1.000 | **0.701** | 1.000 | **0.782** | 1.000 | *1.059* |
| | 3600.00 | 2611447 | 1.000 | **0.884** | 1.000 | **0.646** | 1.000 | 0.998 | 1.000 | 1.002 | 1.000 | 1.001 | 1.000 | **0.640** | 1.000 | **0.559** | 1.000 | **0.930** |
| | 3600.00 | 2305045 | 1.000 | **0.937** | 1.000 | **0.723** | 1.000 | 0.998 | 1.000 | 0.999 | 1.000 | 0.999 | 1.000 | **0.723** | 1.000 | **0.640** | 1.000 | 0.989 |
| | 3600.00 | 2226238 | 1.000 | **0.833** | 1.000 | **0.733** | 1.000 | 0.991 | 1.000 | 0.999 | 1.000 | 0.997 | 1.000 | **0.732** | 1.000 | **0.645** | 1.000 | **0.923** |
| macrophage | 125.50 | 5163 | 1.001 | 1.000 | 1.006 | 1.000 | 1.001 | 1.000 | 1.005 | 1.000 | 1.001 | 1.000 | 1.001 | 1.000 | *1.774* | *1.158* | **0.911** | **0.701** |
| | 157.18 | 7241 | 0.999 | 1.000 | 0.999 | 1.000 | 0.998 | 1.000 | 1.000 | 1.000 | 0.999 | 1.000 | 0.998 | 1.000 | *1.535* | *1.184* | **0.948** | **0.852** |
| | 209.11 | 9131 | 1.000 | 1.000 | 0.999 | 1.000 | 1.002 | 1.000 | 0.998 | 1.000 | 0.999 | 1.000 | 0.999 | 1.000 | *1.438* | *1.250* | **0.821** | *1.053* |
| | 162.53 | 7927 | 1.000 | 1.000 | 1.002 | 1.000 | 1.002 | 1.000 | 1.000 | 1.000 | 0.998 | 1.000 | 1.000 | 1.000 | *1.876* | *1.088* | *1.083* | **0.840** |
| | 121.37 | 5239 | 0.999 | 1.000 | 0.997 | 1.000 | 0.998 | 1.000 | 1.000 | 1.000 | 1.001 | 1.000 | 0.997 | 1.000 | *1.987* | *1.355* | *1.059* | *1.071* |
| manna81 | 0.50 | 1 | 1.000 | 1.000 | 1.000 | 1.000 | 1.000 | 1.000 | 1.000 | 1.000 | 1.000 | 1.000 | 1.000 | 1.000 | 1.007 | 1.000 | 1.020 | 1.000 |
| | 0.50 | 1 | 1.000 | 1.000 | 1.000 | 1.000 | 1.000 | 1.000 | 1.000 | 1.000 | 1.000 | 1.000 | 1.000 | 1.000 | 1.013 | 1.000 | 1.013 | 1.000 |
| | 0.50 | 1 | 1.000 | 1.000 | 1.000 | 1.000 | 1.000 | 1.000 | 1.000 | 1.000 | 1.000 | 1.000 | 1.000 | 1.000 | 1.013 | 1.000 | 1.033 | 1.000 |
| | 0.50 | 1 | 1.000 | 1.000 | 1.000 | 1.000 | 1.000 | 1.000 | 1.000 | 1.000 | 1.000 | 1.000 | 1.000 | 1.000 | 1.000 | 1.000 | 1.013 | 1.000 |
| | 0.50 | 1 | 1.000 | 1.000 | 1.000 | 1.000 | 1.000 | 1.000 | 1.000 | 1.000 | 1.000 | 1.000 | 1.000 | 1.000 | 1.007 | 1.000 | 1.033 | 1.000 |
| markshare1 | 3600.00 | 19010451 | 1.000 | **0.931** | 1.000 | **0.786** | 1.000 | 0.996 | 1.000 | 1.000 | 1.000 | 0.996 | 1.000 | **0.797** | 1.000 | **0.924** | 1.000 | *1.083* |
| | 3600.00 | 20192578 | 1.000 | **0.862** | 1.000 | **0.893** | 1.000 | 0.999 | 1.000 | 0.997 | 1.000 | 1.000 | 1.000 | **0.890** | 1.000 | **0.856** | 1.000 | 1.017 |
| | 3600.00 | 19529435 | 1.000 | **0.919** | 1.000 | **0.907** | 1.000 | 1.000 | 1.000 | 0.998 | 1.000 | 0.997 | 1.000 | **0.919** | 1.000 | **0.918** | 1.000 | *1.219* |
| | 3600.00 | 20381558 | 1.000 | **0.898** | 1.000 | **0.865** | 1.000 | 0.998 | 1.000 | 1.001 | 1.000 | 1.002 | 1.000 | **0.818** | 1.000 | **0.903** | 1.000 | 1.017 |
| | 3600.00 | 21608999 | 1.000 | **0.797** | 1.000 | **0.808** | 1.000 | 1.004 | 1.000 | 1.000 | 1.000 | 0.996 | 1.000 | **0.806** | 1.000 | **0.795** | 1.000 | *1.090* |
| markshare2 | 3600.00 | 4050504 | 1.000 | **0.843** | 1.000 | *1.137* | 1.000 | 0.995 | 1.000 | 0.991 | 1.000 | 0.995 | 1.000 | *1.143* | 1.000 | **0.846** | 1.000 | 0.958 |

cont. on next page

Table 4: Detailed computational results on `MMM-IP` with five random seeds. Relative changes by at least 5% are highlighted in bold and blue (**improvement**) or italic and red (*deterioration*).

| | default | | degeneracy | | dualbound | | leaves | | obj | | nsols | | sblps | | nochecks | | onlyroot | |
|---|---|---|---|---|---|---|---|---|---|---|---|---|---|---|---|---|---|---|---|---|
| Instance | time | nodes | time$_Q$ | nodes$_Q$ | time$_Q$ | nodes$_Q$ | time$_Q$ | nodes$_Q$ | time$_Q$ | nodes$_Q$ | time$_Q$ | nodes$_Q$ | time$_Q$ | nodes$_Q$ | time$_Q$ | nodes$_Q$ | time$_Q$ | nodes$_Q$ |
| | 3600.00 | 5646401 | 1.000 | **0.697** | 1.000 | **0.783** | 1.000 | 1.001 | 1.000 | 1.000 | 1.000 | 1.000 | 1.000 | **0.787** | 1.000 | **0.697** | 1.000 | **0.636** |
| | 3600.00 | 3589038 | 1.000 | 1.026 | 1.000 | 1.002 | 1.000 | 1.002 | 1.000 | 0.996 | 1.000 | 0.998 | 1.000 | 1.002 | 1.000 | 1.026 | 1.000 | *1.083* |
| | 3600.00 | 4244067 | 1.000 | *1.069* | 1.000 | *1.098* | 1.000 | 1.001 | 1.000 | 1.003 | 1.000 | 1.001 | 1.000 | *1.106* | 1.000 | *1.069* | 1.000 | **0.891** |
| | 3600.00 | 5150867 | 1.000 | **0.830** | 1.000 | **0.737** | 1.000 | 1.002 | 1.000 | 0.996 | 1.000 | 0.999 | 1.000 | **0.736** | 1.000 | **0.832** | 1.000 | **0.948** |
| mcsched | 285.92 | 13689 | 0.999 | 1.000 | 0.997 | 1.000 | 0.998 | 1.000 | 1.012 | 1.000 | 0.998 | 1.000 | 1.000 | 1.000 | *1.093* | **0.678** | **0.820** | **0.732** |
| | 271.88 | 11373 | 0.999 | 1.000 | 0.998 | 1.000 | 1.002 | 1.000 | 0.999 | 1.000 | 1.001 | 1.000 | 1.000 | 1.000 | *3.196* | *1.747* | *1.459* | *1.710* |
| | 233.97 | 10059 | 0.996 | 1.000 | 0.997 | 1.000 | 1.005 | 1.000 | 0.998 | 1.000 | 0.997 | 1.000 | 1.001 | 1.000 | *2.468* | *1.296* | *1.287* | *1.310* |
| | 286.25 | 12496 | 0.997 | 1.000 | 0.998 | 1.000 | 0.995 | 1.000 | 0.998 | 1.000 | 0.998 | 1.000 | 0.997 | 1.000 | *2.450* | *1.245* | *1.055* | *1.146* |
| | 270.69 | 10373 | 1.005 | 1.000 | 1.000 | 1.000 | 1.005 | 1.000 | 1.000 | 1.000 | 1.001 | 1.000 | 1.001 | 1.000 | *1.542* | 0.996 | 1.007 | 0.992 |
| mine-166-5 | 58.69 | 1569 | 0.995 | 1.000 | 0.995 | 1.000 | 0.997 | 1.000 | 0.993 | 1.000 | 1.020 | 1.000 | 1.001 | 1.000 | **0.776** | 1.022 | **0.851** | *2.359* |
| | 57.78 | 1191 | 0.990 | 1.000 | 0.992 | 1.000 | 0.996 | 1.000 | 0.990 | 1.000 | 0.990 | 1.000 | 0.991 | 1.000 | **0.920** | **0.930** | *1.110* | **0.901** |
| | 46.45 | 2244 | 1.001 | 1.000 | 1.004 | 1.000 | 1.006 | 1.000 | 1.007 | 1.000 | 1.001 | 1.000 | 1.003 | 1.000 | *1.103* | *1.067* | *1.132* | **0.650** |
| | 64.93 | 139 | 1.006 | 1.000 | 1.004 | 1.000 | 1.005 | 1.000 | 1.001 | 1.000 | 1.005 | 1.000 | 1.002 | 1.000 | **0.711** | *5.397* | **0.619** | *6.515* |
| | 42.28 | 2863 | 1.000 | 1.000 | 1.000 | 1.000 | 1.001 | 1.000 | 1.000 | 1.000 | 1.001 | 1.000 | 0.998 | 1.000 | *1.165* | **0.555** | 0.994 | **0.873** |
| mine-90-10 | 276.07 | 51926 | 1.002 | 1.000 | 1.001 | 1.000 | *1.157* | *1.090* | 0.999 | 1.000 | 1.000 | 1.000 | 1.003 | 1.000 | 1.033 | **0.750** | **0.767** | **0.596** |
| | 125.32 | 22195 | 0.998 | 1.000 | 1.005 | 1.000 | 0.998 | 1.000 | 0.999 | 1.000 | 1.003 | 1.000 | 1.001 | 1.000 | *4.652* | *3.660* | *1.749* | *1.438* |
| | 173.25 | 25979 | 1.001 | 1.000 | 1.004 | 1.000 | 0.999 | 1.000 | 0.999 | 1.000 | 1.004 | 1.000 | 0.999 | 1.000 | *1.687* | *1.060* | *2.725* | *2.264* |
| | 500.83 | 92879 | 1.003 | 1.000 | 1.000 | 1.000 | *1.324* | *1.061* | 0.999 | 1.000 | 0.999 | 1.000 | 0.999 | 1.000 | **0.809** | **0.459** | **0.556** | **0.457** |
| | 236.18 | 40409 | 1.005 | 1.000 | 1.003 | 1.000 | **0.907** | *1.242* | 1.004 | 1.000 | 1.002 | 1.000 | 1.013 | 1.000 | *2.162* | *1.767* | **0.762** | **0.592** |
| misc03 | 1.39 | 398 | 1.004 | 1.000 | 1.004 | 1.000 | 1.000 | 1.000 | 1.013 | 1.000 | 1.004 | 1.000 | 0.983 | 1.000 | **0.799** | **0.267** | **0.711** | **0.263** |
| | 0.71 | 15 | 1.000 | 1.000 | 0.994 | 1.000 | 1.012 | 1.000 | 1.006 | 1.000 | *1.076* | *1.226* | 1.006 | 1.000 | *1.427* | *1.226* | *1.094* | *1.226* |
| | 0.70 | 19 | 1.018 | 1.000 | 0.994 | 1.000 | 1.006 | 1.000 | 1.018 | 1.000 | *1.165* | *1.487* | 0.976 | 1.000 | *1.512* | *1.471* | *1.176* | *1.487* |
| | 1.00 | 111 | 0.990 | 1.000 | 1.010 | 1.000 | 1.010 | 1.000 | 1.000 | 1.000 | 0.985 | 1.000 | 1.000 | 1.000 | 1.015 | **0.763** | **0.855** | **0.706** |
| | 0.66 | 15 | 0.988 | 1.000 | 1.000 | 1.000 | 0.988 | 1.000 | 0.988 | 1.000 | *1.187* | *1.609* | 0.988 | 1.000 | *1.608* | *1.435* | *1.151* | *1.609* |
| misc07 | 38.05 | 28289 | 1.006 | 1.000 | 1.027 | 1.000 | 0.998 | 1.000 | 0.998 | 1.000 | *1.089* | *1.188* | 1.000 | 1.000 | *1.496* | **0.758** | *1.087* | *1.188* |
| | 36.66 | 29782 | 0.995 | 1.000 | 0.996 | 1.000 | 0.998 | 1.000 | 0.999 | 1.000 | 1.001 | 1.000 | 1.004 | 1.000 | *2.177* | 1.017 | **0.946** | *1.080* |
| | 21.56 | 17632 | 0.995 | 1.000 | 0.993 | 1.000 | 0.998 | 1.000 | 1.004 | 1.000 | 1.001 | 1.000 | 0.993 | 1.000 | *3.660* | *1.598* | *2.274* | *2.333* |
| | 45.95 | 38975 | 0.997 | 1.000 | 0.994 | 1.000 | 0.996 | 1.000 | 1.004 | 1.000 | 1.005 | 1.000 | 1.005 | 1.000 | *1.340* | **0.588** | **0.812** | **0.806** |
| | 30.94 | 29476 | 0.988 | 1.000 | 1.001 | 1.000 | 0.998 | 1.000 | 0.990 | 1.000 | *1.082* | **0.937** | 0.994 | 1.000 | *1.266* | **0.557** | *1.084* | **0.937** |
| mitre | 10.20 | 1 | *1.390* | 1.000 | 0.996 | 1.000 | 1.012 | 1.000 | *1.054* | 1.000 | 0.991 | 1.000 | *1.071* | 1.000 | *1.381* | 1.000 | *1.374* | 1.000 |
| | 10.28 | 1 | *1.469* | 1.000 | 0.989 | 1.000 | 1.006 | 1.000 | 0.991 | 1.000 | *1.481* | 1.000 | 0.989 | 1.000 | *1.449* | 1.000 | *1.456* | 1.000 |
| | 10.94 | 1 | *1.384* | 1.000 | 0.955 | 1.000 | 0.961 | 1.000 | **0.946** | 1.000 | *1.357* | 1.000 | 0.955 | 1.000 | *1.389* | 1.000 | *1.447* | 1.000 |
| | 10.45 | 1 | *1.412* | 1.000 | 0.994 | 1.000 | 0.990 | 1.000 | 0.987 | 1.000 | 0.991 | 1.000 | 0.979 | 1.000 | *1.397* | 1.000 | *1.400* | 1.000 |
| | 13.12 | 1 | 1.019 | 1.000 | 0.996 | 1.000 | 0.984 | 1.000 | 0.982 | 1.000 | 1.008 | 1.000 | 0.994 | 1.000 | 1.046 | 1.000 | *1.055* | 1.000 |
| mod008 | 0.50 | 1 | 1.000 | 1.000 | 1.000 | 1.000 | 1.000 | 1.000 | 1.000 | 1.000 | 1.000 | 1.000 | 1.000 | 1.000 | 1.000 | 1.000 | 1.000 | 1.000 |
| | 0.54 | 2 | 0.994 | 1.000 | 1.006 | 1.000 | 1.000 | 1.000 | 0.994 | 1.000 | 0.994 | 1.000 | 1.000 | 1.000 | 1.026 | 1.000 | 1.032 | 1.000 |
| | 0.52 | 2 | 1.013 | 1.000 | 1.013 | 1.000 | 1.007 | 1.000 | 0.993 | 1.000 | 0.987 | 1.000 | 1.007 | 1.000 | 1.026 | 0.990 | 1.020 | 0.990 |
| | 0.52 | 1 | 0.993 | 1.000 | 1.000 | 1.000 | 1.000 | 1.000 | 1.000 | 1.000 | 1.000 | 1.000 | 1.000 | 1.000 | 1.020 | 1.000 | 1.033 | 1.000 |
| | 0.50 | 1 | 1.000 | 1.000 | 1.000 | 1.000 | 1.000 | 1.000 | 1.000 | 1.000 | 1.000 | 1.000 | 1.000 | 1.000 | 1.000 | 1.000 | 1.000 | 1.000 |

cont. on next page

Table 4: Detailed computational results on `MMM-IP` with five random seeds. Relative changes by at least 5% are highlighted in bold and blue (**improvement**) or italic and red (*deterioration*).

| | default | | degeneracy | | dualbound | | leaves | | obj | | nsols | | sblps | | nochecks | | onlyroot | |
|---|---|---|---|---|---|---|---|---|---|---|---|---|---|---|---|---|---|---|
| Instance | time | nodes | $time_Q$ | $nodes_Q$ | $time_Q$ | $nodes_Q$ | $time_Q$ | $nodes_Q$ | $time_Q$ | $nodes_Q$ | $time_Q$ | $nodes_Q$ | $time_Q$ | $nodes_Q$ | $time_Q$ | $nodes_Q$ | $time_Q$ | $nodes_Q$ |
| mod010 | 0.50 | 2 | 1.000 | 1.000 | 1.000 | 1.000 | 1.000 | 1.000 | 1.000 | 1.000 | 1.000 | 1.000 | 1.000 | 1.000 | *1.460* | 1.000 | *1.273* | 1.000 |
| | 0.50 | 1 | 1.000 | 1.000 | 1.000 | 1.000 | 1.000 | 1.000 | 1.000 | 1.000 | 1.000 | 1.000 | 1.000 | 1.000 | *1.193* | 1.000 | *1.180* | 1.000 |
| | 0.50 | 1 | 1.000 | 1.000 | 1.000 | 1.000 | 1.000 | 1.000 | 1.000 | 1.000 | 1.000 | 1.000 | 1.000 | 1.000 | *1.233* | 1.010 | *1.273* | 1.010 |
| | 0.52 | 2 | 1.007 | 1.000 | 1.020 | 1.000 | 0.993 | 1.000 | 0.993 | 1.000 | 1.007 | 1.000 | 1.007 | 1.000 | *1.329* | 1.000 | *1.197* | 1.000 |
| | 0.50 | 1 | 1.000 | 1.000 | 1.000 | 1.000 | 1.000 | 1.000 | 1.000 | 1.000 | 1.000 | 1.000 | 1.000 | 1.000 | *1.347* | 1.010 | *1.320* | 1.010 |
| mzzv11 | 332.43 | 2022 | 0.999 | 1.000 | 1.000 | 1.000 | 0.995 | 1.000 | 1.005 | 1.000 | 0.995 | 1.000 | 0.996 | 1.000 | *1.490* | **0.583** | **0.932** | **0.330** |
| | 229.41 | 959 | 1.002 | 1.000 | 0.999 | 1.000 | 1.001 | 1.000 | 1.005 | 1.000 | 1.002 | 1.000 | 1.003 | 1.000 | *1.453* | **0.818** | *1.573* | *1.220* |
| | 310.05 | 1421 | 1.004 | 1.000 | 1.000 | 1.000 | 1.001 | 1.000 | 0.995 | 1.000 | 0.999 | 1.000 | 0.999 | 1.000 | **0.910** | **0.415** | *1.175* | **0.525** |
| | 327.03 | 1579 | 0.987 | 1.000 | 1.001 | 1.000 | 0.995 | 1.000 | 0.995 | 1.000 | 0.994 | 1.000 | 0.993 | 1.000 | 1.000 | **0.378** | 0.964 | **0.392** |
| | 186.14 | 627 | 1.000 | 1.000 | 1.003 | 1.000 | 1.003 | 1.000 | 1.008 | 1.000 | 1.007 | 1.000 | 1.007 | 1.000 | *1.746* | **0.791** | *1.458* | **0.770** |
| mzzv42z | 152.40 | 265 | 0.992 | 1.000 | 0.996 | 1.000 | 0.999 | 1.000 | 0.993 | 1.000 | 0.991 | 1.000 | 0.999 | 1.000 | *2.086* | *1.477* | 0.984 | **0.753** |
| | 131.25 | 125 | 0.999 | 1.000 | 0.993 | 1.000 | 0.993 | 1.000 | 0.983 | 1.000 | 0.992 | 1.000 | *1.107* | *1.053* | *2.096* | *1.636* | *1.692* | *1.667* |
| | 197.62 | 200 | 1.006 | 1.000 | 1.003 | 1.000 | 1.006 | 1.000 | 1.005 | 1.000 | 1.002 | 1.000 | 1.016 | **0.923** | *2.262* | *1.597* | *1.436* | *1.117* |
| | 213.83 | 456 | 1.002 | 1.000 | 1.001 | 1.000 | 1.001 | 1.000 | 1.006 | 1.000 | 0.996 | 1.000 | **0.945** | **0.736** | *1.504* | **0.845** | *1.677* | 0.968 |
| | 167.65 | 300 | 0.999 | 1.000 | 0.998 | 1.000 | 0.998 | 1.000 | 1.000 | 1.000 | 1.007 | 1.000 | 1.003 | 1.000 | *1.347* | **0.657** | *1.143* | **0.695** |
| neos-1109824 | 37.52 | 1796 | 0.996 | 1.000 | 1.037 | 1.000 | 0.995 | 1.000 | 0.988 | 1.000 | *1.572* | *2.657* | 1.018 | 1.000 | *2.154* | *2.246* | *1.553* | *2.657* |
| | 46.67 | 3438 | 0.998 | 1.000 | 1.027 | *1.072* | 0.996 | 1.000 | 0.998 | 1.000 | 0.971 | 0.962 | **0.782** | **0.594** | *1.320* | **0.810** | 0.969 | 0.962 |
| | 35.16 | 1913 | 0.995 | 1.000 | 1.002 | 1.000 | 0.997 | 1.000 | 1.002 | 1.000 | 0.995 | 1.000 | *1.079* | 0.966 | *1.478* | 0.997 | *1.373* | *1.906* |
| | 77.65 | 8290 | 0.996 | 1.000 | 1.001 | 1.000 | 0.997 | 1.000 | 1.000 | 1.000 | **0.693** | **0.534** | 1.000 | 1.000 | 0.984 | **0.445** | **0.695** | **0.534** |
| | 39.64 | 2247 | 0.990 | 1.000 | 0.992 | 1.000 | 0.991 | 1.000 | 0.994 | 1.000 | *1.121* | *1.077* | 1.030 | **0.933** | *1.926* | *1.546* | *1.120* | *1.077* |
| neos-1337307 | 3600.00 | 299554 | 1.000 | 1.013 | 1.000 | 1.008 | 1.000 | 0.990 | 1.000 | 1.001 | 1.000 | *1.175* | 1.000 | 0.993 | 1.000 | **0.167** | 1.000 | *1.142* |
| | 3600.00 | 268930 | 1.000 | 0.999 | 1.000 | 1.001 | 1.000 | 0.998 | 1.000 | 0.999 | 1.000 | *1.085* | 1.000 | 0.999 | 1.000 | **0.169** | 1.000 | **0.837** |
| | 3600.00 | 350002 | 1.000 | 1.003 | 1.000 | 0.996 | 1.000 | 1.004 | 1.000 | 1.005 | 1.000 | 1.013 | 1.000 | 0.999 | 1.000 | **0.153** | 1.000 | *1.101* |
| | 3600.00 | 300182 | 1.000 | 1.003 | 1.000 | 0.998 | 1.000 | 0.999 | 1.000 | 1.003 | 1.000 | *1.136* | 1.000 | 0.998 | 1.000 | **0.152** | 1.000 | *1.124* |
| | 3600.00 | 367948 | 1.000 | 1.015 | 1.000 | 1.008 | 1.000 | 1.001 | 1.000 | 1.005 | 1.000 | **0.749** | 1.000 | 1.004 | 1.000 | **0.112** | 1.000 | 0.986 |
| neos-1601936 | 2034.60 | 1453 | 0.999 | 1.000 | *1.714* | *3.615* | 0.998 | 1.000 | 0.999 | 1.000 | *1.769* | **0.679** | 1.002 | 1.000 | *1.769* | *1.201* | *1.769* | **0.679** |
| | 3600.00 | 2861 | 1.000 | 1.004 | 1.000 | 1.039 | 1.000 | 1.004 | 1.000 | 0.997 | 1.000 | **0.393** | 1.000 | 1.002 | 1.000 | **0.356** | 1.000 | **0.609** |
| | 1348.55 | 2661 | 0.993 | 1.000 | *2.668* | *5.411* | 1.001 | 1.000 | 1.001 | 1.000 | *2.668* | **0.548** | 1.000 | 1.000 | *2.668* | **0.541** | *2.668* | **0.528** |
| | 3600.00 | 6917 | 1.000 | 1.000 | 1.000 | *1.168* | 1.000 | 1.000 | 1.000 | 1.004 | 1.000 | **0.429** | 1.000 | 1.000 | 1.000 | **0.831** | **0.898** | **0.158** |
| | 3600.00 | 6139 | 1.000 | 1.000 | 1.000 | **0.686** | 1.000 | 1.000 | 1.000 | 1.000 | **0.419** | **0.176** | 1.000 | 1.000 | 1.000 | **0.441** | **0.417** | **0.176** |
| neos-686190 | 83.42 | 4906 | 0.998 | 1.000 | 0.998 | 1.000 | 0.992 | 1.000 | 0.998 | 1.000 | 1.001 | 1.000 | 1.000 | 1.000 | *7.579* | *3.397* | *1.703* | *2.367* |
| | 86.19 | 6047 | 1.004 | 1.000 | 1.001 | 1.000 | 1.000 | 1.000 | 0.999 | 1.000 | 1.001 | 1.000 | 0.999 | 1.000 | *8.670* | *3.407* | *4.366* | *6.493* |
| | 114.74 | 9156 | 1.005 | 1.000 | 1.000 | 1.000 | 0.999 | 1.000 | 0.998 | 1.000 | 0.999 | 1.000 | 1.001 | 1.000 | *3.742* | **0.651** | **0.736** | **0.632** |
| | 116.53 | 9998 | 1.006 | 1.000 | 1.005 | 1.000 | 1.006 | 1.000 | 1.000 | 1.000 | 0.999 | 1.000 | 1.004 | 1.000 | *3.255* | **0.460** | **0.681** | **0.564** |
| | 121.04 | 9806 | 0.998 | 1.000 | 1.002 | 1.000 | 0.999 | 1.000 | 0.998 | 1.000 | 1.004 | 1.000 | 0.997 | 1.000 | *2.484* | **0.493** | **0.627** | **0.572** |
| neos-849702 | 114.99 | 1680 | *1.802* | **0.595** | **0.862** | **0.835** | **0.866** | **0.835** | *1.801* | **0.595** | *1.802* | **0.595** | *1.731* | *4.759* | *1.802* | **0.595** | *9.601* | *11.858* |
| | 353.12 | 77703 | *2.900* | **0.364** | *2.216* | *4.411* | *2.214* | *4.411* | *2.910* | **0.364** | *2.907* | **0.364** | *9.463* | *19.062* | *2.919* | **0.364** | *2.830* | *2.084* |
| | 357.57 | 9988 | *2.795* | **0.359** | **0.511** | **0.342** | **0.514** | **0.342** | *2.802* | **0.359** | *2.797* | **0.359** | *2.369* | *2.377* | *2.799* | **0.359** | *2.111* | **0.080** |
| | 261.49 | 16745 | 1.021 | **0.048** | *1.101* | 1.032 | *1.100* | 1.032 | 1.019 | **0.048** | 1.016 | **0.048** | *1.557* | *1.159* | 1.023 | **0.048** | *5.019* | *1.101* |

cont. on next page

Table 4: Detailed computational results on `MMM-IP` with five random seeds. Relative changes by at least 5% are highlighted in bold and blue (**improvement**) or italic and red (*deterioration*).

| | default | | degeneracy | | dualbound | | leaves | | obj | | nsols | | sblps | | nochecks | | onlyroot | |
|---|---|---|---|---|---|---|---|---|---|---|---|---|---|---|---|---|---|---|
| Instance | time | nodes | $time_Q$ | $nodes_Q$ | $time_Q$ | $nodes_Q$ | $time_Q$ | $nodes_Q$ | $time_Q$ | $nodes_Q$ | $time_Q$ | $nodes_Q$ | $time_Q$ | $nodes_Q$ | $time_Q$ | $nodes_Q$ | $time_Q$ | $nodes_Q$ |
| | 143.57 | 25014 | *1.981* | **0.335** | *2.290* | *4.955* | *2.283* | *4.955* | *1.981* | **0.335** | *1.985* | **0.335** | *6.512* | *16.346* | *1.980* | **0.335** | *5.285* | *5.390* |
| neos-934278 | 3600.00 | 1175 | 1.000 | 0.993 | 1.000 | 0.991 | 1.000 | 0.993 | 1.000 | 0.993 | 1.000 | 0.993 | 1.000 | 1.000 | 1.000 | **0.335** | 1.000 | **0.340** |
| | 3600.00 | 267 | 1.000 | 1.005 | 1.000 | 1.005 | 1.000 | 1.005 | 1.000 | 1.005 | 1.000 | 1.014 | 1.000 | 1.005 | 1.000 | *1.101* | 1.000 | *1.153* |
| | 3600.00 | 266 | 1.000 | 1.000 | 1.000 | 1.003 | 1.000 | 1.008 | 1.000 | 0.989 | 1.000 | 1.000 | 1.000 | 1.008 | 1.000 | *1.124* | 1.000 | *1.182* |
| | 3600.00 | 399 | 1.000 | 1.020 | 1.000 | 0.953 | 1.000 | 1.016 | 1.000 | 1.014 | 1.000 | 1.014 | 1.000 | 1.012 | 1.000 | **0.835** | 1.000 | **0.857** |
| | 3600.00 | 523 | 1.000 | 1.002 | 1.000 | 0.998 | 1.000 | 1.003 | 1.000 | 1.000 | 1.000 | 1.003 | 1.000 | 1.000 | 1.000 | **0.680** | 1.000 | **0.717** |
| neos18 | 24.17 | 1378 | **0.363** | **0.121** | 0.993 | 1.000 | 0.993 | 1.000 | 0.981 | 1.047 | 0.997 | 1.000 | **0.575** | **0.287** | **0.354** | **0.121** | **0.915** | **0.918** |
| | 21.26 | 1604 | *1.162* | **0.653** | *1.403* | *1.457* | 1.028 | *1.199* | 0.996 | 1.000 | 1.000 | 1.000 | *1.429* | *1.072* | 0.967 | **0.467** | 1.031 | *1.431* |
| | 23.36 | 1646 | *1.229* | **0.684** | 0.993 | 1.000 | 0.994 | 1.000 | 0.998 | 1.000 | 0.995 | 1.000 | **0.840** | **0.248** | *1.218* | **0.684** | **0.914** | **0.824** |
| | 19.21 | 996 | **0.812** | **0.443** | **0.840** | **0.849** | 0.997 | 1.000 | 0.997 | 1.000 | 1.021 | 1.000 | **0.789** | **0.366** | **0.810** | **0.443** | *1.059* | *1.091* |
| | 22.07 | 1049 | 0.987 | **0.925** | **0.855** | **0.634** | 1.004 | 1.000 | 1.001 | 1.000 | 0.997 | 1.000 | **0.906** | **0.657** | **0.948** | **0.795** | *1.174* | *2.047* |
| ns1208400 | 2439.01 | 780 | **0.385** | **0.219** | **0.396** | **0.178** | **0.146** | **0.062** | 0.994 | 1.000 | **0.485** | **0.271** | **0.296** | **0.090** | **0.485** | **0.271** | **0.482** | **0.404** |
| | 1985.63 | 5751 | *1.378* | *1.450* | **0.948** | **0.886** | *1.673* | **0.733** | 1.003 | 1.000 | *1.395* | *1.451* | **0.415** | **0.213** | *1.390* | *1.450* | **0.486** | **0.380** |
| | 1467.85 | 1351 | **0.363** | **0.355** | **0.810** | **0.504** | **0.415** | **0.283** | 1.000 | 1.000 | **0.368** | **0.355** | **0.713** | **0.156** | **0.367** | **0.355** | **0.673** | **0.768** |
| | 3600.00 | 2617 | **0.071** | **0.122** | **0.493** | *1.157* | **0.274** | **0.476** | **0.309** | **0.607** | **0.070** | **0.122** | **0.244** | **0.332** | **0.071** | **0.122** | **0.147** | **0.400** |
| | 2388.79 | 3984 | **0.551** | **0.313** | **0.915** | *1.198* | **0.649** | **0.464** | 0.999 | 1.000 | **0.054** | **0.069** | **0.434** | **0.229** | **0.051** | **0.068** | **0.533** | **0.656** |
| ns1688347 | 55.09 | 818 | **0.939** | 0.975 | 0.998 | 1.013 | *1.209* | *1.666* | *1.084* | *1.283* | 1.040 | **0.935** | 1.000 | 1.000 | 0.999 | *1.103* | 1.038 | **0.935** |
| | 45.64 | 422 | **0.931** | **0.518** | 0.996 | 1.000 | *1.062* | **0.892** | 0.986 | 0.969 | *1.391* | *1.591* | 0.998 | 1.000 | *1.058* | **0.839** | *1.396* | *1.591* |
| | 106.90 | 1690 | *1.509* | *1.626* | 1.000 | 1.000 | 0.999 | **0.911** | *1.152* | **0.855** | *1.478* | *1.240* | 1.001 | 1.000 | **0.884** | **0.571** | *1.492* | *1.240* |
| | 139.57 | 1219 | 0.957 | 1.024 | 1.018 | *1.095* | 0.994 | 1.000 | **0.945** | *1.088* | *1.130* | 0.999 | 0.998 | 1.000 | **0.682** | **0.329** | **0.771** | **0.500** |
| | 55.82 | 414 | **0.924** | **0.591** | 0.999 | 1.000 | **0.719** | **0.608** | 0.996 | 1.000 | *1.169* | *1.244* | 1.002 | 1.000 | **0.850** | **0.550** | *1.169* | *1.244* |
| ns1830653 | 132.81 | 5720 | 1.002 | 1.000 | 1.000 | 1.000 | **0.810** | **0.829** | 1.001 | 1.000 | *1.130* | 1.028 | 0.998 | 1.000 | *1.142* | **0.663** | *1.132* | 1.028 |
| | 123.49 | 5512 | 1.003 | 1.000 | 1.002 | 1.000 | *1.248* | *1.054* | 1.001 | 1.000 | *1.221* | *1.285* | 1.003 | 1.000 | *1.548* | *1.133* | *1.223* | *1.285* |
| | 161.13 | 10467 | 0.999 | 1.000 | 1.002 | 1.000 | 0.953 | *1.090* | 1.000 | 1.000 | *1.078* | 0.985 | **0.803** | **0.886** | *1.229* | **0.479** | *1.093* | 0.985 |
| | 143.05 | 7519 | 1.001 | 1.000 | **0.826** | **0.713** | *1.065* | 1.011 | 0.999 | 1.000 | *1.108* | *1.177* | **0.785** | **0.589** | *1.193* | **0.833** | *1.104* | *1.177* |
| | 146.20 | 5836 | 0.999 | 1.000 | 0.999 | 1.000 | **0.779** | **0.639** | 0.999 | 1.000 | *1.158* | **0.738** | **0.851** | **0.746** | 1.028 | **0.406** | *1.162* | **0.738** |
| nsrand-ipx | 422.05 | 50124 | 0.995 | 1.000 | 1.011 | 1.000 | 0.999 | 1.000 | 0.996 | 1.000 | 1.002 | 1.000 | *1.302* | **0.946** | *6.327* | **0.938** | *1.897* | *1.892* |
| | 643.65 | 86869 | 0.992 | 1.000 | 1.002 | 1.000 | 0.996 | 1.000 | 0.998 | 1.000 | 0.996 | 1.000 | *1.275* | **0.923** | *4.263* | **0.891** | **0.503** | **0.378** |
| | 568.28 | 56575 | 1.002 | 1.000 | 1.011 | 1.000 | 1.001 | 1.000 | 1.000 | 1.000 | 1.002 | 1.000 | *1.063* | 0.996 | *6.326* | **0.715** | *1.101* | 1.016 |
| | 582.50 | 61056 | 1.000 | 1.000 | 0.998 | 1.000 | 1.001 | 1.000 | 1.003 | 1.000 | 1.008 | 1.000 | *2.283* | **0.885** | *5.237* | *1.122* | **0.898** | **0.893** |
| | 367.03 | 38012 | 0.991 | 1.000 | 0.995 | 1.000 | 0.993 | 1.000 | 0.999 | 1.000 | 0.995 | 1.000 | *1.142* | 1.032 | *9.405* | *1.958* | *1.867* | *2.088* |
| opt1217 | 0.50 | 1 | 1.000 | 1.000 | 1.000 | 1.000 | 1.000 | 1.000 | 1.000 | 1.000 | 1.000 | 1.000 | 1.000 | 1.000 | 1.000 | 1.000 | 1.000 | 1.000 |
| | 0.50 | 1 | 1.000 | 1.000 | 1.000 | 1.000 | 1.000 | 1.000 | 1.000 | 1.000 | 1.000 | 1.000 | 1.000 | 1.000 | 1.000 | 1.000 | 1.000 | 1.000 |
| | 0.50 | 1 | 1.000 | 1.000 | 1.000 | 1.000 | 1.000 | 1.000 | 1.000 | 1.000 | 1.000 | 1.000 | 1.000 | 1.000 | 1.000 | 1.000 | 1.000 | 1.000 |
| | 0.50 | 1 | 1.000 | 1.000 | 1.000 | 1.000 | 1.000 | 1.000 | 1.000 | 1.000 | 1.000 | 1.000 | 1.000 | 1.000 | 1.000 | 1.000 | 1.000 | 1.000 |
| | 0.50 | 1 | 1.000 | 1.000 | 1.000 | 1.000 | 1.000 | 1.000 | 1.000 | 1.000 | 1.000 | 1.000 | 1.000 | 1.000 | 1.000 | 1.000 | 1.000 | 1.000 |
| p0033 | 0.50 | 1 | 1.000 | 1.000 | 1.000 | 1.000 | 1.000 | 1.000 | 1.000 | 1.000 | 1.000 | 1.000 | 1.000 | 1.000 | 1.000 | 1.000 | 1.000 | 1.000 |
| | 0.50 | 1 | 1.000 | 1.000 | 1.000 | 1.000 | 1.000 | 1.000 | 1.000 | 1.000 | 1.000 | 1.000 | 1.000 | 1.000 | 1.000 | 1.000 | 1.000 | 1.000 |
| | 0.50 | 1 | 1.000 | 1.000 | 1.000 | 1.000 | 1.000 | 1.000 | 1.000 | 1.000 | 1.000 | 1.000 | 1.000 | 1.000 | 1.000 | 1.000 | 1.000 | 1.000 |

cont. on next page

Table 4: Detailed computational results on `MMM-IP` with five random seeds. Relative changes by at least 5% are highlighted in bold and blue (**improvement**) or italic and red (*deterioration*).

| | `default` | | `degeneracy` | | `dualbound` | | `leaves` | | `obj` | | `nsols` | | `sblps` | | `nochecks` | | `onlyroot` | |
|---|---|---|---|---|---|---|---|---|---|---|---|---|---|---|---|---|---|---|
| Instance | **time** | **nodes** | **time**$_Q$ | **nodes**$_Q$ | **time**$_Q$ | **nodes**$_Q$ | **time**$_Q$ | **nodes**$_Q$ | **time**$_Q$ | **nodes**$_Q$ | **time**$_Q$ | **nodes**$_Q$ | **time**$_Q$ | **nodes**$_Q$ | **time**$_Q$ | **nodes**$_Q$ | **time**$_Q$ | **nodes**$_Q$ |
| | 0.50 | 1 | 1.000 | 1.000 | 1.000 | 1.000 | 1.000 | 1.000 | 1.000 | 1.000 | 1.000 | 1.000 | 1.000 | 1.000 | 1.000 | 1.000 | 1.000 | 1.000 |
| | 0.50 | 1 | 1.000 | 1.000 | 1.000 | 1.000 | 1.000 | 1.000 | 1.000 | 1.000 | 1.000 | 1.000 | 1.000 | 1.000 | 1.000 | 1.000 | 1.000 | 1.000 |
| p0201 | 0.55 | 5 | 0.994 | 1.000 | 1.006 | 1.000 | 0.987 | 1.000 | 1.000 | 1.000 | 0.981 | 1.000 | 1.000 | 1.000 | 0.968 | *1.057* | 0.968 | *1.057* |
| | 0.71 | 9 | 1.018 | 1.000 | 1.006 | 1.000 | 1.012 | 1.000 | 1.012 | 1.000 | 1.000 | 1.000 | 1.023 | 1.000 | **0.924** | 0.982 | **0.924** | 0.982 |
| | 1.59 | 5 | 1.019 | 1.000 | 1.008 | 1.000 | 1.000 | 1.000 | 1.000 | 1.000 | 1.015 | 1.000 | 1.023 | 1.000 | **0.749** | *1.210* | **0.819** | *2.419* |
| | 0.90 | 13 | 1.000 | 1.000 | 0.995 | 1.000 | 1.011 | 1.000 | 1.005 | 1.000 | 1.005 | 1.000 | 0.989 | 1.000 | *1.058* | *1.088* | **0.911** | *1.088* |
| | 0.91 | 59 | 0.990 | 1.000 | 0.990 | 1.000 | 1.000 | 1.000 | 0.990 | 1.000 | 0.995 | 1.000 | 0.995 | 1.000 | **0.890** | **0.761** | **0.812** | **0.761** |
| p0282 | 0.50 | 2 | 1.000 | 1.000 | 1.000 | 1.000 | 1.000 | 1.000 | 1.000 | 1.000 | 1.000 | 1.000 | 1.000 | 1.000 | 1.000 | 1.000 | 1.000 | 1.000 |
| | 0.50 | 2 | 1.000 | 1.000 | 1.000 | 1.000 | 1.000 | 1.000 | 1.000 | 1.000 | 1.000 | 1.000 | 1.000 | 1.000 | 1.000 | 1.000 | 1.000 | 1.000 |
| | 0.50 | 2 | 1.000 | 1.000 | 1.000 | 1.000 | 1.000 | 1.000 | 1.000 | 1.000 | 1.000 | 1.000 | 1.000 | 1.000 | 1.000 | 1.000 | 1.000 | 1.000 |
| | 0.50 | 3 | 1.000 | 1.000 | 1.000 | 1.000 | 1.000 | 1.000 | 1.007 | 1.000 | 1.000 | 1.000 | 1.013 | 1.000 | 1.000 | 0.990 | 1.033 | 1.000 |
| | 0.50 | 2 | 1.000 | 1.000 | 1.000 | 1.000 | 1.000 | 1.000 | 1.000 | 1.000 | 1.000 | 1.000 | 1.000 | 1.000 | 1.000 | 1.000 | 1.000 | 1.000 |
| p0548 | 0.50 | 1 | 1.000 | 1.000 | 1.000 | 1.000 | 1.000 | 1.000 | 1.000 | 1.000 | 1.000 | 1.000 | 1.000 | 1.000 | 1.000 | 1.000 | 1.000 | 1.000 |
| | 0.50 | 1 | 1.000 | 1.000 | 1.000 | 1.000 | 1.000 | 1.000 | 1.000 | 1.000 | 1.000 | 1.000 | 1.000 | 1.000 | 1.000 | 1.000 | 1.000 | 1.000 |
| | 0.50 | 1 | 1.000 | 1.000 | 1.000 | 1.000 | 1.000 | 1.000 | 1.000 | 1.000 | 1.000 | 1.000 | 1.000 | 1.000 | 1.000 | 1.000 | 1.000 | 1.000 |
| | 0.50 | 1 | 1.000 | 1.000 | 1.000 | 1.000 | 1.000 | 1.000 | 1.000 | 1.000 | 1.000 | 1.000 | 1.000 | 1.000 | 1.000 | 1.000 | 1.000 | 1.000 |
| | 0.50 | 1 | 1.000 | 1.000 | 1.000 | 1.000 | 1.000 | 1.000 | 1.000 | 1.000 | 1.000 | 1.000 | 1.000 | 1.000 | 1.000 | 1.000 | 1.000 | 1.000 |
| p2756 | 2.46 | 5 | 1.014 | 1.000 | 1.014 | 1.000 | 1.000 | 1.000 | 1.006 | 1.000 | 1.012 | 1.000 | 1.014 | 1.000 | **0.942** | *1.067* | **0.861** | 0.990 |
| | 3.33 | 5 | 1.000 | 1.000 | 1.005 | 1.000 | 1.005 | 1.000 | 1.005 | 1.000 | 1.002 | 1.000 | 1.012 | 1.000 | *1.102* | 1.010 | 0.998 | 0.990 |
| | 2.16 | 3 | 0.991 | 1.000 | 0.987 | 1.000 | 0.994 | 1.000 | 0.981 | 1.000 | 0.991 | 1.000 | 0.978 | 1.000 | *1.073* | 1.010 | 0.981 | 1.000 |
| | 2.60 | 3 | 0.983 | 1.000 | 0.994 | 1.000 | 0.994 | 1.000 | 0.983 | 1.000 | 0.983 | 1.000 | 0.978 | 1.000 | 0.969 | 1.010 | 0.964 | 1.019 |
| | 1.22 | 2 | 1.027 | 1.000 | 1.041 | 1.000 | 1.023 | 1.000 | 1.027 | 1.000 | 1.023 | 1.000 | 1.014 | 1.000 | *1.077* | 1.010 | *1.063* | 1.010 |
| protfold | 3600.00 | 4825 | 1.000 | **0.685** | 1.000 | 1.007 | 1.000 | 1.009 | 1.000 | 0.990 | 1.000 | 1.008 | 1.000 | 1.002 | 1.000 | **0.575** | 1.000 | **0.530** |
| | 3600.00 | 8262 | 1.000 | **0.315** | 1.000 | 1.001 | 1.000 | 1.011 | 1.000 | 0.994 | 1.000 | 0.997 | 1.000 | 1.001 | 1.000 | **0.303** | 1.000 | **0.289** |
| | 3600.00 | 5818 | 1.000 | **0.481** | 1.000 | 0.999 | 1.000 | 0.997 | 1.000 | 0.999 | 1.000 | 0.995 | 1.000 | 0.998 | 1.000 | **0.558** | 1.000 | **0.358** |
| | 3600.00 | 7269 | 1.000 | **0.262** | 1.000 | 1.002 | 1.000 | 1.001 | 1.000 | 1.001 | 1.000 | 1.006 | 1.000 | 1.001 | 1.000 | **0.264** | 1.000 | **0.231** |
| | 3600.00 | 4748 | 1.000 | **0.344** | 1.000 | 0.998 | 1.000 | **0.848** | 1.000 | 1.000 | 1.000 | 1.000 | 1.000 | 1.003 | 1.000 | **0.415** | 1.000 | **0.299** |
| pw-myciel4 | 3600.00 | 519342 | 1.000 | **0.695** | 1.000 | **0.895** | 1.000 | 0.999 | **0.662** | **0.587** | 1.000 | 1.000 | 1.000 | **0.858** | 1.000 | *1.087* | 1.000 | **0.868** |
| | 3600.00 | 442617 | 1.000 | **0.826** | **0.800** | **0.789** | 1.000 | 0.999 | **0.932** | 0.954 | 1.000 | 1.000 | 1.000 | **0.858** | 1.000 | **0.636** | 1.000 | 1.040 |
| | 3600.00 | 389396 | 1.000 | **0.826** | 1.000 | *1.146* | 1.000 | 0.994 | 1.000 | 0.974 | 1.000 | 0.992 | 1.000 | *1.053* | 1.000 | **0.780** | 1.000 | **0.910** |
| | 2920.77 | 395762 | *1.232* | **0.905** | *1.232* | **0.938** | 1.006 | 1.000 | *1.232* | **0.872** | 1.000 | 1.000 | *1.232* | **0.943** | *1.232* | *1.125* | *1.232* | 0.984 |
| | 3600.00 | 405310 | **0.848** | **0.934** | 1.000 | 0.991 | 1.000 | 1.000 | **0.939** | 1.044 | 1.000 | 1.001 | 1.000 | 1.031 | 0.971 | **0.834** | 1.000 | *1.204* |
| qnet1 | 4.52 | 33 | 1.004 | 1.000 | 1.005 | 1.000 | 1.013 | 1.000 | 0.998 | 1.000 | 1.005 | 1.000 | 1.007 | 1.000 | **0.455** | **0.797** | **0.455** | **0.797** |
| | 3.09 | 3 | 1.007 | 1.000 | 1.005 | 1.000 | 1.005 | 1.000 | 1.010 | 1.000 | 1.010 | 1.000 | 0.998 | 1.000 | *1.086* | 0.990 | **0.804** | 1.000 |
| | 3.94 | 29 | 0.996 | 1.000 | 1.002 | 1.000 | 0.994 | 1.000 | 1.004 | 1.000 | 1.000 | 1.000 | 1.002 | 1.000 | *1.138* | 1.000 | **0.949** | 1.000 |
| | 2.23 | 3 | 1.009 | 1.000 | 1.000 | 1.000 | 1.006 | 1.000 | 1.003 | 1.000 | 1.009 | 1.000 | 1.000 | 1.000 | *1.161* | *1.087* | *1.071* | 1.010 |
| | 3.69 | 34 | 1.004 | 1.000 | 0.989 | 1.000 | 0.998 | 1.000 | 0.994 | 1.000 | 0.989 | 1.000 | 1.000 | 1.000 | *1.198* | 0.970 | **0.887** | 1.022 |
| qnet1_o | 3.93 | 27 | 0.988 | 1.000 | 0.992 | 1.000 | 0.996 | 1.000 | 1.002 | 1.000 | 0.996 | 1.000 | 0.988 | 1.000 | **0.919** | **0.929** | **0.836** | **0.929** |
| | 3.35 | 24 | 1.007 | 1.000 | 1.005 | 1.000 | 1.000 | 1.000 | 1.002 | 1.000 | 0.998 | 1.000 | 1.009 | 1.000 | *1.094* | 0.992 | **0.897** | 0.992 |

cont. on next page

Table 4: Detailed computational results on `MMM-IP` with five random seeds. Relative changes by at least 5% are highlighted in bold and blue (**improvement**) or italic and red (*deterioration*).

| | default | | degeneracy | | dualbound | | leaves | | obj | | nsols | | sblps | | nochecks | | onlyroot | |
|---|---|---|---|---|---|---|---|---|---|---|---|---|---|---|---|---|---|---|
| Instance | **time** | **nodes** | **time$_Q$** | **nodes$_Q$** | **time$_Q$** | **nodes$_Q$** | **time$_Q$** | **nodes$_Q$** | **time$_Q$** | **nodes$_Q$** | **time$_Q$** | **nodes$_Q$** | **time$_Q$** | **nodes$_Q$** | **time$_Q$** | **nodes$_Q$** | **time$_Q$** | **nodes$_Q$** |
| | 3.37 | 15 | 1.005 | 1.000 | 1.009 | 1.000 | 0.995 | 1.000 | 1.011 | 1.000 | 1.014 | 1.000 | 1.002 | 1.000 | **0.712** | **0.904** | **0.616** | **0.887** |
| | 1.66 | 2 | 1.004 | 1.000 | 0.996 | 1.000 | 1.000 | 1.000 | 0.992 | 1.000 | 1.004 | 1.000 | 1.008 | 1.000 | *1.143* | 1.000 | 0.996 | 1.000 |
| | 1.83 | 3 | 1.004 | 1.000 | 1.007 | 1.000 | 1.007 | 1.000 | 1.000 | 1.000 | 1.007 | 1.000 | 1.004 | 1.000 | *1.085* | 1.000 | 1.018 | 1.000 |
| reblock67 | 222.59 | 49530 | 1.002 | 1.000 | 0.999 | 1.000 | 1.007 | 1.000 | 0.999 | 1.000 | 0.997 | 1.000 | 0.995 | 1.000 | *2.701* | *1.225* | **0.803** | **0.598** |
| | 161.89 | 40625 | 1.000 | 1.000 | 0.997 | 1.000 | 1.000 | 1.000 | 0.996 | 1.000 | 0.998 | 1.000 | 1.000 | 1.000 | *2.270* | **0.943** | *1.272* | *1.214* |
| | 194.41 | 42542 | 1.001 | 1.000 | 1.002 | 1.000 | **0.857** | **0.853** | 0.997 | 1.000 | 1.003 | 1.000 | 1.001 | 1.000 | *2.191* | 0.951 | **0.902** | **0.876** |
| | 199.06 | 52613 | 1.006 | 1.000 | 1.003 | 1.000 | 0.992 | 1.032 | 1.006 | 1.000 | 1.005 | 1.000 | 1.001 | 1.000 | *3.355* | *1.638* | **0.868** | **0.767** |
| | 239.84 | 57566 | 0.997 | 1.000 | 0.997 | 1.000 | 0.994 | 1.000 | 0.997 | 1.000 | 0.996 | 1.000 | 0.997 | 1.000 | *1.707* | *1.170* | **0.716** | **0.795** |
| rmine6 | 971.86 | 111047 | 1.000 | 1.000 | 1.006 | 1.000 | 1.002 | 1.000 | 1.000 | 1.000 | 1.001 | 1.000 | 0.999 | 1.000 | *1.724* | *1.068* | *1.441* | **0.906** |
| | 883.48 | 89213 | 0.999 | 1.000 | 1.000 | 1.000 | 0.990 | 1.000 | 0.999 | 1.000 | 1.002 | 1.000 | 1.002 | 1.000 | *2.249* | *1.574* | **0.945** | **0.921** |
| | 775.32 | 91530 | 1.002 | 1.000 | 1.002 | 1.000 | 1.000 | 1.000 | 1.001 | 1.000 | 1.002 | 1.000 | 1.000 | 1.000 | *2.090* | *1.268* | 1.037 | 1.014 |
| | 740.58 | 85429 | 1.002 | 1.000 | 1.000 | 1.000 | 0.999 | 1.000 | 1.005 | 1.000 | 0.999 | 1.000 | 0.992 | 1.000 | *1.817* | *1.249* | *1.310* | *1.072* |
| | 771.48 | 76017 | 1.002 | 1.000 | 1.001 | 1.000 | 1.002 | 1.000 | 1.001 | 1.000 | 0.998 | 1.000 | 1.001 | 1.000 | *1.581* | *1.168* | *1.134* | *1.144* |
| rococoC10-001000 | 601.83 | 46983 | 1.003 | 1.000 | 1.000 | 1.000 | 0.998 | 1.000 | 0.998 | 1.000 | 0.998 | 1.000 | 1.000 | 1.000 | *2.365* | *1.898* | *2.192* | *1.817* |
| | 833.40 | 75156 | 1.003 | 1.000 | 1.000 | 1.000 | 1.000 | 1.000 | 1.000 | 1.000 | 1.000 | 1.000 | 1.002 | 1.000 | *2.401* | *1.509* | *1.679* | *1.377* |
| | 517.62 | 37178 | 0.998 | 1.000 | 1.002 | 1.000 | 1.003 | 1.000 | *1.092* | *1.105* | 1.002 | 1.000 | 1.001 | 1.000 | *1.872* | *1.710* | *1.428* | *1.458* |
| | 892.04 | 51410 | 1.001 | 1.000 | 1.002 | 1.000 | 0.999 | 1.000 | 1.001 | 1.000 | 1.003 | 1.000 | 1.003 | 1.000 | *1.737* | *1.994* | *1.076* | *1.344* |
| | 1047.39 | 66016 | 1.001 | 1.000 | 1.000 | 1.000 | 1.001 | 1.000 | 1.001 | 1.000 | 1.003 | 1.000 | 1.003 | 1.000 | *2.255* | *1.863* | *1.367* | *1.434* |
| seymour | 3600.00 | 102303 | 1.000 | 0.996 | 1.000 | 0.993 | 1.000 | 0.997 | 1.000 | 0.999 | 1.000 | 0.997 | 1.000 | 1.001 | 1.000 | **0.717** | 1.000 | **0.768** |
| | 3600.00 | 120875 | 1.000 | 0.999 | 1.000 | 1.001 | 1.000 | 1.005 | 1.000 | 0.997 | 1.000 | 1.005 | 1.000 | 1.001 | 1.000 | **0.626** | 1.000 | **0.656** |
| | 3600.00 | 123371 | 1.000 | 1.006 | 1.000 | 1.006 | 1.000 | 1.008 | 1.000 | 1.010 | 1.000 | 1.010 | 1.000 | 1.006 | 1.000 | **0.737** | 1.000 | **0.791** |
| | 3600.00 | 107935 | 1.000 | 1.005 | 1.000 | 1.003 | 1.000 | 1.004 | 1.000 | 1.003 | 1.000 | 1.002 | 1.000 | 1.011 | 1.000 | **0.836** | 1.000 | **0.943** |
| | 3600.00 | 97288 | 1.000 | 1.003 | 1.000 | 1.001 | 1.000 | 1.006 | 1.000 | 1.001 | 1.000 | 1.003 | 1.000 | 1.002 | 1.000 | 1.002 | 1.000 | *1.069* |
| sp98ir | 57.90 | 3948 | 1.003 | 1.000 | 1.009 | 1.000 | 1.003 | 1.000 | 1.020 | 1.000 | 1.010 | 1.000 | 1.008 | 1.000 | *2.980* | *1.355* | *1.637* | *1.467* |
| | 65.87 | 4178 | 1.001 | 1.000 | 1.003 | 1.000 | 1.003 | 1.000 | 1.003 | 1.000 | 1.004 | 1.000 | 1.007 | 1.000 | *2.428* | *1.260* | *1.115* | *1.144* |
| | 67.64 | 4735 | 1.004 | 1.000 | 1.002 | 1.000 | 1.001 | 1.000 | 1.003 | 1.000 | 1.001 | 1.000 | 1.008 | 1.000 | *3.056* | *1.357* | *1.306* | *1.249* |
| | 72.68 | 4603 | 0.997 | 1.000 | 1.002 | 1.000 | 1.003 | 1.000 | 0.999 | 1.000 | 1.001 | 1.000 | 1.005 | 1.000 | *2.275* | **0.892** | **0.890** | **0.882** |
| | 70.05 | 4780 | 0.997 | 1.000 | 0.997 | 1.000 | 0.995 | 1.000 | 0.999 | 1.000 | 1.005 | 1.000 | 1.004 | 1.000 | *2.109* | *1.134* | *1.088* | *1.248* |
| stein27 | 1.12 | 2749 | 1.024 | 1.000 | 0.995 | 1.000 | 1.028 | 1.000 | 1.024 | 1.000 | 1.033 | 1.000 | 0.991 | 1.000 | *1.491* | *1.081* | *1.094* | *1.081* |
| | 1.19 | 2699 | 1.018 | 1.000 | 1.009 | 1.000 | 1.005 | 1.000 | 1.009 | 1.000 | 1.027 | 1.000 | 1.023 | 1.000 | *1.502* | 1.025 | 1.014 | 1.009 |
| | 1.25 | 2889 | 0.973 | 1.000 | 0.991 | 1.000 | 0.982 | 1.000 | 1.000 | 1.000 | 0.991 | 1.000 | 0.964 | 1.000 | *1.382* | 0.990 | 0.978 | 0.990 |
| | 1.14 | 2715 | 0.995 | *1.056* | 0.972 | 1.000 | 0.967 | 1.000 | 0.981 | 1.000 | 0.972 | 1.000 | 0.981 | 1.000 | *1.486* | *1.056* | 0.958 | *1.056* |
| | 1.13 | 2911 | 0.977 | 0.972 | 1.000 | 1.000 | 0.986 | 1.000 | 1.000 | 1.000 | 0.981 | 1.000 | 0.977 | 1.000 | *1.394* | 0.951 | 0.977 | 0.972 |
| stein45 | 12.80 | 40723 | 0.985 | 0.962 | 0.995 | 1.000 | 0.981 | 1.000 | 0.986 | 1.000 | 0.979 | 1.000 | 0.981 | 1.000 | *1.630* | 0.971 | 0.964 | **0.916** |
| | 12.94 | 39012 | 0.968 | 0.995 | 0.996 | 1.000 | 0.984 | 1.000 | 0.972 | 1.000 | 0.989 | 1.000 | 0.972 | 1.000 | *1.698* | 0.966 | **0.935** | 0.976 |
| | 13.04 | 41930 | 1.000 | 0.985 | 0.999 | 1.000 | 0.988 | 1.000 | 0.994 | 1.000 | 0.999 | 1.000 | 1.006 | 1.000 | *1.628* | **0.931** | 0.994 | 0.975 |
| | 12.20 | 37543 | 1.022 | 1.015 | 0.993 | 1.000 | 1.001 | 1.000 | 0.994 | 1.000 | 0.992 | 1.000 | 0.992 | 1.000 | *1.665* | *1.055* | 1.047 | 1.015 |
| | 13.29 | 40503 | 0.981 | 1.006 | 0.994 | 1.000 | 0.978 | 1.000 | 0.992 | 1.000 | 0.992 | 1.000 | 0.999 | 1.000 | *1.666* | 0.984 | 0.985 | 1.006 |
| tanglegram2 | 4.09 | 3 | 1.004 | 1.000 | 1.010 | 1.000 | 1.010 | 1.000 | 1.012 | 1.000 | 1.006 | 1.000 | 1.014 | 1.000 | 1.026 | 1.000 | 1.020 | 1.000 |

cont. on next page

Table 4: Detailed computational results on `MMM-IP` with five random seeds. Relative changes by at least 5% are highlighted in bold and blue (**improvement**) or italic and red (*deterioration*).

| | default | | degeneracy | | dualbound | | leaves | | obj | | nsols | | sblps | | nochecks | | onlyroot | |
|---|---|---|---|---|---|---|---|---|---|---|---|---|---|---|---|---|---|---|
| Instance | **time** | **nodes** | **time**$_{\mathbf{Q}}$ | **nodes**$_{\mathbf{Q}}$ | **time**$_{\mathbf{Q}}$ | **nodes**$_{\mathbf{Q}}$ | **time**$_{\mathbf{Q}}$ | **nodes**$_{\mathbf{Q}}$ | **time**$_{\mathbf{Q}}$ | **nodes**$_{\mathbf{Q}}$ | **time**$_{\mathbf{Q}}$ | **nodes**$_{\mathbf{Q}}$ | **time**$_{\mathbf{Q}}$ | **nodes**$_{\mathbf{Q}}$ | **time**$_{\mathbf{Q}}$ | **nodes**$_{\mathbf{Q}}$ | **time**$_{\mathbf{Q}}$ | **nodes**$_{\mathbf{Q}}$ |
| | 3.85 | 3 | 1.010 | 1.000 | 1.008 | 1.000 | 1.010 | 1.000 | 1.008 | 1.000 | 1.008 | 1.000 | 1.002 | 1.000 | 1.029 | 1.000 | 1.023 | 1.000 |
| | 4.25 | 3 | 0.990 | 1.000 | 1.000 | 1.000 | 0.998 | 1.000 | 1.008 | 1.000 | 1.008 | 1.000 | 1.008 | 1.000 | 1.011 | 1.000 | 1.023 | 1.000 |
| | 4.11 | 3 | 1.002 | 1.000 | 1.002 | 1.000 | 1.000 | 1.000 | 0.992 | 1.000 | 0.992 | 1.000 | 0.998 | 1.000 | 1.016 | 1.000 | 1.010 | 1.000 |
| | 4.23 | 3 | 0.994 | 1.000 | 1.000 | 1.000 | 0.985 | 1.000 | 0.990 | 1.000 | 0.998 | 1.000 | 0.994 | 1.000 | 1.010 | 1.000 | 1.011 | 1.000 |

Table 5: Detailed computational results on `MMM-IP` with five random seeds. Relative changes by at least 5% are highlighted in bold and blue (**improvement**) or italic and red (*deterioration*).

| | `default` | | `degeneracy` + `leaves` | | `leaves` + `obj` | | `all6criterion` | |
|---|---|---|---|---|---|---|---|---|
| Instance | **time** | **nodes** | **time$_Q$** | **nodes$_Q$** | **time$_Q$** | **nodes$_Q$** | **time$_Q$** | **nodes$_Q$** |
| 10teams | 2.14 | 1 | *1.293* | 1.000 | 1.022 | 1.000 | *1.248* | 1.000 |
| | 20.00 | 1367 | **0.244** | **0.072** | *1.230* | 1.001 | **0.246** | **0.072** |
| | 15.57 | 287 | *1.636* | *2.729* | 0.998 | 1.000 | *1.508* | *1.835* |
| | 4.16 | 5 | *2.723* | *1.219* | 1.006 | 1.000 | *2.725* | *1.219* |
| | 4.94 | 5 | **0.912** | 1.000 | 1.012 | 1.000 | **0.911** | 1.000 |
| 30n20b8 | 314.13 | 195 | **0.604** | **0.553** | 0.992 | 1.000 | **0.344** | **0.346** |
| | 157.33 | 71 | **0.711** | **0.596** | 1.000 | 1.000 | **0.808** | **0.596** |
| | 87.25 | 4 | 1.001 | 1.000 | 0.997 | 1.000 | *1.193* | 1.000 |
| | 213.06 | 102 | *1.202* | *1.312* | 0.976 | 0.950 | 1.048 | **0.782** |
| | 470.09 | 290 | **0.813** | **0.723** | 1.004 | 1.000 | **0.243** | **0.421** |
| acc-tight5 | 125.35 | 802 | *1.093* | *1.171* | **0.426** | **0.302** | *1.087* | *1.171* |
| | 84.38 | 508 | *1.786* | *2.650* | **0.426** | **0.252** | *1.790* | *2.650* |
| | 50.45 | 195 | 0.974 | **0.929** | *1.740* | *1.966* | 0.978 | **0.929** |
| | 85.88 | 545 | **0.339** | **0.189** | **0.420** | **0.274** | **0.343** | **0.189** |
| | 73.38 | 627 | **0.418** | **0.161** | *1.969* | *2.048* | **0.418** | **0.161** |
| air03 | 1.70 | 1 | 1.037 | 1.000 | 1.037 | 1.000 | 1.007 | 1.000 |
| | 1.71 | 1 | 1.004 | 1.000 | 1.007 | 1.000 | 1.004 | 1.000 |
| | 1.74 | 1 | 0.996 | 1.000 | 1.004 | 1.000 | 0.993 | 1.000 |
| | 1.75 | 1 | 1.000 | 1.000 | 0.993 | 1.000 | 1.002 | 1.000 |
| | 1.75 | 1 | 1.000 | 1.000 | 1.000 | 1.000 | 0.996 | 1.000 |
| air04 | 49.46 | 110 | 1.015 | 1.000 | 1.003 | 1.000 | 1.007 | 1.000 |
| | 36.09 | 39 | 1.005 | 1.000 | 0.999 | 1.000 | 1.003 | 1.000 |
| | 52.88 | 118 | 0.997 | 1.000 | 0.999 | 1.000 | 0.989 | 1.000 |
| | 36.56 | 31 | 1.010 | 1.000 | 1.015 | 1.000 | 1.001 | 1.000 |
| | 48.35 | 122 | 0.998 | 1.000 | 1.007 | 1.000 | 1.001 | 1.000 |
| air05 | 24.82 | 189 | 0.997 | 1.000 | 1.001 | 1.000 | 0.991 | 1.000 |
| | 26.54 | 239 | 1.004 | 1.000 | 1.001 | 1.000 | 1.009 | 1.000 |
| | 31.37 | 438 | 1.001 | 1.000 | 1.000 | 1.000 | 0.997 | 1.000 |
| | 27.36 | 283 | 1.005 | 1.000 | 1.004 | 1.000 | 1.010 | 1.000 |
| | 27.14 | 279 | 1.001 | 1.000 | 0.999 | 1.000 | 1.000 | 1.000 |
| blend2 | 0.69 | 714 | 1.036 | 1.000 | 1.047 | 1.000 | 1.036 | 1.000 |
| | 1.08 | 973 | 0.995 | 1.000 | 1.010 | 1.000 | 1.010 | 1.000 |
| | 1.05 | 912 | 1.000 | 1.000 | 0.985 | 1.000 | 0.990 | 1.000 |
| | 0.67 | 543 | 0.982 | 1.000 | 1.000 | 1.000 | 1.006 | 1.000 |
| | 0.66 | 709 | 1.042 | 1.000 | 1.048 | 1.000 | 1.042 | 1.000 |
| bley_xl1 | 164.49 | 1 | 1.009 | 1.000 | 1.000 | 1.000 | 1.007 | 1.000 |
| | 187.04 | 3 | 1.013 | *1.146* | 1.007 | 1.000 | 1.017 | *1.146* |
| | 182.51 | 13 | **0.887** | **0.903** | 1.006 | 1.000 | **0.896** | **0.903** |
| | 133.90 | 1 | 1.002 | 1.000 | 1.011 | 1.000 | 1.004 | 1.000 |
| | 135.56 | 1 | 1.012 | 1.000 | 1.008 | 1.000 | 1.007 | 1.000 |
| bnatt350 | 378.09 | 3383 | **0.052** | **0.030** | **0.107** | **0.054** | **0.052** | **0.030** |
| | 521.76 | 4947 | **0.033** | **0.022** | **0.292** | **0.147** | **0.033** | **0.022** |
| | 397.10 | 4029 | **0.047** | **0.026** | **0.123** | **0.053** | **0.048** | **0.026** |
| | 155.46 | 1495 | **0.109** | **0.067** | *2.608* | *2.974* | **0.109** | **0.067** |

cont. on next page

Table 5: Detailed computational results on `MMM-IP` with five random seeds. Relative changes by at least 5% are highlighted in bold and blue (**improvement**) or italic and red (*deterioration*).

| | `default` | | `degeneracy + leaves` | | `leaves + obj` | | `all6criterion` | |
|---|---|---|---|---|---|---|---|---|
| Instance | **time** | **nodes** | **time$_Q$** | **nodes$_Q$** | **time$_Q$** | **nodes$_Q$** | **time$_Q$** | **nodes$_Q$** |
| | 994.16 | 10057 | **0.044** | **0.021** | **0.330** | **0.259** | **0.044** | **0.021** |
| cap6000 | 2.69 | 1830 | 1.011 | 1.000 | 1.011 | 1.000 | 0.997 | 1.000 |
| | 2.67 | 1837 | 0.995 | 1.000 | 0.995 | 1.000 | 1.000 | 1.000 |
| | 2.61 | 1959 | 0.994 | 1.000 | 1.017 | 1.000 | 1.014 | 1.000 |
| | 3.17 | 1924 | 1.005 | 1.000 | 1.000 | 1.000 | 0.998 | 1.000 |
| | 2.82 | 1762 | 1.008 | 1.000 | 1.005 | 1.000 | 0.997 | 1.000 |
| cov1075 | 3600.00 | 1450527 | 0.991 | *1.189* | 1.000 | 0.996 | 0.991 | *1.189* |
| | 3600.00 | 1790225 | 1.000 | 0.974 | 1.000 | 1.002 | 1.000 | 0.978 |
| | 3600.00 | 1531529 | 1.000 | *1.179* | 1.000 | 1.000 | 1.000 | *1.151* |
| | 3600.00 | 1798178 | 1.000 | 0.984 | 1.000 | 1.008 | 1.000 | 0.984 |
| | 3600.00 | 1204894 | 1.000 | *1.507* | 1.000 | 0.998 | 1.000 | *1.505* |
| csched010 | 3600.00 | 237077 | 1.000 | 0.961 | 1.000 | 0.963 | **0.812** | **0.726** |
| | 3600.00 | 202876 | 1.000 | **0.863** | 1.000 | **0.861** | 1.000 | *1.145* |
| | 2891.69 | 177527 | *1.145* | *1.206* | *1.144* | *1.206* | *1.245* | *1.177* |
| | 3600.00 | 238611 | 0.989 | **0.820** | 0.986 | **0.820** | 1.000 | **0.927** |
| | 3168.37 | 173932 | *1.136* | *1.300* | *1.136* | *1.304* | *1.136* | *1.369* |
| eil33-2 | 69.44 | 583 | 1.013 | 1.000 | 1.007 | 1.000 | 1.006 | 1.000 |
| | 81.92 | 739 | 1.000 | 1.000 | 0.998 | 1.000 | 1.003 | 1.000 |
| | 84.87 | 847 | 1.004 | 1.000 | 1.005 | 1.000 | 1.001 | 1.000 |
| | 73.60 | 543 | 1.005 | 1.000 | 1.003 | 1.000 | 1.005 | 1.000 |
| | 70.01 | 583 | 1.001 | 1.000 | 1.001 | 1.000 | 0.999 | 1.000 |
| eilB101 | 378.99 | 10582 | 1.002 | 1.000 | 1.015 | 1.000 | 0.993 | 1.000 |
| | 365.27 | 8427 | 0.994 | 1.000 | 0.994 | 1.000 | 0.996 | 1.000 |
| | 401.63 | 9358 | 1.004 | 1.000 | 0.996 | 1.000 | 1.004 | 1.000 |
| | 397.46 | 11465 | 1.000 | 1.000 | 1.000 | 1.000 | 1.001 | 1.000 |
| | 365.79 | 10006 | 1.005 | 1.000 | 0.997 | 1.000 | 0.999 | 1.000 |
| enigma | 0.50 | 809 | 1.000 | **0.825** | 1.000 | 0.998 | 1.000 | **0.542** |
| | 0.50 | 449 | 1.000 | **0.497** | 1.000 | 1.000 | 1.000 | **0.470** |
| | 0.50 | 541 | 1.000 | **0.254** | 1.000 | 0.977 | 1.000 | **0.237** |
| | 0.50 | 679 | *1.060* | **0.537** | 1.000 | 1.000 | 1.000 | **0.151** |
| | 0.50 | 445 | 1.000 | **0.543** | 1.000 | 1.006 | 1.033 | **0.543** |
| enlight13 | 0.50 | 1 | 1.000 | 1.000 | 1.000 | 1.000 | 1.000 | 1.000 |
| | 0.50 | 1 | 1.000 | 1.000 | 1.000 | 1.000 | 1.000 | 1.000 |
| | 0.50 | 1 | 1.000 | 1.000 | 1.000 | 1.000 | 1.000 | 1.000 |
| | 0.50 | 1 | 1.000 | 1.000 | 1.000 | 1.000 | 1.000 | 1.000 |
| | 0.50 | 1 | 1.000 | 1.000 | 1.000 | 1.000 | 1.000 | 1.000 |
| enlight14 | 0.50 | 1 | 1.000 | 1.000 | 1.000 | 1.000 | 1.000 | 1.000 |
| | 0.50 | 1 | 1.000 | 1.000 | 1.000 | 1.000 | 1.000 | 1.000 |
| | 0.50 | 1 | 1.000 | 1.000 | 1.000 | 1.000 | 1.000 | 1.000 |
| | 0.50 | 1 | 1.000 | 1.000 | 1.000 | 1.000 | 1.000 | 1.000 |
| | 0.50 | 1 | 1.000 | 1.000 | 1.000 | 1.000 | 1.000 | 1.000 |
| fiber | 1.00 | 2 | 1.005 | 1.000 | 1.025 | 1.000 | 1.005 | 1.000 |
| | 1.42 | 4 | 1.000 | 1.000 | 0.988 | 1.000 | 1.000 | 1.000 |

cont. on next page

Table 5: Detailed computational results on `MMM-IP` with five random seeds. Relative changes by at least 5% are highlighted in bold and blue (**improvement**) or italic and red (*deterioration*).

| | `default` | | `degeneracy + leaves` | | `leaves + obj` | | `all6criterion` | |
|---|---|---|---|---|---|---|---|---|
| Instance | **time** | **nodes** | **time$_{Q}$** | **nodes$_{Q}$** | **time$_{Q}$** | **nodes$_{Q}$** | **time$_{Q}$** | **nodes$_{Q}$** |
| | 1.74 | 4 | 0.996 | 1.000 | 1.004 | 1.000 | 1.007 | 1.000 |
| | 1.53 | 4 | 1.000 | 1.000 | 1.008 | 1.000 | 1.000 | 1.000 |
| | 1.00 | 3 | 1.000 | 1.000 | 0.990 | 1.000 | 0.990 | 1.000 |
| flugpl | 0.50 | 1 | 1.000 | 1.000 | 1.000 | 1.000 | 1.000 | 1.000 |
| | 0.50 | 1 | 1.000 | 1.000 | 1.000 | 1.000 | 1.000 | 1.000 |
| | 0.50 | 1 | 1.000 | 1.000 | 1.000 | 1.000 | 1.000 | 1.000 |
| | 0.50 | 1 | 1.000 | 1.000 | 1.000 | 1.000 | 1.000 | 1.000 |
| | 0.50 | 1 | 1.000 | 1.000 | 1.000 | 1.000 | 1.000 | 1.000 |
| gt2 | 0.50 | 1 | 1.000 | 1.000 | 1.000 | 1.000 | 1.000 | 1.000 |
| | 0.50 | 1 | 1.000 | 1.000 | 1.000 | 1.000 | 1.000 | 1.000 |
| | 0.50 | 1 | 1.000 | 1.000 | 1.000 | 1.000 | 1.000 | 1.000 |
| | 0.50 | 1 | 1.000 | 1.000 | 1.000 | 1.000 | 1.000 | 1.000 |
| | 0.50 | 1 | 1.000 | 1.000 | 1.000 | 1.000 | 1.000 | 1.000 |
| harp2 | 691.56 | 1491878 | 1.008 | 1.000 | 0.999 | 1.000 | 1.005 | 1.000 |
| | 1159.66 | 2857677 | 0.998 | 1.000 | 0.998 | 1.000 | 0.997 | 1.000 |
| | 336.73 | 595854 | 0.999 | 1.000 | 1.006 | 1.000 | 0.999 | 1.000 |
| | 937.54 | 2046181 | 0.999 | 1.000 | 1.002 | 1.000 | 1.002 | 1.000 |
| | 425.68 | 923780 | 1.013 | 1.000 | 1.007 | 1.000 | 1.008 | 1.000 |
| iis-100-0-cov | 629.92 | 103302 | 0.998 | 1.000 | 1.021 | 1.000 | 0.992 | 1.000 |
| | 545.20 | 92393 | 0.995 | 1.000 | 0.997 | 1.000 | 1.001 | 1.000 |
| | 520.78 | 81639 | 1.001 | 1.000 | 1.000 | 1.000 | 1.000 | 1.000 |
| | 506.53 | 84800 | 0.999 | 1.000 | 0.995 | 1.000 | 1.000 | 1.000 |
| | 544.56 | 92313 | 1.007 | 1.000 | 1.003 | 1.000 | 1.004 | 1.000 |
| iis-bupa-cov | 1915.72 | 157619 | 1.000 | 1.000 | 1.003 | 1.000 | 1.001 | 1.000 |
| | 1944.83 | 164079 | 0.999 | 1.000 | 0.996 | 1.000 | 0.999 | 1.000 |
| | 1881.83 | 156888 | 0.991 | 1.000 | 0.992 | 1.000 | 1.008 | 1.000 |
| | 1959.19 | 161969 | 0.994 | 1.000 | 0.999 | 1.000 | 0.997 | 1.000 |
| | 1820.52 | 158284 | 0.994 | 1.000 | 0.993 | 1.000 | 1.014 | 1.000 |
| iis-pima-cov | 293.79 | 4824 | 1.004 | 1.000 | 0.999 | 1.000 | 0.999 | 1.000 |
| | 437.19 | 12189 | 1.007 | 1.000 | 1.001 | 1.000 | 1.013 | 1.000 |
| | 246.99 | 4668 | 1.002 | 1.000 | 1.000 | 1.000 | 1.002 | 1.000 |
| | 284.88 | 6607 | 1.002 | 1.000 | 1.000 | 1.000 | 1.003 | 1.000 |
| | 556.78 | 17821 | 0.997 | 1.000 | 0.997 | 1.000 | 1.001 | 1.000 |
| l152lav | 5.51 | 23 | 1.011 | 1.000 | 1.014 | 1.000 | 0.992 | 1.000 |
| | 1.67 | 31 | 0.989 | 1.000 | 0.993 | 1.000 | 1.000 | 1.000 |
| | 2.93 | 50 | 0.995 | 1.000 | 1.003 | 1.000 | 1.000 | 1.000 |
| | 2.34 | 20 | 0.985 | 1.000 | 0.994 | 1.000 | 0.991 | 1.000 |
| | 3.02 | 52 | 1.000 | 1.000 | 1.012 | 1.000 | 1.002 | 1.000 |
| lectsched-4-obj | 27.04 | 805 | **0.750** | **0.215** | *1.052* | **0.838** | **0.749** | **0.215** |
| | 23.16 | 796 | *1.193* | **0.288** | 1.008 | **0.823** | *1.198* | **0.288** |
| | 86.08 | 2419 | **0.285** | **0.206** | **0.626** | **0.461** | **0.265** | **0.121** |
| | 26.83 | 1231 | **0.327** | **0.082** | **0.933** | **0.540** | **0.331** | **0.082** |
| | 22.73 | 745 | **0.903** | **0.253** | 0.992 | 0.954 | **0.911** | **0.253** |

cont. on next page

Table 5: Detailed computational results on MMM-IP with five random seeds. Relative changes by at least 5% are highlighted in bold and blue (**improvement**) or italic and red (*deterioration*).

| | default | | degeneracy + leaves | | leaves + obj | | all6criterion | |
|---|---|---|---|---|---|---|---|---|
| Instance | **time** | **nodes** | **time$_Q$** | **nodes$_Q$** | **time$_Q$** | **nodes$_Q$** | **time$_Q$** | **nodes$_Q$** |
| lseu | 0.75 | 268 | 0.994 | 1.000 | 1.029 | 1.000 | 1.006 | 1.000 |
| | 0.50 | 4 | 1.000 | 1.000 | 1.000 | 1.000 | 1.000 | 1.000 |
| | 0.51 | 205 | 0.993 | 1.000 | 0.993 | 1.000 | 1.000 | 1.000 |
| | 0.50 | 91 | 1.000 | 1.000 | 1.000 | 1.000 | 1.000 | 1.000 |
| | 0.50 | 191 | 1.000 | 1.000 | 1.000 | 1.000 | 1.000 | 1.000 |
| m100n500k4r1 | 3600.00 | 2170459 | 1.000 | *1.118* | 1.000 | 1.002 | 1.000 | **0.781** |
| | 3600.00 | 2278212 | 1.000 | 1.035 | 1.000 | 0.996 | 1.000 | **0.783** |
| | 3600.00 | 2617105 | 1.000 | **0.884** | 1.000 | 1.000 | 1.000 | **0.560** |
| | 3600.00 | 2309660 | 1.000 | **0.935** | 1.000 | 0.998 | 1.000 | **0.638** |
| | 3600.00 | 2245563 | 1.000 | **0.834** | 1.000 | 0.996 | 1.000 | **0.647** |
| macrophage | 125.50 | 5163 | 1.000 | 1.000 | 1.006 | 1.000 | 1.001 | 1.000 |
| | 157.18 | 7241 | 0.998 | 1.000 | 0.998 | 1.000 | 0.999 | 1.000 |
| | 209.11 | 9131 | 1.000 | 1.000 | 1.000 | 1.000 | 1.001 | 1.000 |
| | 162.53 | 7927 | 0.998 | 1.000 | 1.000 | 1.000 | 1.002 | 1.000 |
| | 121.37 | 5239 | 0.995 | 1.000 | 0.996 | 1.000 | 0.997 | 1.000 |
| manna81 | 0.50 | 1 | 1.000 | 1.000 | 1.000 | 1.000 | 1.000 | 1.000 |
| | 0.50 | 1 | 1.000 | 1.000 | 1.000 | 1.000 | 1.000 | 1.000 |
| | 0.50 | 1 | 1.000 | 1.000 | 1.000 | 1.000 | 1.000 | 1.000 |
| | 0.50 | 1 | 1.000 | 1.000 | 1.000 | 1.000 | 1.000 | 1.000 |
| | 0.50 | 1 | 1.000 | 1.000 | 1.000 | 1.000 | 1.000 | 1.000 |
| markshare1 | 3600.00 | 19093205 | 1.000 | **0.927** | 1.000 | 0.997 | 1.000 | **0.924** |
| | 3600.00 | 20203360 | 1.000 | **0.859** | 1.000 | 0.998 | 1.000 | **0.856** |
| | 3600.00 | 19520584 | 1.000 | **0.919** | 1.000 | 0.999 | 1.000 | **0.918** |
| | 3600.00 | 20427204 | 1.000 | **0.901** | 1.000 | 1.000 | 1.000 | **0.903** |
| | 3600.00 | 21513467 | 1.000 | **0.795** | 1.000 | 1.000 | 1.000 | **0.796** |
| markshare2 | 3600.00 | 4072360 | 1.000 | **0.843** | 1.000 | 0.995 | 1.000 | **0.847** |
| | 3600.00 | 5640455 | 1.000 | **0.699** | 1.000 | 1.000 | 1.000 | **0.699** |
| | 3600.00 | 3582521 | 1.000 | 1.025 | 1.000 | 0.998 | 1.000 | 1.023 |
| | 3600.00 | 4239467 | 1.000 | *1.067* | 1.000 | 1.006 | 1.000 | *1.070* |
| | 3600.00 | 5142521 | 1.000 | **0.828** | 1.000 | 1.004 | 1.000 | **0.836** |
| mcsched | 285.92 | 13689 | 0.999 | 1.000 | 1.000 | 1.000 | 0.997 | 1.000 |
| | 271.88 | 11373 | 0.999 | 1.000 | 0.999 | 1.000 | 1.001 | 1.000 |
| | 233.97 | 10059 | 1.003 | 1.000 | 0.995 | 1.000 | 1.000 | 1.000 |
| | 286.25 | 12496 | 0.995 | 1.000 | 1.000 | 1.000 | 0.997 | 1.000 |
| | 270.69 | 10373 | 1.004 | 1.000 | 0.999 | 1.000 | 1.000 | 1.000 |
| mine-166-5 | 58.69 | 1569 | 0.995 | 1.000 | 1.001 | 1.000 | 0.998 | 1.000 |
| | 57.78 | 1191 | 0.993 | 1.000 | 0.994 | 1.000 | 0.992 | 1.000 |
| | 46.45 | 2244 | 0.997 | 1.000 | 1.004 | 1.000 | 1.005 | 1.000 |
| | 64.93 | 139 | 0.998 | 1.000 | 1.007 | 1.000 | 1.000 | 1.000 |
| | 42.28 | 2863 | 1.001 | 1.000 | 1.001 | 1.000 | 1.000 | 1.000 |
| mine-90-10 | 276.07 | 47643 | *1.162* | *1.090* | *1.162* | *1.090* | *1.157* | *1.090* |
| | 125.32 | 22195 | 0.999 | 1.000 | 0.997 | 1.000 | 1.002 | 1.000 |
| | 173.25 | 25979 | 1.001 | 1.000 | 0.999 | 1.000 | 1.000 | 1.000 |
| | 500.83 | 87496 | *1.322* | *1.061* | *1.321* | *1.061* | *1.322* | *1.061* |

cont. on next page

Table 5: Detailed computational results on MMM-IP with five random seeds. Relative changes by at least 5% are highlighted in bold and blue (**improvement**) or italic and red (*deterioration*).

| | default | | degeneracy + leaves | | leaves + obj | | all6criterion | |
|---|---|---|---|---|---|---|---|---|
| Instance | time | nodes | $time_Q$ | $nodes_Q$ | $time_Q$ | $nodes_Q$ | $time_Q$ | $nodes_Q$ |
| | 236.18 | 32520 | **0.907** | *1.242* | **0.908** | *1.242* | **0.904** | *1.242* |
| misc03 | 1.39 | 398 | 1.008 | 1.000 | 0.996 | 1.000 | 0.992 | 1.000 |
| | 0.71 | 15 | 0.994 | 1.000 | 1.018 | 1.000 | *1.094* | *1.226* |
| | 0.70 | 19 | 1.006 | 1.000 | 1.006 | 1.000 | *1.200* | *1.370* |
| | 1.00 | 111 | 1.000 | 1.000 | 0.990 | 1.000 | 0.995 | 1.000 |
| | 0.66 | 15 | 1.018 | 1.000 | 0.994 | 1.000 | *1.193* | *1.609* |
| misc07 | 38.05 | 28289 | 0.998 | 1.000 | 1.026 | 1.000 | *1.090* | *1.188* |
| | 36.66 | 29782 | 0.999 | 1.000 | 0.999 | 1.000 | 0.999 | 1.000 |
| | 21.56 | 17632 | 0.999 | 1.000 | 0.996 | 1.000 | 0.994 | 1.000 |
| | 45.95 | 38975 | 1.006 | 1.000 | 1.007 | 1.000 | 1.004 | 1.000 |
| | 30.94 | 29476 | 0.993 | 1.000 | 0.996 | 1.000 | *1.088* | **0.937** |
| mitre | 10.20 | 1 | *1.376* | 1.000 | 0.992 | 1.000 | *1.354* | 1.000 |
| | 10.28 | 1 | *1.444* | 1.000 | 1.009 | 1.000 | *1.473* | 1.000 |
| | 10.94 | 1 | *1.380* | 1.000 | 0.960 | 1.000 | *1.391* | 1.000 |
| | 10.45 | 1 | *1.401* | 1.000 | 0.994 | 1.000 | *1.399* | 1.000 |
| | 13.12 | 1 | 1.038 | 1.000 | 0.991 | 1.000 | *1.051* | 1.000 |
| mod008 | 0.50 | 1 | 1.000 | 1.000 | 1.000 | 1.000 | 1.000 | 1.000 |
| | 0.54 | 2 | 1.000 | 1.000 | 1.000 | 1.000 | 1.000 | 1.000 |
| | 0.52 | 2 | 1.000 | 1.000 | 0.993 | 1.000 | 1.000 | 1.000 |
| | 0.52 | 1 | 0.993 | 1.000 | 0.993 | 1.000 | 1.000 | 1.000 |
| | 0.50 | 1 | 1.000 | 1.000 | 1.000 | 1.000 | 1.000 | 1.000 |
| mod010 | 0.50 | 2 | 1.000 | 1.000 | 1.000 | 1.000 | 1.000 | 1.000 |
| | 0.50 | 1 | 1.000 | 1.000 | 1.000 | 1.000 | 1.000 | 1.000 |
| | 0.50 | 1 | 1.000 | 1.000 | 1.000 | 1.000 | 1.000 | 1.000 |
| | 0.52 | 2 | 1.000 | 1.000 | 0.993 | 1.000 | 0.993 | 1.000 |
| | 0.50 | 1 | 1.000 | 1.000 | 1.000 | 1.000 | 1.000 | 1.000 |
| mzzv11 | 332.43 | 2022 | 0.999 | 1.000 | 0.994 | 1.000 | 1.000 | 1.000 |
| | 229.41 | 959 | 1.002 | 1.000 | 1.003 | 1.000 | 1.001 | 1.000 |
| | 310.05 | 1421 | 0.998 | 1.000 | 1.000 | 1.000 | 0.998 | 1.000 |
| | 327.03 | 1579 | 0.994 | 1.000 | 0.995 | 1.000 | 0.994 | 1.000 |
| | 186.14 | 627 | 1.001 | 1.000 | 1.002 | 1.000 | 1.000 | 1.000 |
| mzzv42z | 152.40 | 265 | 0.994 | 1.000 | 0.992 | 1.000 | 0.993 | 1.000 |
| | 131.25 | 125 | 0.993 | 1.000 | 0.993 | 1.000 | *1.105* | *1.053* |
| | 197.62 | 200 | 1.004 | 1.000 | 1.006 | 1.000 | 1.024 | **0.923** |
| | 213.83 | 456 | 0.999 | 1.000 | 1.001 | 1.000 | **0.933** | **0.736** |
| | 167.65 | 300 | 1.000 | 1.000 | 1.000 | 1.000 | 1.004 | 1.000 |
| neos-1109824 | 37.52 | 1796 | 0.998 | 1.000 | 0.995 | 1.000 | *1.662* | *2.859* |
| | 46.67 | 3438 | 1.000 | 1.000 | 0.998 | 1.000 | 1.009 | 0.969 |
| | 35.16 | 1913 | 0.999 | 1.000 | 0.998 | 1.000 | *1.088* | 0.966 |
| | 77.65 | 8290 | 1.003 | 1.000 | 1.003 | 1.000 | **0.696** | **0.534** |
| | 39.64 | 2247 | 0.993 | 1.000 | 0.991 | 1.000 | *1.156* | *1.084* |
| neos-1337307 | 3600.00 | 302483 | 1.000 | 1.004 | 1.000 | 1.006 | 1.000 | *1.178* |
| | 3600.00 | 269510 | 1.000 | 1.001 | 1.000 | 1.003 | 1.000 | *1.084* |

cont. on next page

Table 5: Detailed computational results on MMM-IP with five random seeds. Relative changes by at least 5% are highlighted in bold and blue (**improvement**) or italic and red (*deterioration*).

| | default | | degeneracy + leaves | | leaves + obj | | all6criterion | |
|---|---|---|---|---|---|---|---|---|
| Instance | time | nodes | $time_Q$ | $nodes_Q$ | $time_Q$ | $nodes_Q$ | $time_Q$ | $nodes_Q$ |
| | 3600.00 | 348643 | 1.000 | 1.002 | 1.000 | 1.000 | 1.000 | 1.012 |
| | 3600.00 | 300587 | 1.000 | 1.011 | 1.000 | 0.993 | 1.000 | *1.131* |
| | 3600.00 | 367746 | 1.000 | 1.002 | 1.000 | 1.007 | 1.000 | **0.748** |
| neos-1601936 | 2034.60 | 1453 | 0.996 | 1.000 | 1.009 | 1.000 | *1.769* | *2.100* |
| | 3600.00 | 2849 | 1.000 | 1.003 | 1.000 | 1.000 | 1.000 | **0.377** |
| | 1348.55 | 2661 | 1.000 | 1.000 | 1.000 | 1.000 | *2.668* | **0.581** |
| | 3600.00 | 6919 | 1.000 | 0.999 | 1.000 | 0.999 | 1.000 | **0.584** |
| | 3600.00 | 6139 | 1.000 | 1.000 | 1.000 | 1.010 | 0.993 | **0.519** |
| neos-686190 | 83.42 | 4906 | 1.003 | 1.000 | 0.997 | 1.000 | 0.995 | 1.000 |
| | 86.19 | 6047 | 1.001 | 1.000 | 1.000 | 1.000 | 1.002 | 1.000 |
| | 114.74 | 9156 | 1.000 | 1.000 | 1.001 | 1.000 | 1.000 | 1.000 |
| | 116.53 | 9998 | 1.008 | 1.000 | 1.002 | 1.000 | 0.999 | 1.000 |
| | 121.04 | 9806 | 1.004 | 1.000 | 1.003 | 1.000 | 0.998 | 1.000 |
| neos-849702 | 114.99 | 2031 | *1.820* | **0.595** | *1.817* | **0.595** | *1.804* | **0.595** |
| | 353.12 | 17540 | *2.901* | **0.364** | *2.924* | **0.364** | *2.903* | **0.364** |
| | 357.57 | 29405 | *2.792* | **0.359** | *2.795* | **0.359** | *2.803* | **0.359** |
| | 261.49 | 16226 | 1.025 | **0.048** | 1.027 | **0.048** | 1.025 | **0.048** |
| | 143.57 | 4968 | *1.988* | **0.335** | *1.991* | **0.335** | *1.994* | **0.335** |
| neos-934278 | 3600.00 | 1184 | 1.000 | 0.991 | 1.000 | 0.993 | 1.000 | 0.993 |
| | 3600.00 | 265 | 1.000 | 1.014 | 1.000 | 1.011 | 1.000 | 1.003 |
| | 3600.00 | 263 | 1.000 | 1.000 | 1.000 | 1.000 | 1.000 | 0.989 |
| | 3600.00 | 391 | 1.000 | 1.014 | 1.000 | 1.014 | 1.000 | 0.963 |
| | 3600.00 | 521 | 1.000 | 1.002 | 1.000 | 0.997 | 1.000 | 0.998 |
| neos18 | 24.17 | 1378 | **0.360** | **0.121** | 1.004 | 1.047 | **0.359** | **0.121** |
| | 21.26 | 1321 | *1.154* | **0.653** | 1.027 | *1.199* | 0.971 | **0.467** |
| | 23.36 | 1646 | *1.227* | **0.684** | 1.000 | 1.000 | *1.229* | **0.684** |
| | 19.21 | 996 | **0.808** | **0.443** | 1.003 | 1.000 | **0.805** | **0.443** |
| | 22.07 | 1049 | 0.989 | **0.925** | 0.997 | 1.000 | 0.951 | **0.795** |
| ns1208400 | 2439.01 | 14168 | **0.485** | **0.271** | **0.146** | **0.062** | **0.484** | **0.271** |
| | 1985.63 | 7882 | *1.390* | *1.450* | *1.673* | **0.733** | *1.392* | *1.450* |
| | 1467.85 | 5026 | **0.367** | **0.355** | **0.416** | **0.283** | **0.367** | **0.355** |
| | 3600.00 | 5611 | **0.071** | **0.122** | **0.274** | **0.476** | **0.071** | **0.122** |
| | 2388.79 | 8711 | **0.050** | **0.068** | **0.650** | **0.464** | **0.051** | **0.068** |
| ns1688347 | 55.09 | 451 | **0.942** | 0.975 | *1.215* | *1.666* | **0.940** | 0.975 |
| | 45.64 | 485 | **0.930** | **0.518** | *1.058* | **0.892** | **0.929** | **0.518** |
| | 106.90 | 1865 | *1.530* | *1.438* | *1.160* | **0.638** | *1.330* | *1.366* |
| | 139.57 | 1219 | 0.956 | 1.024 | **0.942** | *1.088* | *1.205* | *1.089* |
| | 55.82 | 745 | **0.858** | **0.527** | **0.724** | **0.608** | **0.862** | **0.527** |
| ns1830653 | 132.81 | 6918 | **0.813** | **0.829** | **0.811** | **0.829** | 0.980 | **0.788** |
| | 123.49 | 5226 | *1.250* | *1.054* | *1.251* | *1.054* | *1.515* | *1.580* |
| | 161.13 | 9597 | **0.948** | *1.090* | **0.946** | *1.090* | **0.863** | **0.656** |
| | 143.05 | 7435 | *1.065* | 1.011 | *1.069* | 1.011 | **0.901** | **0.816** |
| | 146.20 | 9194 | **0.783** | **0.639** | **0.782** | **0.639** | *1.160* | **0.740** |

cont. on next page

Table 5: Detailed computational results on `MMM-IP` with five random seeds. Relative changes by at least 5% are highlighted in bold and blue (**improvement**) or italic and red (*deterioration*).

| Instance | default | | degeneracy + leaves | | leaves + obj | | all6criterion | |
|---|---|---|---|---|---|---|---|---|
| | **time** | **nodes** | $\textbf{time}_Q$ | $\textbf{nodes}_Q$ | $\textbf{time}_Q$ | $\textbf{nodes}_Q$ | $\textbf{time}_Q$ | $\textbf{nodes}_Q$ |
| nsrand-ipx | 422.05 | 50124 | 0.998 | 1.000 | 1.003 | 1.000 | *1.297* | **0.946** |
| | 643.65 | 86869 | 0.999 | 1.000 | 0.999 | 1.000 | *1.274* | **0.923** |
| | 568.28 | 56575 | 1.004 | 1.000 | 1.006 | 1.000 | *1.058* | 0.996 |
| | 582.50 | 61056 | 0.996 | 1.000 | 1.006 | 1.000 | *2.291* | **0.885** |
| | 367.03 | 38012 | 0.995 | 1.000 | 0.994 | 1.000 | *1.147* | 1.032 |
| opt1217 | 0.50 | 1 | 1.000 | 1.000 | 1.000 | 1.000 | 1.000 | 1.000 |
| | 0.50 | 1 | 1.000 | 1.000 | 1.000 | 1.000 | 1.000 | 1.000 |
| | 0.50 | 1 | 1.000 | 1.000 | 1.000 | 1.000 | 1.000 | 1.000 |
| | 0.50 | 1 | 1.000 | 1.000 | 1.000 | 1.000 | 1.000 | 1.000 |
| | 0.50 | 1 | 1.000 | 1.000 | 1.000 | 1.000 | 1.000 | 1.000 |
| p0033 | 0.50 | 1 | 1.000 | 1.000 | 1.000 | 1.000 | 1.000 | 1.000 |
| | 0.50 | 1 | 1.000 | 1.000 | 1.000 | 1.000 | 1.000 | 1.000 |
| | 0.50 | 1 | 1.000 | 1.000 | 1.000 | 1.000 | 1.000 | 1.000 |
| | 0.50 | 1 | 1.000 | 1.000 | 1.000 | 1.000 | 1.000 | 1.000 |
| | 0.50 | 1 | 1.000 | 1.000 | 1.000 | 1.000 | 1.000 | 1.000 |
| p0201 | 0.55 | 5 | 0.981 | 1.000 | 0.987 | 1.000 | 0.987 | 1.000 |
| | 0.71 | 9 | 1.018 | 1.000 | 1.023 | 1.000 | 0.994 | 1.000 |
| | 1.59 | 5 | 1.015 | 1.000 | 1.015 | 1.000 | 1.004 | 1.000 |
| | 0.90 | 13 | 1.005 | 1.000 | 1.000 | 1.000 | 1.005 | 1.000 |
| | 0.91 | 59 | 0.984 | 1.000 | 0.984 | 1.000 | 0.984 | 1.000 |
| p0282 | 0.50 | 2 | 1.000 | 1.000 | 1.000 | 1.000 | 1.000 | 1.000 |
| | 0.50 | 2 | 1.000 | 1.000 | 1.000 | 1.000 | 1.000 | 1.000 |
| | 0.50 | 2 | 1.000 | 1.000 | 1.000 | 1.000 | 1.000 | 1.000 |
| | 0.50 | 3 | 1.007 | 1.000 | 1.000 | 1.000 | 1.007 | 1.000 |
| | 0.50 | 2 | 1.000 | 1.000 | 1.000 | 1.000 | 1.000 | 1.000 |
| p0548 | 0.50 | 1 | 1.000 | 1.000 | 1.000 | 1.000 | 1.000 | 1.000 |
| | 0.50 | 1 | 1.000 | 1.000 | 1.000 | 1.000 | 1.000 | 1.000 |
| | 0.50 | 1 | 1.000 | 1.000 | 1.000 | 1.000 | 1.000 | 1.000 |
| | 0.50 | 1 | 1.000 | 1.000 | 1.000 | 1.000 | 1.000 | 1.000 |
| | 0.50 | 1 | 1.000 | 1.000 | 1.000 | 1.000 | 1.000 | 1.000 |
| p2756 | 2.46 | 5 | 1.006 | 1.000 | 1.012 | 1.000 | 1.000 | 1.000 |
| | 3.33 | 5 | 1.005 | 1.000 | 1.018 | 1.000 | 0.995 | 1.000 |
| | 2.16 | 3 | 0.994 | 1.000 | 0.987 | 1.000 | 0.994 | 1.000 |
| | 2.60 | 3 | 0.994 | 1.000 | 0.992 | 1.000 | 0.975 | 1.000 |
| | 1.22 | 2 | 1.036 | 1.000 | 1.027 | 1.000 | 1.041 | 1.000 |
| protfold | 3600.00 | 4783 | 1.000 | **0.685** | 1.000 | 1.008 | 1.000 | **0.685** |
| | 3600.00 | 8172 | 1.000 | **0.314** | 1.000 | 1.003 | 1.000 | **0.315** |
| | 3600.00 | 5838 | 1.000 | **0.478** | 1.000 | 1.004 | 1.000 | **0.480** |
| | 3600.00 | 7263 | 1.000 | **0.262** | 1.000 | 1.000 | 1.000 | **0.263** |
| | 3600.00 | 5620 | 1.000 | **0.348** | 1.000 | **0.844** | 1.000 | **0.345** |
| pw-myciel4 | 3600.00 | 519843 | 1.000 | **0.696** | **0.656** | **0.587** | 1.000 | **0.901** |
| | 3600.00 | 443116 | 1.000 | **0.824** | **0.930** | 0.954 | 1.000 | **0.688** |
| | 3600.00 | 391642 | 1.000 | **0.830** | 1.000 | 0.968 | 1.000 | 0.969 |
| | 2920.77 | 395762 | *1.232* | **0.907** | *1.232* | **0.872** | *1.232* | *1.067* |

cont. on next page

Table 5: Detailed computational results on `MMM-IP` with five random seeds. Relative changes by at least 5% are highlighted in bold and blue (**improvement**) or italic and red (*deterioration*).

| | `default` | | `degeneracy + leaves` | | `leaves + obj` | | `all6criterion` | |
|---|---|---|---|---|---|---|---|---|
| Instance | time | nodes | $time_Q$ | $nodes_Q$ | $time_Q$ | $nodes_Q$ | $time_Q$ | $nodes_Q$ |
| | 3600.00 | 405235 | **0.847** | **0.934** | **0.938** | 1.044 | 1.000 | **0.856** |
| qnet1 | 4.52 | 33 | 1.004 | 1.000 | 1.009 | 1.000 | 1.007 | 1.000 |
| | 3.09 | 3 | 1.000 | 1.000 | 1.002 | 1.000 | 1.000 | 1.000 |
| | 3.94 | 29 | 1.002 | 1.000 | 1.004 | 1.000 | 1.000 | 1.000 |
| | 2.23 | 3 | 0.994 | 1.000 | 0.994 | 1.000 | 1.006 | 1.000 |
| | 3.69 | 34 | 0.998 | 1.000 | 0.998 | 1.000 | 0.989 | 1.000 |
| qnet1_o | 3.93 | 27 | 1.008 | 1.000 | 0.986 | 1.000 | 1.004 | 1.000 |
| | 3.35 | 24 | 1.009 | 1.000 | 1.005 | 1.000 | 1.002 | 1.000 |
| | 3.37 | 15 | 1.000 | 1.000 | 1.007 | 1.000 | 1.000 | 1.000 |
| | 1.66 | 2 | 1.008 | 1.000 | 1.004 | 1.000 | 0.992 | 1.000 |
| | 1.83 | 3 | 1.004 | 1.000 | 1.000 | 1.000 | 1.004 | 1.000 |
| reblock67 | 222.59 | 49530 | 0.997 | 1.000 | 0.997 | 1.000 | 1.001 | 1.000 |
| | 161.89 | 40625 | 0.993 | 1.000 | 1.002 | 1.000 | 1.000 | 1.000 |
| | 194.41 | 49882 | **0.858** | **0.853** | **0.866** | **0.853** | **0.859** | **0.853** |
| | 199.06 | 50996 | 0.991 | 1.032 | 0.994 | 1.032 | 0.998 | 1.032 |
| | 239.84 | 57566 | 0.993 | 1.000 | 0.996 | 1.000 | 0.997 | 1.000 |
| rmine6 | 971.86 | 111047 | 1.004 | 1.000 | 1.000 | 1.000 | 0.999 | 1.000 |
| | 883.48 | 89213 | 1.001 | 1.000 | 0.990 | 1.000 | 1.001 | 1.000 |
| | 775.32 | 91530 | 1.005 | 1.000 | 0.999 | 1.000 | 1.011 | 1.000 |
| | 740.58 | 85429 | 1.003 | 1.000 | 0.993 | 1.000 | 1.006 | 1.000 |
| | 771.48 | 76017 | 0.999 | 1.000 | 1.000 | 1.000 | 1.004 | 1.000 |
| rococoC10-001000 | 601.83 | 46983 | 0.999 | 1.000 | 0.999 | 1.000 | 0.997 | 1.000 |
| | 833.40 | 75156 | 1.008 | 1.000 | 1.004 | 1.000 | 0.996 | 1.000 |
| | 517.62 | 37178 | 1.003 | 1.000 | *1.090* | *1.105* | *1.089* | *1.105* |
| | 892.04 | 51410 | 1.002 | 1.000 | 1.002 | 1.000 | 1.002 | 1.000 |
| | 1047.39 | 66016 | 1.000 | 1.000 | 1.001 | 1.000 | 1.001 | 1.000 |
| seymour | 3600.00 | 102645 | 1.000 | 0.996 | 1.000 | 0.999 | 1.000 | 0.997 |
| | 3600.00 | 120322 | 1.000 | 0.999 | 1.000 | 0.997 | 1.000 | 0.998 |
| | 3600.00 | 122337 | 1.000 | 1.006 | 1.000 | 1.004 | 1.000 | 1.007 |
| | 3600.00 | 107537 | 1.000 | 0.998 | 1.000 | 1.000 | 1.000 | 0.996 |
| | 3600.00 | 96726 | 1.000 | 1.007 | 1.000 | 1.001 | 1.000 | 0.999 |
| sp98ir | 57.90 | 3948 | 1.004 | 1.000 | 1.006 | 1.000 | 1.001 | 1.000 |
| | 65.87 | 4178 | 1.002 | 1.000 | 1.004 | 1.000 | 0.996 | 1.000 |
| | 67.64 | 4735 | 1.005 | 1.000 | 1.003 | 1.000 | 1.001 | 1.000 |
| | 72.68 | 4603 | 0.995 | 1.000 | 0.995 | 1.000 | 1.002 | 1.000 |
| | 70.05 | 4780 | 1.004 | 1.000 | 1.000 | 1.000 | 1.000 | 1.000 |
| stein27 | 1.12 | 2749 | 1.009 | 1.000 | 1.024 | 1.000 | 1.024 | 1.000 |
| | 1.19 | 2699 | 1.009 | 1.000 | 0.991 | 1.000 | 1.000 | 1.000 |
| | 1.25 | 2889 | 0.991 | 1.000 | 0.982 | 1.000 | 0.996 | 1.000 |
| | 1.14 | 2715 | 0.977 | *1.056* | 0.963 | 1.000 | 0.991 | *1.056* |
| | 1.13 | 2911 | 0.981 | 0.972 | 0.991 | 1.000 | 1.023 | 0.972 |
| stein45 | 12.80 | 40723 | 0.967 | 0.962 | 0.987 | 1.000 | 1.041 | 1.002 |
| | 12.94 | 39012 | 0.963 | 0.995 | 0.999 | 1.000 | 0.991 | 0.995 |

cont. on next page

Table 5: Detailed computational results on `MMM-IP` with five random seeds. Relative changes by at least 5% are highlighted in bold and blue (**improvement**) or italic and red (*deterioration*).

| Instance | default | | degeneracy + leaves | | leaves + obj | | all6criterion | |
|---|---|---|---|---|---|---|---|---|
| | **time** | **nodes** | **time$_Q$** | **nodes$_Q$** | **time$_Q$** | **nodes$_Q$** | **time$_Q$** | **nodes$_Q$** |
| | 13.04 | 41930 | 1.010 | 0.985 | 0.998 | 1.000 | 1.015 | 0.985 |
| | 12.20 | 37543 | 1.038 | 1.015 | 0.998 | 1.000 | 1.044 | 1.015 |
| | 13.29 | 40503 | 0.984 | 1.006 | 0.985 | 1.000 | 1.003 | 1.006 |
| tanglegram2 | 4.09 | 3 | 1.012 | 1.000 | 1.006 | 1.000 | 1.010 | 1.000 |
| | 3.85 | 3 | 1.000 | 1.000 | 1.006 | 1.000 | 0.996 | 1.000 |
| | 4.25 | 3 | 0.998 | 1.000 | 0.996 | 1.000 | 1.050 | 1.000 |
| | 4.11 | 3 | 0.994 | 1.000 | 0.990 | 1.000 | 1.004 | 1.000 |
| | 4.23 | 3 | 0.990 | 1.000 | 0.992 | 1.000 | 0.998 | 1.000 |

Table 6: Detailed computational results on `MMM-IP` with five random seeds. Relative changes by at least 5% are highlighted in bold and blue (**improvement**) or italic and red (*deterioration*).

| | `degeneracy` | | `variablebounds` | | `conflicts` | | `infervals` | | `primsols` | |
|---|---|---|---|---|---|---|---|---|---|---|
| Instance | **time** | **nodes** | **time$_Q$** | **nodes$_Q$** | **time$_Q$** | **nodes$_Q$** | **time$_Q$** | **nodes$_Q$** | **time$_Q$** | **nodes$_Q$** |
| 10teams | 3.01 | 1 | 0.990 | 1.000 | **0.825** | 1.000 | 0.995 | 1.000 | 0.988 | 1.000 |
| | 4.07 | 5 | 1.006 | 1.000 | *1.266* | 1.000 | *1.162* | 1.000 | 1.026 | 1.000 |
| | 25.62 | 956 | 1.012 | 1.000 | **0.668** | **0.159** | *1.091* | *1.195* | **0.441** | **0.118** |
| | 13.04 | 28 | 0.998 | 1.000 | **0.331** | **0.820** | **0.428** | **0.820** | 1.002 | 1.000 |
| | 4.40 | 5 | 1.004 | 1.000 | *4.026* | *2.943* | *1.278* | 1.000 | 1.007 | 1.000 |
| 30n20b8 | 189.53 | 63 | *1.521* | *1.245* | *1.573* | *1.607* | *1.458* | *1.503* | 1.000 | 1.000 |
| | 111.78 | 2 | *2.089* | *1.873* | *1.975* | *1.676* | 1.023 | 1.000 | 0.998 | 1.000 |
| | 87.35 | 4 | 1.001 | 1.000 | 0.993 | 1.000 | 0.999 | 1.000 | 1.001 | 1.000 |
| | 257.61 | 165 | **0.841** | **0.762** | 0.994 | 1.000 | 1.009 | 1.000 | 1.004 | 1.000 |
| | 381.26 | 182 | 0.982 | **0.926** | *1.061* | **0.858** | 1.005 | 1.000 | 0.999 | 1.000 |
| acc-tight5 | 139.72 | 956 | *1.284* | *1.562* | **0.205** | **0.116** | **0.433** | **0.267** | *1.560* | *1.788* |
| | 110.82 | 1160 | **0.516** | **0.318** | *1.515* | *1.136* | **0.513** | **0.239** | *1.140* | *1.150* |
| | 49.06 | 174 | *2.796* | *4.310* | *1.596* | *2.558* | **0.623** | **0.442** | *1.712* | *2.423* |
| | 28.52 | 22 | *2.037* | *1.828* | *2.899* | *4.975* | *1.400* | *1.434* | *5.993* | *11.180* |
| | 30.00 | 17 | 1.019 | 1.000 | **0.936** | *1.325* | *2.597* | *3.556* | 1.012 | *1.103* |
| air03 | 1.73 | 1 | 1.000 | 1.000 | 1.007 | 1.000 | 0.993 | 1.000 | 1.011 | 1.000 |
| | 1.72 | 1 | 1.015 | 1.000 | 1.000 | 1.000 | 1.011 | 1.000 | 0.996 | 1.000 |
| | 1.72 | 1 | 1.011 | 1.000 | 1.007 | 1.000 | 1.000 | 1.000 | 1.007 | 1.000 |
| | 1.75 | 1 | 1.000 | 1.000 | 0.993 | 1.000 | 0.993 | 1.000 | 1.004 | 1.000 |
| | 1.74 | 1 | 1.026 | 1.000 | 1.007 | 1.000 | 0.989 | 1.000 | 1.000 | 1.000 |
| air04 | 49.60 | 110 | 1.009 | 1.000 | 1.008 | 1.000 | 1.008 | 1.000 | 1.008 | 1.000 |
| | 36.38 | 39 | 1.001 | 1.000 | 0.995 | 1.000 | 0.997 | 1.000 | 0.995 | 1.000 |
| | 53.17 | 118 | 0.999 | 1.000 | 0.988 | 1.000 | 0.998 | 1.000 | 1.002 | 1.000 |
| | 37.09 | 31 | 0.995 | 1.000 | 0.993 | 1.000 | 0.992 | 1.000 | 0.996 | 1.000 |
| | 48.66 | 122 | 0.996 | 1.000 | 0.994 | 1.000 | 1.002 | 1.000 | 0.999 | 1.000 |

cont. on next page

Table 6: Detailed computational results on MMM-IP with five random seeds. Relative changes by at least 5% are highlighted in bold and blue (**improvement**) or italic and red (*deterioration*).

| | degeneracy | | variablebounds | | conflicts | | infervals | | primsols | |
|---|---|---|---|---|---|---|---|---|---|---|
| Instance | **time** | **nodes** | **time$_Q$** | **nodes$_Q$** | **time$_Q$** | **nodes$_Q$** | **time$_Q$** | **nodes$_Q$** | **time$_Q$** | **nodes$_Q$** |
| air05 | 24.56 | 189 | 1.002 | 1.000 | 1.002 | 1.000 | 1.001 | 1.000 | 1.002 | 1.000 |
| | 26.64 | 239 | 1.000 | 1.000 | 0.998 | 1.000 | 0.998 | 1.000 | 1.009 | 1.000 |
| | 31.10 | 438 | 1.002 | 1.000 | 1.005 | 1.000 | 1.046 | 1.000 | 1.003 | 1.000 |
| | 27.57 | 283 | 1.000 | 1.000 | 1.000 | 1.000 | 0.999 | 1.000 | 0.998 | 1.000 |
| | 27.23 | 279 | 1.000 | 1.000 | 0.996 | 1.000 | 0.998 | 1.000 | 0.993 | 1.000 |
| blend2 | 0.76 | 714 | 0.989 | 1.000 | 0.994 | 1.000 | 0.972 | 1.000 | 1.006 | 1.000 |
| | 1.05 | 973 | 1.015 | 1.000 | 0.980 | 1.000 | 1.005 | 1.000 | 1.005 | 1.000 |
| | 1.06 | 912 | 0.985 | 1.000 | 0.985 | 1.000 | 0.990 | 1.000 | 1.005 | 1.000 |
| | 0.66 | 543 | 1.018 | 1.000 | 0.994 | 1.000 | 0.994 | 1.000 | 1.012 | 1.000 |
| | 0.71 | 709 | 0.994 | 1.000 | 0.988 | 1.000 | 1.000 | 1.000 | 1.000 | 1.000 |
| bley_xl1 | 163.94 | 1 | 0.996 | 1.000 | 1.004 | 1.000 | 1.005 | 1.000 | 1.004 | 1.000 |
| | 189.31 | 18 | 0.999 | 1.000 | 1.002 | 1.000 | *1.052* | *1.314* | *1.107* | *1.458* |
| | 162.64 | 2 | 0.984 | 1.000 | 0.999 | 1.000 | 1.000 | 1.000 | 1.003 | 1.000 |
| | 134.67 | 1 | 0.997 | 1.000 | 1.001 | 1.000 | 0.998 | 1.000 | 0.999 | 1.000 |
| | 136.43 | 1 | 1.003 | 1.000 | 1.001 | 1.000 | 1.003 | 1.000 | 1.003 | 1.000 |
| bnatt350 | 18.90 | 6 | **0.897** | 1.000 | *6.188* | *7.047* | *2.522* | *1.330* | *19.530* | *40.189* |
| | 16.29 | 10 | *1.098* | 1.000 | *27.910* | *46.845* | *2.398* | *1.309* | *13.634* | *15.082* |
| | 18.02 | 6 | **0.834** | 1.000 | *15.298* | *20.085* | *3.298* | *3.075* | *13.075* | *20.151* |
| | 16.08 | 7 | 0.999 | 1.000 | *12.361* | *12.720* | *1.444* | 1.047 | *17.894* | *26.383* |
| | 43.21 | 112 | **0.732** | **0.623** | *9.481* | *21.113* | *2.919* | *3.425* | *10.242* | *27.854* |
| cap6000 | 2.70 | 1830 | 1.008 | 1.000 | 1.000 | 1.000 | 1.003 | 1.000 | 1.005 | 1.000 |
| | 2.62 | 1837 | 0.992 | 1.000 | 1.017 | 1.000 | 0.978 | 1.000 | 1.003 | 1.000 |
| | 2.58 | 1959 | 1.025 | 1.000 | 1.028 | 1.000 | 1.011 | 1.000 | 1.014 | 1.000 |
| | 3.20 | 1924 | 0.976 | 1.000 | 0.983 | 1.000 | 0.974 | 1.000 | 0.993 | 1.000 |
| | 2.85 | 1762 | 0.982 | 1.000 | 0.992 | 1.000 | 0.992 | 1.000 | 0.984 | 1.000 |

cont. on next page

Table 6: Detailed computational results on `MMM-IP` with five random seeds. Relative changes by at least 5% are highlighted in bold and blue (**improvement**) or italic and red (*deterioration*).

| | degeneracy | | variablebounds | | conflicts | | infervals | | primsols | |
|---|---|---|---|---|---|---|---|---|---|---|
| Instance | **time** | **nodes** | **time$_Q$** | **nodes$_Q$** | **time$_Q$** | **nodes$_Q$** | **time$_Q$** | **nodes$_Q$** | **time$_Q$** | **nodes$_Q$** |
| cov1075 | 3560.96 | 1725170 | 1.001 | 1.000 | 1.001 | 1.000 | 1.011 | **0.827** | 1.001 | 1.000 |
| | 3600.00 | 1745291 | 1.000 | 1.000 | 1.000 | 1.006 | 1.000 | 1.009 | 1.000 | 1.003 |
| | 3600.00 | 1834965 | 1.000 | 0.998 | 1.000 | 1.001 | 1.000 | **0.823** | 1.000 | 0.990 |
| | 3600.00 | 1771419 | 1.000 | 1.000 | 1.000 | 0.999 | 1.000 | 1.007 | 1.000 | 1.000 |
| | 3600.00 | 1817274 | 1.000 | 0.999 | 1.000 | 0.999 | 1.000 | **0.657** | 1.000 | 0.999 |
| csched010 | 3600.00 | 237694 | 1.000 | 0.997 | 1.000 | 0.997 | 1.000 | 0.996 | 1.000 | 0.997 |
| | 3600.00 | 203138 | 1.000 | 1.000 | 1.000 | 1.002 | 1.000 | 0.999 | 1.000 | 1.000 |
| | 2893.02 | 177527 | 0.998 | 1.000 | 1.001 | 1.000 | 1.000 | 1.000 | 0.999 | 1.000 |
| | 3600.00 | 238053 | 1.000 | 1.003 | 1.000 | 1.001 | 1.000 | 0.999 | 1.000 | 0.998 |
| | 3171.41 | 173932 | 0.999 | 1.000 | 0.999 | 1.000 | 0.998 | 1.000 | 1.000 | 1.000 |
| eil33-2 | 69.62 | 583 | 1.004 | 1.000 | 1.001 | 1.000 | 0.996 | 1.000 | 1.003 | 1.000 |
| | 82.52 | 739 | 0.994 | 1.000 | 0.994 | 1.000 | 0.999 | 1.000 | 0.981 | 1.000 |
| | 83.95 | 847 | 1.012 | 1.000 | 1.013 | 1.000 | 1.016 | 1.000 | 1.010 | 1.000 |
| | 72.71 | 543 | 1.015 | 1.000 | 1.019 | 1.000 | 1.016 | 1.000 | 1.005 | 1.000 |
| | 69.08 | 583 | 1.016 | 1.000 | 1.016 | 1.000 | 1.015 | 1.000 | 1.017 | 1.000 |
| eilB101 | 378.34 | 10582 | 0.995 | 1.000 | 0.990 | 1.000 | 1.012 | 1.000 | 1.000 | 1.000 |
| | 359.35 | 8427 | 1.009 | 1.000 | 1.001 | 1.000 | 1.017 | 1.000 | 1.011 | 1.000 |
| | 399.77 | 9358 | 1.003 | 1.000 | 0.999 | 1.000 | 0.999 | 1.000 | 1.005 | 1.000 |
| | 396.04 | 11465 | 1.006 | 1.000 | 0.995 | 1.000 | 1.012 | 1.000 | 0.997 | 1.000 |
| | 361.93 | 10006 | 1.009 | 1.000 | 1.011 | 1.000 | 1.013 | 1.000 | 1.006 | 1.000 |
| enigma | 0.50 | 650 | 1.000 | 0.997 | 1.033 | *768.000* | *1.113* | *1.135* | 1.000 | 1.033 |
| | 0.50 | 173 | 1.000 | *1.681* | 1.000 | *1.905* | 1.000 | *1.495* | 1.000 | 1.000 |
| | 0.50 | 63 | 1.000 | *1.374* | 1.000 | *1.123* | 1.027 | *4.491* | 1.000 | 1.000 |
| | 0.61 | 318 | **0.932** | **0.289** | **0.932** | *1.122* | 1.031 | *1.218* | 0.975 | 1.000 |
| | 0.53 | 196 | 1.013 | *1.588* | 0.980 | **0.405** | 0.980 | *1.493* | 1.000 | 1.000 |

cont. on next page

Table 6: Detailed computational results on `MMM-IP` with five random seeds. Relative changes by at least 5% are highlighted in bold and blue (**improvement**) or italic and red (*deterioration*).

| | degeneracy | | variablebounds | | conflicts | | infervals | | primsols | |
|---|---|---|---|---|---|---|---|---|---|---|
| Instance | **time** | **nodes** | **time$_Q$** | **nodes$_Q$** | **time$_Q$** | **nodes$_Q$** | **time$_Q$** | **nodes$_Q$** | **time$_Q$** | **nodes$_Q$** |
| enlight13 | 0.50 | 1 | 1.000 | 1.000 | 1.000 | 1.000 | 1.000 | 1.000 | 1.000 | 1.000 |
| | 0.50 | 1 | 1.000 | 1.000 | 1.000 | 1.000 | 1.000 | 1.000 | 1.000 | 1.000 |
| | 0.50 | 1 | 1.000 | 1.000 | 1.000 | 1.000 | 1.000 | 1.000 | 1.000 | 1.000 |
| | 0.50 | 1 | 1.000 | 1.000 | 1.000 | 1.000 | 1.000 | 1.000 | 1.000 | 1.000 |
| | 0.50 | 1 | 1.000 | 1.000 | 1.000 | 1.000 | 1.000 | 1.000 | 1.000 | 1.000 |
| enlight14 | 0.50 | 1 | 1.000 | 1.000 | 1.000 | 1.000 | 1.000 | 1.000 | 1.000 | 1.000 |
| | 0.50 | 1 | 1.000 | 1.000 | 1.000 | 1.000 | 1.000 | 1.000 | 1.000 | 1.000 |
| | 0.50 | 1 | 1.000 | 1.000 | 1.000 | 1.000 | 1.000 | 1.000 | 1.000 | 1.000 |
| | 0.50 | 1 | 1.000 | 1.000 | 1.000 | 1.000 | 1.000 | 1.000 | 1.000 | 1.000 |
| | 0.50 | 1 | 1.000 | 1.000 | 1.000 | 1.000 | 1.000 | 1.000 | 1.000 | 1.000 |
| fiber | 1.02 | 2 | 0.995 | 1.000 | 0.980 | 1.000 | 1.000 | 1.000 | 1.000 | 1.000 |
| | 1.39 | 4 | 1.004 | 1.000 | 1.008 | 1.000 | 1.013 | 1.000 | 1.004 | 1.000 |
| | 1.76 | 4 | 0.982 | 1.000 | 0.996 | 1.000 | 0.978 | 1.000 | 0.993 | 1.000 |
| | 1.55 | 4 | 1.000 | 1.000 | 0.988 | 1.000 | 0.988 | 1.000 | 0.988 | 1.000 |
| | 1.00 | 3 | 0.995 | 1.000 | 1.010 | 1.000 | 1.005 | 1.000 | 1.005 | 1.000 |
| flugpl | 0.50 | 1 | 1.000 | 1.000 | 1.000 | 1.000 | 1.000 | 1.000 | 1.000 | 1.000 |
| | 0.50 | 1 | 1.000 | 1.000 | 1.000 | 1.000 | 1.000 | 1.000 | 1.000 | 1.000 |
| | 0.50 | 1 | 1.000 | 1.000 | 1.000 | 1.000 | 1.000 | 1.000 | 1.000 | 1.000 |
| | 0.50 | 1 | 1.000 | 1.000 | 1.000 | 1.000 | 1.000 | 1.000 | 1.000 | 1.000 |
| | 0.50 | 1 | 1.000 | 1.000 | 1.000 | 1.000 | 1.000 | 1.000 | 1.000 | 1.000 |
| gt2 | 0.50 | 1 | 1.000 | 1.000 | 1.000 | 1.000 | 1.000 | 1.000 | 1.000 | 1.000 |
| | 0.50 | 1 | 1.000 | 1.000 | 1.000 | 1.000 | 1.000 | 1.000 | 1.000 | 1.000 |
| | 0.50 | 1 | 1.000 | 1.000 | 1.000 | 1.000 | 1.000 | 1.000 | 1.000 | 1.000 |
| | 0.50 | 1 | 1.000 | 1.000 | 1.000 | 1.000 | 1.000 | 1.000 | 1.000 | 1.000 |
| | 0.50 | 1 | 1.000 | 1.000 | 1.000 | 1.000 | 1.000 | 1.000 | 1.000 | 1.000 |

cont. on next page

Table 6: Detailed computational results on `MMM-IP` with five random seeds. Relative changes by at least 5% are highlighted in bold and blue (**improvement**) or italic and red (*deterioration*).

| | `degeneracy` | | `variablebounds` | | `conflicts` | | `infervals` | | `primsols` | |
|---|---|---|---|---|---|---|---|---|---|---|
| Instance | **time** | **nodes** | **time$_Q$** | **nodes$_Q$** | **time$_Q$** | **nodes$_Q$** | **time$_Q$** | **nodes$_Q$** | **time$_Q$** | **nodes$_Q$** |
| harp2 | 694.94 | 1491878 | 1.003 | 1.000 | 0.999 | 1.000 | 0.996 | 1.000 | 0.989 | 1.000 |
| | 1163.82 | 2857677 | 0.992 | 1.000 | 0.989 | 1.000 | 0.993 | 1.000 | 0.993 | 1.000 |
| | 333.99 | 595854 | 1.008 | 1.000 | 1.002 | 1.000 | 1.011 | 1.000 | 1.009 | 1.000 |
| | 934.82 | 2046181 | 1.006 | 1.000 | 0.999 | 1.000 | 1.001 | 1.000 | 1.003 | 1.000 |
| | 433.24 | 923780 | 0.998 | 1.000 | 0.985 | 1.000 | 0.991 | 1.000 | 0.996 | 1.000 |
| iis-100-0-cov | 625.29 | 103302 | 0.999 | 1.000 | 1.002 | 1.000 | 0.998 | 1.000 | 0.998 | 1.000 |
| | 544.56 | 92393 | 1.001 | 1.000 | 1.025 | 1.000 | 0.999 | 1.000 | 1.001 | 1.000 |
| | 520.20 | 81639 | 1.022 | 1.000 | 1.018 | 1.000 | *1.069* | 1.000 | 1.007 | 1.000 |
| | 507.53 | 84800 | 0.995 | 1.000 | 1.008 | 1.000 | 1.001 | 1.000 | 1.025 | 1.000 |
| | 542.59 | 92313 | *1.055* | 1.000 | 1.036 | 1.000 | 1.011 | 1.000 | 1.016 | 1.000 |
| iis-bupa-cov | 1917.51 | 157619 | 0.991 | 1.000 | 0.994 | 1.000 | 0.997 | 1.000 | 0.996 | 1.000 |
| | 1939.94 | 164079 | 1.010 | 1.000 | 1.019 | 1.000 | 1.014 | 1.000 | 1.020 | 1.000 |
| | 1878.21 | 156888 | *1.053* | 1.000 | 0.998 | 1.000 | 1.024 | 1.000 | 1.003 | 1.000 |
| | 1938.34 | 161969 | 1.002 | 1.000 | 1.022 | 1.000 | 1.026 | 1.000 | 1.002 | 1.000 |
| | 1802.15 | 158284 | 1.030 | 1.000 | 1.048 | 1.000 | 1.003 | 1.000 | 1.030 | 1.000 |
| iis-pima-cov | 292.21 | 4824 | 1.006 | 1.000 | 1.002 | 1.000 | 1.002 | 1.000 | 1.006 | 1.000 |
| | 434.93 | 12189 | 1.020 | 1.000 | 1.000 | 1.000 | 1.025 | 1.000 | 1.039 | 1.000 |
| | 246.99 | 4668 | 1.007 | 1.000 | 1.019 | 1.000 | 1.009 | 1.000 | 0.995 | 1.000 |
| | 285.55 | 6607 | 1.005 | 1.000 | 1.006 | 1.000 | 1.011 | 1.000 | 1.015 | 1.000 |
| | 558.43 | 17821 | 1.017 | 1.000 | 0.998 | 1.000 | 1.002 | 1.000 | 1.048 | 1.000 |
| l152lav | 5.47 | 23 | 0.994 | 1.000 | 0.998 | 1.000 | 1.009 | 1.000 | 1.000 | 1.000 |
| | 1.63 | 31 | 1.011 | 1.000 | 1.015 | 1.000 | 0.996 | 1.000 | 1.011 | 1.000 |
| | 2.93 | 50 | 0.992 | 1.000 | 0.987 | 1.000 | 1.000 | 1.000 | 1.003 | 1.000 |
| | 2.29 | 20 | 1.006 | 1.000 | 1.012 | 1.000 | 1.009 | 1.000 | 1.009 | 1.000 |
| | 3.00 | 52 | 1.012 | 1.000 | 1.010 | 1.000 | 1.015 | 1.000 | 1.012 | 1.000 |

cont. on next page

Table 6: Detailed computational results on `MMM-IP` with five random seeds. Relative changes by at least 5% are highlighted in bold and blue (**improvement**) or italic and red (*deterioration*).

| Instance | degeneracy time | degeneracy nodes | variablebounds $time_Q$ | variablebounds $nodes_Q$ | conflicts $time_Q$ | conflicts $nodes_Q$ | infervals $time_Q$ | infervals $nodes_Q$ | primsols $time_Q$ | primsols $nodes_Q$ |
|---|---|---|---|---|---|---|---|---|---|---|
| lectsched-4-obj | 20.11 | 95 | 0.957 | *1.113* | *1.624* | *3.487* | 1.008 | *1.913* | *2.428* | *4.513* |
| | 28.02 | 158 | *1.120* | *1.093* | *1.242* | *2.430* | 0.989 | **0.930** | *2.457* | *2.988* |
| | 24.11 | 565 | **0.645** | **0.415** | *1.559* | *1.989* | *1.552* | *1.788* | *2.656* | *1.714* |
| | 8.09 | 9 | 1.022 | 1.000 | **0.865** | 1.000 | **0.818** | 0.963 | *9.364* | *11.716* |
| | 20.57 | 114 | *1.340* | *1.178* | *1.510* | *3.664* | *1.381* | *2.271* | *7.981* | *9.061* |
| lseu | 0.77 | 268 | 1.011 | 1.000 | 0.977 | 1.000 | 1.000 | 1.000 | 0.994 | 1.000 |
| | 0.50 | 4 | 1.000 | 1.000 | 1.000 | 1.000 | 1.000 | 1.000 | 1.000 | 1.000 |
| | 0.50 | 205 | 1.000 | 1.000 | 1.020 | 1.000 | 1.000 | 1.000 | 1.000 | 1.000 |
| | 0.50 | 91 | 1.000 | 1.000 | 1.000 | 1.000 | 1.000 | 1.000 | 1.000 | 1.000 |
| | 0.51 | 191 | 0.993 | 1.000 | 1.007 | 1.000 | 0.993 | 1.000 | 1.020 | 1.000 |
| m100n500k4r1 | 3600.00 | 2426435 | 1.000 | 1.032 | 1.000 | 0.997 | 1.000 | 0.957 | 1.000 | **0.883** |
| | 3600.00 | 2358830 | 1.000 | 0.997 | 1.000 | 0.993 | 1.000 | *1.414* | 1.000 | **0.897** |
| | 3600.00 | 2313256 | 1.000 | 1.028 | 1.000 | 0.982 | 1.000 | 1.039 | 1.000 | 0.982 |
| | 3600.00 | 2163125 | 1.000 | 1.040 | 1.000 | 0.997 | 1.000 | *1.441* | 1.000 | **0.950** |
| | 3600.00 | 1871050 | 1.000 | *1.054* | 1.000 | 0.975 | 1.000 | *1.115* | 1.000 | 1.032 |
| macrophage | 125.68 | 5163 | 1.001 | 1.000 | 0.999 | 1.000 | 1.002 | 1.000 | 1.001 | 1.000 |
| | 157.07 | 7241 | 1.007 | 1.000 | 1.007 | 1.000 | 1.007 | 1.000 | 1.005 | 1.000 |
| | 209.04 | 9131 | 1.010 | 1.000 | 0.998 | 1.000 | 1.003 | 1.000 | 1.006 | 1.000 |
| | 162.46 | 7927 | 1.039 | 1.000 | 1.002 | 1.000 | 1.013 | 1.000 | 1.008 | 1.000 |
| | 121.23 | 5239 | 0.997 | 1.000 | 1.015 | 1.000 | 0.998 | 1.000 | 1.014 | 1.000 |
| manna81 | 0.50 | 1 | 1.000 | 1.000 | 1.000 | 1.000 | 1.000 | 1.000 | 1.000 | 1.000 |
| | 0.50 | 1 | 1.000 | 1.000 | 1.000 | 1.000 | 1.000 | 1.000 | 1.000 | 1.000 |
| | 0.50 | 1 | 1.000 | 1.000 | 1.000 | 1.000 | 1.000 | 1.000 | 1.000 | 1.000 |
| | 0.50 | 1 | 1.000 | 1.000 | 1.000 | 1.000 | 1.000 | 1.000 | 1.000 | 1.000 |
| | 0.50 | 1 | 1.000 | 1.000 | 1.000 | 1.000 | 1.000 | 1.000 | 1.000 | 1.000 |

cont. on next page

Table 6: Detailed computational results on `MMM-IP` with five random seeds. Relative changes by at least 5% are highlighted in bold and blue (**improvement**) or italic and red (*deterioration*).

| | degeneracy | | variablebounds | | conflicts | | infervals | | primsols | |
|---|---|---|---|---|---|---|---|---|---|---|
| Instance | **time** | **nodes** | **time**$_Q$ | **nodes**$_Q$ | **time**$_Q$ | **nodes**$_Q$ | **time**$_Q$ | **nodes**$_Q$ | **time**$_Q$ | **nodes**$_Q$ |
| markshare1 | 3600.00 | 17773362 | 1.000 | *1.150* | 1.000 | 0.990 | 1.000 | 1.038 | 1.000 | **0.655** |
| | 3600.00 | 17405834 | 1.000 | *1.184* | 1.000 | 0.998 | 1.000 | *1.126* | 1.000 | 1.000 |
| | 3600.00 | 17945211 | 1.000 | *1.149* | 1.000 | 1.001 | 1.000 | 0.954 | 1.000 | **0.482** |
| | 3600.00 | 18348094 | 1.000 | *1.131* | 1.000 | 1.003 | 1.000 | 0.980 | 1.000 | **0.853** |
| | 3600.00 | 17151823 | 1.000 | *1.130* | 1.000 | 0.997 | 1.000 | *1.050* | 1.000 | **0.889** |
| markshare2 | 3600.00 | 3433363 | 1.000 | *1.178* | 1.000 | 1.003 | 1.000 | *1.352* | 1.000 | *1.064* |
| | 3600.00 | 3929356 | 1.000 | **0.901** | 1.000 | 1.003 | 1.000 | *1.127* | 1.000 | **0.877** |
| | 3600.00 | 3675025 | 1.000 | 1.036 | 1.000 | 1.000 | 1.000 | 0.975 | 1.000 | 1.004 |
| | 3600.00 | 4532415 | 1.000 | **0.819** | 1.000 | 0.999 | 1.000 | 1.032 | 1.000 | **0.789** |
| | 3600.00 | 4267105 | 1.000 | **0.860** | 1.000 | 1.009 | 1.000 | **0.887** | 1.000 | *1.101* |
| mcsched | 285.70 | 13689 | 1.000 | 1.000 | 0.998 | 1.000 | 1.004 | 1.000 | 1.000 | 1.000 |
| | 271.59 | 11373 | 0.998 | 1.000 | 1.015 | 1.000 | 1.022 | 1.000 | 0.999 | 1.000 |
| | 233.00 | 10059 | 1.014 | 1.000 | 1.036 | 1.000 | 1.006 | 1.000 | 1.001 | 1.000 |
| | 285.38 | 12496 | 1.011 | 1.000 | 1.003 | 1.000 | 0.998 | 1.000 | 0.996 | 1.000 |
| | 272.13 | 10373 | 0.997 | 1.000 | 1.009 | 1.000 | 1.027 | 1.000 | 1.016 | 1.000 |
| mine-166-5 | 58.40 | 1569 | 0.997 | 1.000 | 0.998 | 1.000 | 0.997 | 1.000 | 0.999 | 1.000 |
| | 57.21 | 1191 | 1.007 | 1.000 | 1.003 | 1.000 | 0.998 | 1.000 | 1.005 | 1.000 |
| | 46.49 | 2244 | *1.064* | 1.000 | 1.010 | 1.000 | 1.007 | 1.000 | 0.999 | 1.000 |
| | 65.30 | 139 | 0.991 | 1.000 | 1.003 | 1.000 | 1.043 | 1.000 | 1.003 | 1.000 |
| | 42.28 | 2863 | *1.061* | 1.000 | 1.002 | 1.000 | 1.021 | 1.000 | 1.019 | 1.000 |
| mine-90-10 | 276.51 | 47643 | 1.001 | 1.000 | 0.997 | 1.000 | 0.999 | 1.000 | 0.997 | 1.000 |
| | 125.11 | 22195 | 1.002 | 1.000 | 1.001 | 1.000 | 1.013 | 1.000 | 1.018 | 1.000 |
| | 173.38 | 25979 | 1.042 | 1.000 | 1.001 | 1.000 | 1.010 | 1.000 | 1.020 | 1.000 |
| | 502.57 | 87496 | 0.998 | 1.000 | 0.995 | 1.000 | 1.005 | 1.000 | 0.997 | 1.000 |
| | 237.48 | 32520 | 1.002 | 1.000 | *1.054* | 1.000 | 0.992 | 1.000 | 0.999 | 1.000 |

cont. on next page

Table 6: Detailed computational results on `MMM-IP` with five random seeds. Relative changes by at least 5% are highlighted in bold and blue (**improvement**) or italic and red (*deterioration*).

| | degeneracy | | variablebounds | | conflicts | | infervals | | primsols | |
|---|---|---|---|---|---|---|---|---|---|---|
| Instance | **time** | **nodes** | **time$_Q$** | **nodes$_Q$** | **time$_Q$** | **nodes$_Q$** | **time$_Q$** | **nodes$_Q$** | **time$_Q$** | **nodes$_Q$** |
| misc03 | 1.40 | 398 | 0.988 | 1.000 | 0.992 | 1.000 | 0.983 | 1.000 | 0.988 | 1.000 |
| | 0.71 | 15 | 1.023 | 1.000 | 1.018 | 1.000 | 1.006 | 1.000 | 1.000 | 1.000 |
| | 0.73 | 19 | 1.000 | 1.000 | 0.983 | 1.000 | 0.988 | 1.000 | 0.977 | 1.000 |
| | 0.98 | 111 | 1.025 | 1.000 | 1.015 | 1.000 | 0.990 | 1.000 | 1.020 | 1.000 |
| | 0.64 | 15 | 1.018 | 1.000 | 1.000 | 1.000 | 1.000 | 1.000 | 1.018 | 1.000 |
| misc07 | 38.28 | 28289 | 0.997 | 1.000 | 0.990 | 1.000 | 0.982 | 1.000 | 0.995 | 1.000 |
| | 36.49 | 29782 | 1.007 | 1.000 | 1.006 | 1.000 | 1.006 | 1.000 | 1.007 | 1.000 |
| | 21.45 | 17632 | 1.003 | 1.000 | 1.010 | 1.000 | 1.004 | 1.000 | 1.003 | 1.000 |
| | 45.80 | 38975 | 1.014 | 1.000 | 1.013 | 1.000 | 1.004 | 1.000 | 1.007 | 1.000 |
| | 30.56 | 29476 | 1.010 | 1.000 | 1.014 | 1.000 | 1.009 | 1.000 | 1.009 | 1.000 |
| mitre | 14.57 | 1 | 1.006 | 1.000 | **0.757** | 1.000 | 0.984 | 1.000 | 0.988 | 1.000 |
| | 15.57 | 1 | 0.992 | 1.000 | **0.717** | 1.000 | 0.993 | 1.000 | 0.999 | 1.000 |
| | 15.52 | 1 | 1.005 | 1.000 | **0.718** | 1.000 | 1.005 | 1.000 | 1.001 | 1.000 |
| | 15.17 | 1 | 0.988 | 1.000 | **0.746** | 1.000 | 0.994 | 1.000 | 1.001 | 1.000 |
| | 13.39 | 1 | 1.037 | 1.000 | 0.994 | 1.000 | 1.027 | 1.000 | 1.028 | 1.000 |
| mod008 | 0.50 | 1 | 1.000 | 1.000 | 1.000 | 1.000 | 1.000 | 1.000 | 1.000 | 1.000 |
| | 0.53 | 2 | 1.000 | 1.000 | 0.987 | 1.000 | 1.013 | 1.000 | 1.007 | 1.000 |
| | 0.54 | 2 | 0.987 | 1.000 | 0.981 | 1.000 | 0.987 | 1.000 | 0.994 | 1.000 |
| | 0.51 | 1 | 1.000 | 1.000 | 1.000 | 1.000 | 1.007 | 1.000 | 1.013 | 1.000 |
| | 0.50 | 1 | 1.000 | 1.000 | 1.000 | 1.000 | 1.000 | 1.000 | 1.000 | 1.000 |
| mod010 | 0.50 | 2 | 1.000 | 1.000 | 1.000 | 1.000 | 1.000 | 1.000 | 1.000 | 1.000 |
| | 0.50 | 1 | 1.000 | 1.000 | 1.000 | 1.000 | 1.000 | 1.000 | 1.000 | 1.000 |
| | 0.50 | 1 | 1.000 | 1.000 | 1.000 | 1.000 | 1.000 | 1.000 | 1.000 | 1.000 |
| | 0.53 | 2 | 0.993 | 1.000 | 0.980 | 1.000 | 1.000 | 1.000 | 0.987 | 1.000 |
| | 0.50 | 1 | 1.000 | 1.000 | 1.000 | 1.000 | 1.000 | 1.000 | 1.000 | 1.000 |

cont. on next page

Table 6: Detailed computational results on `MMM-IP` with five random seeds. Relative changes by at least 5% are highlighted in bold and blue (**improvement**) or italic and red (*deterioration*).

| | degeneracy | | variablebounds | | conflicts | | infervals | | primsols | |
|---|---|---|---|---|---|---|---|---|---|---|
| Instance | **time** | **nodes** | **time$_Q$** | **nodes$_Q$** | **time$_Q$** | **nodes$_Q$** | **time$_Q$** | **nodes$_Q$** | **time$_Q$** | **nodes$_Q$** |
| mzzv11 | 332.03 | 2022 | 0.995 | 1.000 | 0.995 | 1.000 | 0.998 | 1.000 | 0.996 | 1.000 |
| | 229.82 | 959 | 1.005 | 1.000 | 0.999 | 1.000 | 0.991 | 1.000 | 0.999 | 1.000 |
| | 311.22 | 1421 | 0.987 | 1.000 | 1.000 | 1.000 | 0.996 | 1.000 | 0.988 | 1.000 |
| | 322.67 | 1579 | 1.009 | 1.000 | 1.014 | 1.000 | 1.009 | 1.000 | 0.998 | 1.000 |
| | 186.23 | 627 | 1.002 | 1.000 | 0.999 | 1.000 | 1.000 | 1.000 | 0.997 | 1.000 |
| mzzv42z | 151.11 | 265 | 0.994 | 1.000 | 1.003 | 1.000 | 1.001 | 1.000 | 1.001 | 1.000 |
| | 131.14 | 125 | 0.994 | 1.000 | 0.995 | 1.000 | 0.998 | 1.000 | 0.995 | 1.000 |
| | 198.89 | 200 | 1.000 | 1.000 | 0.996 | 1.000 | 0.995 | 1.000 | 0.999 | 1.000 |
| | 214.35 | 456 | 0.997 | 1.000 | 1.005 | 1.000 | 0.998 | 1.000 | 0.990 | 1.000 |
| | 167.53 | 300 | 1.005 | 1.000 | 0.998 | 1.000 | 1.000 | 1.000 | 1.001 | 1.000 |
| neos-1109824 | 37.38 | 1796 | 0.998 | 1.000 | 0.997 | 1.000 | 0.996 | 1.000 | 0.997 | 1.000 |
| | 46.57 | 3438 | 1.012 | 1.000 | 1.006 | 1.000 | 1.008 | 1.000 | 1.035 | 1.000 |
| | 34.97 | 1913 | 1.013 | 1.000 | 1.009 | 1.000 | 1.021 | 1.000 | 1.022 | 1.000 |
| | 77.32 | 8290 | 1.043 | 1.000 | 1.005 | 1.000 | *1.055* | 1.000 | 1.036 | 1.000 |
| | 39.22 | 2247 | 1.008 | 1.000 | 1.010 | 1.000 | 1.010 | 1.000 | 1.038 | 1.000 |
| neos-1337307 | 3600.00 | 306352 | 1.000 | 0.991 | 1.000 | 0.993 | 1.000 | 0.994 | 1.000 | 0.993 |
| | 3600.00 | 269284 | 1.000 | 0.998 | 1.000 | 0.979 | 1.000 | 1.002 | 1.000 | 0.952 |
| | 3600.00 | 349832 | 1.000 | 0.970 | 1.000 | 0.996 | 1.000 | 0.996 | 1.000 | 0.991 |
| | 3600.00 | 301363 | 1.000 | 0.989 | 1.000 | 0.978 | 1.000 | **0.932** | 1.000 | 1.006 |
| | 3600.00 | 373217 | 1.000 | 0.995 | 1.000 | 0.991 | 1.000 | 0.983 | 1.000 | 0.986 |
| neos-1601936 | 2033.50 | 1453 | 0.999 | 1.000 | 0.994 | 1.000 | 0.999 | 1.000 | 0.996 | 1.000 |
| | 3600.00 | 2861 | 1.000 | 0.996 | 1.000 | 1.001 | 1.000 | 1.000 | 1.000 | 1.000 |
| | 1338.93 | 2661 | 1.013 | 1.000 | 1.006 | 1.000 | 1.008 | 1.000 | 1.009 | 1.000 |
| | 3600.00 | 6917 | 1.000 | 1.004 | 1.000 | 1.000 | 1.000 | 0.990 | 1.000 | 0.999 |
| | 3600.00 | 6139 | 1.000 | 1.000 | 1.000 | 0.977 | 1.000 | 1.000 | 1.000 | 1.010 |

cont. on next page

Table 6: Detailed computational results on `MMM-IP` with five random seeds. Relative changes by at least 5% are highlighted in bold and blue (**improvement**) or italic and red (*deterioration*).

| | degeneracy | | variablebounds | | conflicts | | infervals | | primsols | |
|---|---|---|---|---|---|---|---|---|---|---|
| Instance | **time** | **nodes** | **time$_Q$** | **nodes$_Q$** | **time$_Q$** | **nodes$_Q$** | **time$_Q$** | **nodes$_Q$** | **time$_Q$** | **nodes$_Q$** |
| neos-686190 | 83.27 | 4906 | 0.996 | 1.000 | 1.002 | 1.000 | 0.995 | 1.000 | 0.996 | 1.000 |
| | 86.55 | 6047 | 1.000 | 1.000 | 1.013 | 1.000 | 0.997 | 1.000 | 1.000 | 1.000 |
| | 115.34 | 9156 | 0.995 | 1.000 | 0.994 | 1.000 | 1.021 | 1.000 | 1.006 | 1.000 |
| | 117.19 | 9998 | 0.995 | 1.000 | 1.006 | 1.000 | 0.992 | 1.000 | 0.996 | 1.000 |
| | 120.84 | 9806 | 1.008 | 1.000 | 0.998 | 1.000 | 1.003 | 1.000 | 1.007 | 1.000 |
| neos-849702 | 207.98 | 1169 | *1.979* | *8.256* | *6.088* | *54.068* | *2.833* | *25.153* | *4.876* | *14.426* |
| | 1025.85 | 6328 | **0.139** | **0.049** | **0.772** | *1.387* | **0.158** | **0.594** | 0.999 | 1.000 |
| | 1001.35 | 10499 | **0.904** | **0.400** | **0.181** | **0.099** | **0.686** | *5.053* | 1.016 | 1.000 |
| | 267.07 | 683 | *2.864* | *2.388* | *2.815* | *9.508* | **0.763** | *13.476* | *5.812* | *84.398* |
| | 285.44 | 1596 | *1.679* | *1.538* | *1.067* | **0.353** | **0.470** | *4.297* | *5.308* | *35.870* |
| neos-934278 | 3600.00 | 1175 | 1.000 | 0.999 | 1.000 | 1.000 | 1.000 | 1.000 | 1.000 | 1.006 |
| | 3600.00 | 267 | 1.000 | 0.995 | 1.000 | 1.008 | 1.000 | 1.000 | 1.000 | 0.995 |
| | 3600.00 | 263 | 1.000 | 1.000 | 1.000 | 1.008 | 1.000 | 0.983 | 1.000 | 1.008 |
| | 3600.00 | 401 | 1.000 | 0.992 | 1.000 | 0.994 | 1.000 | 0.992 | 1.000 | 0.994 |
| | 3600.00 | 522 | 1.000 | 0.998 | 1.000 | 0.998 | 1.000 | 1.000 | 1.000 | 0.997 |
| neos18 | 8.13 | 79 | *1.234* | *1.078* | *2.113* | *3.598* | *1.604* | *2.369* | *2.852* | *4.073* |
| | 24.87 | 828 | 1.043 | *1.054* | *1.347* | *1.815* | *1.232* | *1.531* | 0.999 | 1.028 |
| | 28.94 | 1095 | *1.089* | 1.000 | **0.883** | 0.965 | **0.680** | **0.362** | **0.689** | **0.585** |
| | 15.42 | 386 | 1.001 | 0.965 | *1.107* | *1.091* | *1.161* | **0.872** | *2.091* | *2.175* |
| | 21.77 | 963 | 0.961 | **0.949** | *1.249* | *1.147* | **0.928** | **0.710** | *1.281* | **0.871** |
| ns1208400 | 938.01 | 3025 | **0.912** | **0.753** | **0.543** | **0.314** | **0.763** | **0.692** | 1.017 | 1.032 |
| | 2737.49 | 11477 | **0.667** | **0.375** | **0.269** | **0.100** | **0.220** | **0.081** | 1.005 | 1.000 |
| | 532.73 | 1719 | *1.828* | *1.341* | **0.710** | **0.800** | *1.579* | **0.825** | *2.439* | *2.687* |
| | 253.07 | 599 | *2.316* | *2.724* | *1.530* | *1.409* | *1.318* | **0.720** | *2.093* | *2.482* |
| | 1316.58 | 2660 | *1.221* | *1.112* | *1.259* | *2.671* | **0.330** | **0.303** | 1.015 | 1.000 |

cont. on next page

Table 6: Detailed computational results on MMM-IP with five random seeds. Relative changes by at least 5% are highlighted in bold and blue (**improvement**) or italic and red (*deterioration*).

| | degeneracy | | variablebounds | | conflicts | | infervals | | primsols | |
|---|---|---|---|---|---|---|---|---|---|---|
| Instance | **time** | **nodes** | **time**$_Q$ | **nodes**$_Q$ | **time**$_Q$ | **nodes**$_Q$ | **time**$_Q$ | **nodes**$_Q$ | **time**$_Q$ | **nodes**$_Q$ |
| ns1688347 | 51.67 | 437 | *1.239* | *1.255* | 1.025 | *1.220* | *1.103* | **0.719** | 0.998 | 1.000 |
| | 42.42 | 203 | *1.423* | *2.165* | *1.332* | *2.370* | *1.406* | *2.670* | *1.886* | *4.353* |
| | 161.82 | 3095 | *1.050* | **0.752** | **0.833** | **0.552** | *1.111* | **0.556** | 1.003 | 1.000 |
| | 133.50 | 1251 | **0.878** | **0.597** | 0.964 | **0.638** | 1.049 | *1.245* | 1.008 | 1.000 |
| | 51.48 | 399 | **0.948** | **0.890** | *1.378* | *2.481* | 0.953 | *1.156* | *1.274* | *1.691* |
| ns1830653 | 133.08 | 6918 | 0.993 | 1.000 | 0.998 | 1.000 | 0.994 | 1.000 | 0.995 | 1.000 |
| | 123.87 | 5226 | 0.999 | 1.000 | 0.997 | 1.000 | *1.061* | 1.000 | 1.007 | 1.000 |
| | 161.00 | 9597 | 1.004 | 1.000 | 1.001 | 1.000 | 1.001 | 1.000 | 1.007 | 1.000 |
| | 143.21 | 7435 | *1.054* | 1.000 | 1.005 | 1.000 | 0.995 | 1.000 | 1.002 | 1.000 |
| | 146.08 | 9194 | 1.000 | 1.000 | 1.004 | 1.000 | 1.007 | 1.000 | 1.002 | 1.000 |
| nsrand-ipx | 419.85 | 50124 | 1.006 | 1.000 | 1.005 | 1.000 | 1.007 | 1.000 | 1.001 | 1.000 |
| | 638.54 | 86869 | 1.002 | 1.000 | 1.005 | 1.000 | 0.999 | 1.000 | 1.005 | 1.000 |
| | 569.66 | 56575 | 0.994 | 1.000 | 1.003 | 1.000 | 1.005 | 1.000 | 1.004 | 1.000 |
| | 582.34 | 61056 | 1.004 | 1.000 | 1.001 | 1.000 | 1.000 | 1.000 | 1.002 | 1.000 |
| | 363.56 | 38012 | 1.007 | 1.000 | 1.002 | 1.000 | 1.002 | 1.000 | 1.000 | 1.000 |
| opt1217 | 0.50 | 1 | 1.000 | 1.000 | 1.000 | 1.000 | 1.000 | 1.000 | 1.000 | 1.000 |
| | 0.50 | 1 | 1.000 | 1.000 | 1.000 | 1.000 | 1.000 | 1.000 | 1.000 | 1.000 |
| | 0.50 | 1 | 1.000 | 1.000 | 1.000 | 1.000 | 1.000 | 1.000 | 1.000 | 1.000 |
| | 0.50 | 1 | 1.000 | 1.000 | 1.000 | 1.000 | 1.000 | 1.000 | 1.000 | 1.000 |
| | 0.50 | 1 | 1.000 | 1.000 | 1.000 | 1.000 | 1.000 | 1.000 | 1.000 | 1.000 |
| p0033 | 0.50 | 1 | 1.000 | 1.000 | 1.000 | 1.000 | 1.000 | 1.000 | 1.000 | 1.000 |
| | 0.50 | 1 | 1.000 | 1.000 | 1.000 | 1.000 | 1.000 | 1.000 | 1.000 | 1.000 |
| | 0.50 | 1 | 1.000 | 1.000 | 1.000 | 1.000 | 1.000 | 1.000 | 1.000 | 1.000 |
| | 0.50 | 1 | 1.000 | 1.000 | 1.000 | 1.000 | 1.000 | 1.000 | 1.000 | 1.000 |
| | 0.50 | 1 | 1.000 | 1.000 | 1.000 | 1.000 | 1.000 | 1.000 | 1.000 | 1.000 |

cont. on next page

Table 6: Detailed computational results on MMM-IP with five random seeds. Relative changes by at least 5% are highlighted in bold and blue (**improvement**) or italic and red (*deterioration*).

| | degeneracy | | variablebounds | | conflicts | | infervals | | primsols | |
|---|---|---|---|---|---|---|---|---|---|---|
| Instance | **time** | **nodes** | **time$_Q$** | **nodes$_Q$** | **time$_Q$** | **nodes$_Q$** | **time$_Q$** | **nodes$_Q$** | **time$_Q$** | **nodes$_Q$** |
| p0201 | 0.54 | 5 | 1.006 | 1.000 | 1.013 | 1.000 | 1.019 | 1.000 | 1.000 | 1.000 |
| | 0.74 | 9 | 0.989 | 1.000 | 1.000 | 1.000 | 0.983 | 1.000 | 1.011 | 1.000 |
| | 1.64 | 5 | 1.000 | 1.000 | 1.004 | 1.000 | 0.985 | 1.000 | 0.996 | 1.000 |
| | 0.90 | 13 | 1.000 | 1.000 | 1.000 | 1.000 | 1.000 | 1.000 | 1.011 | 1.000 |
| | 0.89 | 59 | 1.005 | 1.000 | 1.005 | 1.000 | 1.011 | 1.000 | 1.016 | 1.000 |
| p0282 | 0.50 | 2 | 1.000 | 1.000 | 1.000 | 1.000 | 1.000 | 1.000 | 1.000 | 1.000 |
| | 0.50 | 2 | 1.000 | 1.000 | 1.000 | 1.000 | 1.000 | 1.000 | 1.000 | 1.000 |
| | 0.50 | 2 | 1.000 | 1.000 | 1.000 | 1.000 | 1.000 | 1.000 | 1.000 | 1.000 |
| | 0.50 | 3 | 1.000 | 1.000 | 1.000 | 1.000 | 1.007 | 1.000 | 1.000 | 1.000 |
| | 0.50 | 2 | 1.000 | 1.000 | 1.000 | 1.000 | 1.000 | 1.000 | 1.000 | 1.000 |
| p0548 | 0.50 | 1 | 1.000 | 1.000 | 1.000 | 1.000 | 1.000 | 1.000 | 1.000 | 1.000 |
| | 0.50 | 1 | 1.000 | 1.000 | 1.000 | 1.000 | 1.000 | 1.000 | 1.000 | 1.000 |
| | 0.50 | 1 | 1.000 | 1.000 | 1.000 | 1.000 | 1.000 | 1.000 | 1.000 | 1.000 |
| | 0.50 | 1 | 1.000 | 1.000 | 1.000 | 1.000 | 1.000 | 1.000 | 1.000 | 1.000 |
| | 0.50 | 1 | 1.000 | 1.000 | 1.000 | 1.000 | 1.000 | 1.000 | 1.000 | 1.000 |
| p2756 | 2.51 | 5 | 0.991 | 1.000 | 0.989 | 1.000 | 1.000 | 1.000 | 0.994 | 1.000 |
| | 3.33 | 5 | 1.009 | 1.000 | 1.005 | 1.000 | 1.009 | 1.000 | 1.007 | 1.000 |
| | 2.13 | 3 | 0.994 | 1.000 | 0.990 | 1.000 | 0.987 | 1.000 | 1.006 | 1.000 |
| | 2.54 | 3 | 1.000 | 1.000 | 1.008 | 1.000 | 1.006 | 1.000 | 1.014 | 1.000 |
| | 1.28 | 2 | 1.000 | 1.000 | 0.996 | 1.000 | 1.004 | 1.000 | 1.009 | 1.000 |
| protfold | 3600.00 | 3244 | 1.000 | **0.823** | 1.000 | **0.632** | 1.000 | *2.305* | 1.000 | *1.328* |
| | 3600.00 | 2504 | 1.000 | **0.715** | 1.000 | 0.988 | 1.000 | *2.045* | 1.000 | 0.950 |
| | 3600.00 | 2755 | 1.000 | **0.858** | 1.000 | *1.160* | 1.000 | *1.907* | 1.000 | **0.733** |
| | 3600.00 | 1831 | 1.000 | *1.278* | 1.000 | *1.551* | 1.000 | *2.081* | 1.000 | *1.256* |
| | 3600.00 | 1867 | 1.000 | **0.814** | 1.000 | *1.325* | 1.000 | *1.584* | 1.000 | **0.831** |

cont. on next page

Table 6: Detailed computational results on MMM-IP with five random seeds. Relative changes by at least 5% are highlighted in bold and blue (**improvement**) or italic and red (*deterioration*).

| | degeneracy | | variablebounds | | conflicts | | infervals | | primsols | |
|---|---|---|---|---|---|---|---|---|---|---|
| Instance | time | nodes | $\text{time}_Q$ | $\text{nodes}_Q$ | $\text{time}_Q$ | $\text{nodes}_Q$ | $\text{time}_Q$ | $\text{nodes}_Q$ | $\text{time}_Q$ | $\text{nodes}_Q$ |
| pw-myciel4 | 3600.00 | 361131 | 1.000 | *1.253* | **0.383** | **0.439** | 1.000 | *1.320* | 1.000 | *1.077* |
| | 3600.00 | 365865 | 1.000 | 1.039 | 1.000 | **0.895** | 1.000 | *1.253* | 1.000 | *1.054* |
| | 3600.00 | 323450 | 1.000 | *1.199* | 1.000 | *1.102* | 1.000 | *1.294* | 1.000 | *1.161* |
| | 3600.00 | 358151 | 1.000 | 1.020 | 1.000 | **0.856** | 1.000 | *1.395* | **0.925** | *1.097* |
| | 3054.26 | 378549 | *1.179* | **0.843** | *1.179* | *1.128* | *1.179* | 1.019 | *1.179* | **0.893** |
| qnet1 | 4.54 | 33 | 1.002 | 1.000 | 0.998 | 1.000 | 1.002 | 1.000 | 1.005 | 1.000 |
| | 3.12 | 3 | 0.995 | 1.000 | 1.007 | 1.000 | 1.002 | 1.000 | 0.993 | 1.000 |
| | 3.92 | 29 | 1.002 | 1.000 | 0.996 | 1.000 | 0.996 | 1.000 | 1.006 | 1.000 |
| | 2.26 | 3 | 0.991 | 1.000 | 0.991 | 1.000 | 0.994 | 1.000 | 0.991 | 1.000 |
| | 3.71 | 34 | 0.994 | 1.000 | 1.002 | 1.000 | 0.996 | 1.000 | 0.987 | 1.000 |
| qnet1_o | 3.87 | 27 | 1.012 | 1.000 | 1.002 | 1.000 | 1.010 | 1.000 | 1.008 | 1.000 |
| | 3.38 | 24 | 1.000 | 1.000 | 0.995 | 1.000 | 0.995 | 1.000 | 1.002 | 1.000 |
| | 3.39 | 15 | 0.995 | 1.000 | 1.000 | 1.000 | 0.993 | 1.000 | 1.000 | 1.000 |
| | 1.67 | 2 | 1.000 | 1.000 | 1.000 | 1.000 | 0.993 | 1.000 | 1.004 | 1.000 |
| | 1.84 | 3 | 1.000 | 1.000 | 1.007 | 1.000 | 1.000 | 1.000 | 1.007 | 1.000 |
| reblock67 | 223.13 | 49530 | 0.998 | 1.000 | 0.995 | 1.000 | 0.995 | 1.000 | 0.989 | 1.000 |
| | 161.91 | 40625 | 1.001 | 1.000 | 1.004 | 1.000 | 1.006 | 1.000 | 1.019 | 1.000 |
| | 194.68 | 49882 | 1.002 | 1.000 | 1.006 | 1.000 | 0.999 | 1.000 | 1.032 | 1.000 |
| | 200.17 | 50996 | 1.006 | 1.000 | 1.001 | 1.000 | 0.992 | 1.000 | 0.995 | 1.000 |
| | 239.14 | 57566 | 1.017 | 1.000 | 0.999 | 1.000 | 1.002 | 1.000 | 0.999 | 1.000 |
| rmine6 | 971.59 | 111047 | 0.997 | 1.000 | 1.002 | 1.000 | 0.999 | 1.000 | 1.003 | 1.000 |
| | 882.29 | 89213 | 0.999 | 1.000 | 1.040 | 1.000 | 1.013 | 1.000 | 1.026 | 1.000 |
| | 776.69 | 91530 | 0.992 | 1.000 | 1.000 | 1.000 | 1.004 | 1.000 | 1.021 | 1.000 |
| | 742.25 | 85429 | 1.010 | 1.000 | 1.001 | 1.000 | 1.005 | 1.000 | 1.000 | 1.000 |
| | 773.39 | 76017 | 0.995 | 1.000 | 1.002 | 1.000 | 1.013 | 1.000 | 0.995 | 1.000 |

cont. on next page

Table 6: Detailed computational results on MMM-IP with five random seeds. Relative changes by at least 5% are highlighted in bold and blue (**improvement**) or italic and red (*deterioration*).

| Instance | degeneracy | | variablebounds | | conflicts | | infervals | | primsols | |
|---|---|---|---|---|---|---|---|---|---|---|
| | **time** | **nodes** | **time$_Q$** | **nodes$_Q$** | **time$_Q$** | **nodes$_Q$** | **time$_Q$** | **nodes$_Q$** | **time$_Q$** | **nodes$_Q$** |
| rococoC10-001000 | 603.64 | 46983 | 0.995 | 1.000 | 0.997 | 1.000 | 0.994 | 1.000 | 0.995 | 1.000 |
| | 836.01 | 75156 | 0.999 | 1.000 | 0.998 | 1.000 | 0.998 | 1.000 | 1.007 | 1.000 |
| | 516.72 | 37178 | *1.052* | 1.000 | 1.006 | 1.000 | 1.003 | 1.000 | 1.006 | 1.000 |
| | 893.05 | 51410 | 0.999 | 1.000 | 1.006 | 1.000 | 1.004 | 1.000 | 1.001 | 1.000 |
| | 1048.40 | 66016 | 1.004 | 1.000 | 1.010 | 1.000 | 1.008 | 1.000 | 1.003 | 1.000 |
| seymour | 3600.00 | 102192 | 1.000 | 1.003 | 1.000 | 1.002 | 1.000 | 1.007 | 1.000 | 1.000 |
| | 3600.00 | 120186 | 1.000 | 1.002 | 1.000 | 0.996 | 1.000 | 1.000 | 1.000 | 1.000 |
| | 3600.00 | 123015 | 1.000 | 1.003 | 1.000 | 0.997 | 1.000 | 0.997 | 1.000 | 0.990 |
| | 3600.00 | 108085 | 1.000 | 1.001 | 1.000 | 0.996 | 1.000 | 0.998 | 1.000 | 0.995 |
| | 3600.00 | 97046 | 1.000 | 1.001 | 1.000 | 0.994 | 1.000 | 0.999 | 1.000 | 1.001 |
| sp98ir | 58.06 | 3948 | 0.999 | 1.000 | 0.997 | 1.000 | 1.000 | 1.000 | 0.999 | 1.000 |
| | 65.95 | 4178 | 1.039 | 1.000 | 1.003 | 1.000 | 1.012 | 1.000 | 1.006 | 1.000 |
| | 67.93 | 4735 | 1.000 | 1.000 | 1.029 | 1.000 | 1.020 | 1.000 | 1.014 | 1.000 |
| | 72.47 | 4603 | 1.003 | 1.000 | 1.005 | 1.000 | 1.003 | 1.000 | 1.011 | 1.000 |
| | 69.87 | 4780 | *1.063* | 1.000 | 0.998 | 1.000 | 1.001 | 1.000 | 1.013 | 1.000 |
| stein27 | 1.17 | 2749 | 0.972 | 1.000 | 1.014 | 1.000 | 1.023 | 1.000 | 0.995 | 1.000 |
| | 1.23 | 2699 | 1.000 | 1.000 | 0.969 | 1.000 | 0.996 | 1.000 | 0.978 | 1.000 |
| | 1.19 | 2889 | 0.995 | 1.000 | 1.018 | 1.000 | 1.014 | 1.000 | 1.023 | 1.000 |
| | 1.13 | 2873 | 0.953 | 1.000 | 0.958 | 1.000 | 0.962 | **0.947** | 1.014 | 1.000 |
| | 1.08 | 2827 | 1.019 | 1.000 | 1.010 | 1.000 | 1.014 | 1.029 | 0.990 | 1.000 |
| stein45 | 12.59 | 39175 | 0.994 | 1.000 | 0.989 | 1.000 | 0.996 | 1.039 | 0.997 | 1.000 |
| | 12.49 | 38805 | 1.004 | 1.000 | 1.005 | 1.000 | 1.024 | 1.005 | 1.003 | 1.000 |
| | 13.04 | 41279 | 1.005 | 1.000 | 1.024 | 1.000 | 0.998 | 1.016 | 1.004 | 1.000 |
| | 12.49 | 38104 | 1.006 | 1.000 | 1.004 | 1.000 | 0.980 | 0.985 | 1.004 | 1.000 |
| | 13.02 | 40755 | 1.008 | 1.000 | 0.997 | 1.000 | 1.009 | 0.994 | 0.994 | 1.000 |

cont. on next page

Table 6: Detailed computational results on `MMM-IP` with five random seeds. Relative changes by at least 5% are highlighted in bold and blue (**improvement**) or italic and red (*deterioration*).

| | degeneracy | | variablebounds | | conflicts | | infervals | | primsols | |
|---|---|---|---|---|---|---|---|---|---|---|
| Instance | **time** | **nodes** | **time$_{\mathbf{Q}}$** | **nodes$_{\mathbf{Q}}$** | **time$_{\mathbf{Q}}$** | **nodes$_{\mathbf{Q}}$** | **time$_{\mathbf{Q}}$** | **nodes$_{\mathbf{Q}}$** | **time$_{\mathbf{Q}}$** | **nodes$_{\mathbf{Q}}$** |
| tanglegram2 | 4.11 | 3 | 1.002 | 1.000 | 1.006 | 1.000 | 0.990 | 1.000 | 1.000 | 1.000 |
| | 3.90 | 3 | 1.022 | 1.000 | 0.994 | 1.000 | 1.014 | 1.000 | 1.016 | 1.000 |
| | 4.20 | 3 | 1.012 | 1.000 | 1.008 | 1.000 | 1.015 | 1.000 | 1.008 | 1.000 |
| | 4.12 | 3 | 1.020 | 1.000 | 1.008 | 1.000 | 1.008 | 1.000 | 1.016 | 1.000 |
| | 4.20 | 3 | 1.000 | 1.000 | 1.000 | 1.000 | 1.038 | 1.000 | 1.000 | 1.000 |